\documentclass[12pt]{article}
\usepackage{amsmath,amsthm,amssymb,color,epic,eepic, amsfonts,graphicx,latexsym,graphics, 
cite
}
\usepackage{mathrsfs}
\usepackage{wrapfig}

\renewcommand{\baselinestretch}{1.2}

\def\baselinestretch{1.4}
\newlength{\minitwocolumn}
\newcommand{\Z}{{\Bbb Z}} 
\newcommand{\R}{{\Bbb R}} 
\newcommand{\C}{{\Bbb C}} 
\newcommand{\N}{{\Bbb N}} 
\newcommand{\FF}{{\Bbb F}} 
\newcommand{\PP}{{\Bbb P}} 
\newcommand{\F}{{\mathcal F}}

\newcommand{\gS}{\mathfrak S}

\newcommand{\cC}{{\mathcal C}}
\newcommand{\cD}{{\mathcal D}}
\newcommand{\cA}{{\mathcal A}}
\newcommand{\cU}{{\mathcal U}}

\newcommand{\cB}{{\mathcal B}}

\newcommand{\cN}{{\mathcal N}}
\newcommand{\cL}{{\mathcal L}}

\newcommand{\cR}{{\mathcal R}}
\newcommand{\cP}{{\mathcal P}}

\newcommand{\cM}{{\mathcal M}}

\newcommand{\cQ}{{\mathcal Q}}

\newcommand{\cV}{{\mathcal V}}
\newcommand{\cW}{{\mathcal W}}

\newcommand{\la}{\lambda}
\newcommand{\La}{\Lambda}

\newcommand{\al}{\alpha}

\newcommand{\ep}{\epsilon}
\newcommand{\vep}{\varepsilon}

\newcommand{\tT}{\widetilde{T}}

\newcommand{\trho}{{\widetilde{\rho}}}

\newcommand{\bep}{\bar{\epsilon}}

\newcommand{\hf}{\widehat{f}}

\newcommand{\hV}{\widehat{V}}

\newcommand{\bfh}{{\bf h}}

\newcommand{\bfd}{{\bf d}}
\newcommand{\bfv}{{\bf v}}
\newcommand{\bfw}{{\bf w}}
\newcommand{\bfx}{{\bf x}}

\newcommand{\bfu}{{\bf u}}
\newcommand{\bu}{{\bf u}}
\newcommand{\bft}{{\bf t}}
\newcommand{\bfz}{{\bf z}}
\newcommand{\bfM}{{\bf M}}
\newcommand{\bfL}{{\bf L}}

\newcommand{\bal}{{\boldsymbol \al}}
\newcommand{\bbe}{{\boldsymbol \beta}}
\newcommand{\bga}{{\boldsymbol \gamma}}
\newcommand{\bla}{{\boldsymbol \la}}
\newcommand{\bmu}{{\boldsymbol \mu}}
\newcommand{\bom}{{\boldsymbol \omega}}
\newcommand{\blla}{{\boldsymbol \la}}
\newcommand{\bxi}{{\boldsymbol \xi}}
\newcommand{\blmu}{{\boldsymbol \mu}}

\newcommand{\blempty}{{\boldsymbol \emptyset}}
\newcommand{\bgz}{{\boldsymbol{\gz}}}

\newcommand{\blLa}{{\boldsymbol \Lambda}}

\newcommand{\Hom}{\operatorname{Hom}}

\newcommand{\nn}{{\nonumber}}
\newcommand{\bea}{\begin{eqnarray}}
\newcommand{\ena}{\end{eqnarray}}
\newcommand{\be}{\begin{eqnarray*}}
\newcommand{\en}{\end{eqnarray*}}

\newcommand{\beit}{\begin{itemize}}
\newcommand{\enit}{\end{itemize}}

\newcommand{\lb}[1]{\label{#1}}

\newcommand{\ds}[1]{{\displaystyle #1 }}

\newcommand{\End}{\mathrm{ End}}

\newcommand{\id}{\mathrm{ id}}

\newcommand{\wt}{\mathrm{ wt}}
\newcommand{\Span}{\mathrm{Span}}
\newcommand{\Sym}{\mathrm{Sym}}

\newcommand{\E}{\mathrm{E}}
\newcommand{\rA}{\mathrm{A}}
\newcommand{\rT}{\mathrm{T}}
\newcommand{\rK}{\mathrm{K}}
\newcommand{\rH}{\mathrm{H}}

\newcommand{\Stab}{\mathrm{Stab}}
\newcommand{\gC}{\mathfrak{C}}

\newcommand{\Lie}{\mathrm{Lie}\ }

\newcommand{\gz}{\mathfrak{z}}

\newcommand{\bra}[1]{\langle #1 |}        
\newcommand{\ket}[1]{{| #1 \rangle}}      

\def\infq4p#1{{(#1;q^4,p)_\infty}}

\newcommand{\tot}{\, \widetilde{\otimes}\, }

\newcommand{\mmatrix}[1]{\begin{matrix} #1 \end{matrix}}

\font\teneufm=eufm10
\font\seveneufm=eufm7
\font\fiveeufm=eufm5
\newfam\eufmfam
\textfont\eufmfam=\teneufm
\scriptfont\eufmfam=\seveneufm
\scriptscriptfont\eufmfam=\fiveeufm
\def\frak#1{{\fam\eufmfam\relax#1}}
\let\goth\mathfrak
\newcommand{\slth}{\widehat{\goth{sl}}_2}

\newcommand{\slnh}{\widehat{\goth{sl}}_N}
\newcommand{\sln}{\goth{sl}_N}
\newcommand{\g}{\goth{g}}

\newcommand{\Bqla}{{{\mathcal B}_{q,\lambda}}}

\newcommand{\glntbig}{{\mbox{\fourteeneufm gl}}_{N,tor}}

\newcommand{\gl}{{\goth{gl}}}

\newcommand{\glnt}{{\goth{gl}_{N,tor}}}

\newcommand{\glnh}{\widehat{\goth{gl}}_N}

\newcommand{\h}{\goth{h}}
\newcommand{\gh}{\widehat{\goth{g}}}

\font\fourteeneufm=eufm10 scaled\magstep2    
   
\makeatletter
\@addtoreset{equation}{section}
\makeatother

\newtheorem{thm}{Theorem}[section]
\newtheorem{prop}[thm]{Proposition}
\newtheorem{lem}[thm]{Lemma}
\newtheorem{cor}[thm]{Corollary}
\newtheorem{conj}[thm]{Conjecture}

\newtheorem{dfn}[thm]{Definition}

\begin{document}

%
%
\begin{center}
{\Large Elliptic Quantum Toroidal Algebra $U_{t_1,t_2,p}(\glntbig)$ and \\  
Elliptic Stable Envelopes for the $A^{(1)}_{N-1}$ 
Quiver Varieties}
\\[10mm]
{\large  Hitoshi Konno and Andrey Smirnov}
\end{center}

\bigskip

\begin{abstract}
\noindent 
We propose a new construction of vertex operators of the elliptic quantum toroidal algebra  $U_{t_1,t_2,p}(\gl_{N,tor})$ by combining representations of the algebra and formulas of the elliptic stable envelopes for the  $A^{(1)}_{N-1}$ quiver variety $\cM(\bfv,\bfw)$. Compositions of the vertex operators turn out  consistent to the shuffle product formula of the elliptic stable envelopes.  Their highest to highest expectation values provide  $\rK$-theoretic vertex functions for $\cM(\bfv,\bfw)$.  
We also derive exchange relation of the vertex operators 
 and construct a $L$-operator satisfying the $RLL=LLR^*$ relation with $R$ and $R^*$ being elliptic dynamical $R$-matrices defined as transition matrices of the elliptic stable envelopes.  
Assuming a universal form of $L$ and defining a comultiplication $\Delta$ in terms of it, we show that our vertex operators are  intertwining operators of the $U_{t_1,t_2,p}(\gl_{N,tor})$-modules w.r.t $\Delta$.   
\end{abstract}
\nopagebreak

\section{Introduction}\lb{sec:intro}
This is a continuation of our work \cite{KS} on a construction of vertex operators 
of the elliptic quantum toroidal algebras  and making a connection between representations of the algebras and geometry of affine quiver varieties\cite{Na94,Na98,Na01}, especially the elliptic stable envelopes\cite{AO,Ok20}. In \cite{KS}, we studied  the toroidal $\gl_1$ case and made a connection to the equivariant $\rK$-theory and the elliptic cohomology of the Jordan quiver varieties $\cM(n,r)$ known as the ADHM instanton moduli spaces. In this paper we discuss its extension to the higher rank toroidal $\gl_N$ case. The corresponding quiver varieties are the affine quiver varieties of $A^{(1)}_{N-1}$ type. 

Let $\g$ be a finite dimensional simple Lie algebra over $\C$. Let $\gh$  be the untwisted affine Lie algebra\cite{Kac} i.e. a one-dimensional canonical central extension of the loop Lie algebra $\C^{\times}\to \g$. Let  $\g_{tor}$ be the toroidal Lie algebra associated with $\g$ i.e.     
 a canonical two-dimensional central extension of the double loop Lie algebra $\C^{\times}\times \C^{\times} \to \g$\cite{GKV}. 
The elliptic quantum groups (EQG)  $U_{q,p}(\gh)$ associated with  $\gh$ were introduced in \cite{K98,JKOS,FKO,Konno09}.  Their various representations were constructed in 
\cite{K98,JKOS,FKO, KK03, Konno08,Konno09,Konno16,Konno17,Konno18}. See also \cite{KonnoBook,Konno24}. 
It's toroidal version,  i.e.  the elliptic quantum toroidal algebras (EQTA)  $U_{t_1,t_2,p}(\g_{tor})$ associated with $\g_{tor}$, were formulated  in \cite{KOgl1,KS,KO24} in the same scheme as $U_{q,p}(\gh)$. Here the parameters $t_1, t_2$ are independent only for $\g_{tor}=\gl_{N,tor}$. For the other cases, $t_1=t_2$. Their various representations were constructed  in \cite{KOgl1,KS,KO24}. The $U_{q,p}(\gh)$ and  $U_{t_1,t_2,p}(\g_{tor})$  
are natural $p$-deformations of the quantum affine algebras $U_q(\gh)$\cite{Drinfeld85,Jimbo85} and the quantum toroidal algebras $U_{t_1,t_2}(\g_{tor})$\cite{GKV} in the Drinfeld realization\cite{Drinfeld}, respectively. 

The notion of the elliptic stable envelopes were introduced in \cite{AO} as certain good elliptic cohomology classes for the torus equivariant elliptic cohomology $\E_\rT(X)$ of the quiver variety $X$. See also \cite{Ok20}. For the cohomology $\rH^*_\rT(X)$ and  the $\rK$-theory $\rK_\rT(X)$ cases, the stable envelopes are identified with the Chern-Schwartz-MacPherson classes and the motivic Chern classes, respectively\cite{FR,FRW,RV,AMSS}. A key property to the geometric representation in \cite{AO} is  that the transition matrices of the elliptic stable envelopes defined with different chambers give elliptic solutions to the dynamical Yang-Baxter equation\cite{Felder,JKOStg}. There the K\"ahler parameters in $\E_\rT(X)$ are identified with the dynamical parameters. 
By making use of such solutions, one may obtain 
  representations of EQG/EQTA associated with $X$ geometrically.  See also \cite{MO} and \cite{Ok}, for the rational case, i.e. the Yangian double, and the trigonometric case, i.e. the quantum group, respectively. To compare such geometric representations to those of 
  $U_{q,p}(\gh)$ and $U_{t_1,t_2,p}(\g_{tor})$ is  our basic  aim.

In the previous works \cite{Konno17, Konno18, KonnoBook, Konno24, KO24, KOgl1, KS}, we have studied several connections between representations of $U_{q,p}(\slnh)$, $U_{t_1,t_1,p}(\gl_{1,tor})$ and geometry of corresponding quiver varieties. We have shown that the vertex operators of them, which are intertwining operators w.r.t. the standard comultiplication, are realized as operator valued integrals, called the screened vertex operators, with the elliptic stable envelopes being their integration kernels. 
Compositions of the vertex operators turns out consistent to the shuffle product formula for the elliptic stable envelopes obtained  geometrically\cite{Botta}. Furthermore the highest to highest expectation values of compositions of the vertex operators give  $\rK$-theoretic vertex functions\cite{AOMosc,AFO,Ok,Ok20,OS} for the corresponding quiver varieties. Such nice consistency between representations of EQG/EQTA and geometry of quiver varieties leads us to  
 the following conjectures.  
\begin{conj} 
\

\vspace{-2mm}
\begin{itemize}
\item[(1)] Vertex operators of  $U_{q,p}(\gh)$ or $U_{t_1,t_2,p}(\g_{tor})$ defined as intertwining 
operators  w.r.t. the standard comultiplication are realized as screened vertex operators with the elliptic stable envelope for $\E_\rT(X)$ being their integration kernels. The corresponding quiver variety $X$ is  the Dynkin quiver variety of type $\g$ for $U_{q,p}(\gh)$ and the affine Dynkin quiver variety of type $\gh$ for $U_{t_1,t_1,p}(\g_{tor})$,respectively.  
In particular, the components of the vertex operators are labeled by the $\rT$-fixed points of $X$. 
\item[(2)] Composition of thus obtained vertex operators gives the  shuffle product  formula for the elliptic stable envelopes as a relation among the integration kernels. 
\item[(3)] The highest to highest expectation values of compositions of vertex operators give the $\rK$-theoretic vertex functions counting quasimaps $\PP^1\dashrightarrow X$.
\end{itemize}
\end{conj}

Examples: 
\begin{itemize}
\item[1)] For $U_{q,p}(\slnh)$,   $X$= the $A_{N-1}$ linear quiver varieties  with dimension vector $\bfv=(v_1, \cdots, v_{N-1})$, $0\leq v_1\leq \cdots \leq v_{N-1}\leq n$ and framing vector $\bfw=(0,\cdots,0,n)$\cite{Konno17,Konno18}.  $X$ is isomorphic to $T^*\hspace{-0.8mm}fl(\bfv)$,  the cotangent bundle to the partial flag variety $fl(\bfv)$ consists of the flags $0\subset \C^{v_1}\subset \cdots\subset \C^{v_{N-1}}\subset \C^n$\cite{Na94}. 
Let $V_{u}=V[u,u^{-1}]$, $V=\oplus_{\mu=1}^N\C v_\mu$ be the evaluation representation associated with the vector representation and $\F_{\omega}$ be the level-1 highest weight representation of  $U_{q,p}(\slnh)$ with highest weight $\omega$. The basic type I vertex operator is the intertwining operator $\Phi(u): \F_{\omega} \to V_{u}\tot \F_{\omega'}$ with 
$\Phi(u)=\sum_{\mu}v_\mu\tot \Phi_\mu(u)$. The composition 
\be
&&\Phi(u_1,\cdots,u_n)=(\id\tot\cdots\tot\id\tot\Phi(u_1))\circ  \cdots\circ (\id\tot\Phi(u_{n-1}))\circ\Phi(u_{n})\\
&&\hspace{5cm} :\ \F_{\omega}\ \to   \ V_{u_{n}}\tot \cdots\tot V_{u_1}\tot\F_{\omega'}
\en
has components  $\Phi_{\mu_1}(u_1)\cdots \Phi_{\mu_n}(u_n)$ w.r.t. the standard bases $v_{\mu_n}\tot\cdots\tot v_{\mu_1}$. There is a 
bijection between $(\mu_1,\cdots,\mu_n)\in [1,N]^n$ and a partition $(I_1,\cdots,I_N)$  
of $[1,n]$, i.e. $I_1\cup \cdots \cup I_N=[1,n]$, $I_i\cap I_j=\emptyset \ (i\not=j)$. Setting $v_i=|I_1\cup \cdots \cup I_i |$, the partitions  $(I_1,\cdots,I_N)$ label the torus fixed points of $fl(\bfv)$\cite{Konno17}. By a bijection between the Gelfand-Tsetlin bases $\xi_{(I_1,\cdots,I_N)}$ of $V_{u_{n}}\tot \cdots\tot V_{u_1}$ and the set of fixed point classes $[(I_1,\cdots,I_N)]$ of $\E_\rT(T^*\hspace{-0.1cm}fl(\bfv))$\cite{Konno18} the level-0 representation $V_{u_{n}}\tot \cdots\tot V_{u_1}$ is equivalent to the  geometric action of $U_{q,p}(\slnh)$ on the fixed point classes of $\bigoplus_{\bfv}\E_\rT(T^*\hspace{-0.1cm}fl(\bfv))$\cite{Konno18,KonnoBook,Konno24}.  The  vertex operator $\Phi_{\mu_1}(u_1)\cdots \Phi_{\mu_n}(u_n)$ is realized as a screened vertex operator with the elliptic stable envelope of the fixed point $(I_1,\cdots,I_N)$ for $\E_\rT(T^*\hspace{-0.1cm}fl(\bfv))$ being its integration kernel.  
In particular, the basic one $\Phi_\mu(u)$ is the screened vertex operator with the 
elliptic stable envelope of the fixed point $(\emptyset,\cdots,I_\mu=\{1\},\cdots,\emptyset)$ for $\E_\rT(T^*\hspace{-0.1cm}fl(\bfv))$ with $\ds{\bfv=(\underbrace{0,\cdots,0}_{\mu-1},\underbrace{1,\cdots,1}_{N-\mu})}$. 
 For the shuffle product see \cite{Konno17} and Remark after Corollary \ref{cor:shuffletStabgen} below.  For the $\rK$-theoretic vertex function, see  \cite{Konno24}. 
\item[2)] For $U_{t_1,t_2,p}(\gl_{1,tor})$,  $X$= the Jordan quiver varieties $\cM(n,r)$\cite{KS}. Let $\F^{(0,-1)}_u=\Span\{\ket{\la}_u, \la:\mbox{partition} \}$, $u\in \C^\times$ be the level-(0,-1) representation  and $\F^{(1,N)}_v$, $N\in \Z$, $v\in \C^\times$  be the level-$(1,N)$  representation of $U_{t_1,t_2,p}(\gl_{1,tor})$. The basic type I vertex operator is the intertwining operator $\Phi(u): \F^{(1,N)}_{v} \to \F^{(0,-1)}_{u}\tot \F^{(1,N+1)}_{-v/u}$ with 
$\Phi(u)=\sum_{\la}\ket{\la}_u\tot \Phi_\la(u)$. The composition 
\be
&&\Phi(u_1,\cdots,u_n)=(\id\tot\cdots\tot\id\tot\Phi(u_1))\circ  \cdots\circ (\id\tot\Phi(u_{r-1}))\circ\Phi(u_{r})\\
&&\hspace{5cm} :\ \F^{(1,N)}_{v}\ \to   \ \F^{(0,-1)}_{u_r}\tot \cdots\tot \F^{(0,-1)}_{u_1}\tot\F^{(1,N+1)}_{(-)^rv/u_1\cdots u_r}
\en
has components $\Phi_{\la_1}(u_1)\cdots \Phi_{\la_r}(u_r)$ w.r.t. the bases $\ket{\la_r}_{u_r}\tot\cdots\tot \ket{\la_1}_{u_1}$  of\\
 $\F^{(0,-1)}_{u_r}\tot \cdots\tot \F^{(0,-1)}_{u_1}$. The $r$-tuple of partitions $(\la_1,\cdots,\la_r)$ labels a torus fixed point of $\cM(n,r)$ with $n=\sum_{i=1}^r|\la_i|$. 
The bases $\ket{\la_r}_{u_r}\tot \cdots \tot \ket{\la_1}_{u_1}$ are conjectured to correspond to the fixed point classes $[(\la_1,\cdots,\la_r)]$ in $\E_\rT(\cM(n,r))$. The vertex operator $\Phi_{\la_1}(u_1)\cdots \Phi_{\la_r}(u_r)$  is realized as a screened vertex operator by taking the elliptic stable envelope of the fixed point $(\la_1,\cdots,\la_r)$ for $\E_\rT(\cM(n,r))$ as the integration kernel\cite{KS}. The $\rK$-theoretic vertex functions obtained from these vertex operators coincide with those in \cite{DinkinsThesis,SD20} obtained geometrically.
\end{itemize}

The aim of this paper is to add the third example, the case $U_{t_1,t_2,p}(\gl_{N,tor})$  and the corresponding  $A^{(1)}_{N-1}$ quiver varieties $\cM(\bfv,\bfw)$, $\bfv, \bfw\in \N^N$.  We propose a  formulation of the vertex operators 
of  $U_{t_1,t_2,p}(\gl_{N,tor})$ by combining its representations given in \cite{KO24} and the formula for the elliptic stable envelopes for the $A^{(1)}_{N-1}$ quiver varieties $\cM(\bfv,\bfw)$\cite{Smirnov18,Dinkins21}.  We construct both the type I, $\Phi_{\bla}(\bfu)$, and the type II dual, $\Psi^*_{\bla}(\bfu)$, vertex operators as screened vertex operators whose integration kernels are given by the elliptic stable envelope for $\E_T(\cM(\bfv,\bfw))$. It turns out that the identity \eqref{IdGeoRep}  connecting representation and geometry ensures the existences of the vertex operators.  
We also show that  compositions of the vertex operators are consistent to the shuffle product formula of the elliptic stable envelopes obtained geometrically following \cite{Botta}. 
This allows us to construct the vertex operators for general $\cM(\bfv,\bfw)$ as compositions of the basic vertex operators $\Phi^{(k)}_\la(u)$ and $\Psi^{*(k)}_\la(u)$\  $(k \in\{0,1,\cdots,N-1\}) $  for $\cM(\bfv,\bfw)$ with one dimensional framing $\bfw=(\delta_{i,k})$.  At the same time this upgrades Conjecture \ref{Nvw} on an equivalence between the level-$(0,-1)$ representation $\F^{(0,-1)(k)}_u$  of $U_{t_1,t_2,p}(\gl_{N,tor})$ and its geometric action on $\bigoplus_{\bfv}\E_\rT(\cM(\bfv,(\delta_{i,k})))$ to Conjecture \ref{level0mw2Mvw}  
between the tensor product of all $\F^{(0,-1)(k)}_{u_j}$, $j=1,\cdots,w_k$, $k=0,1,\cdots,N-1$ and a geometric action 
of the same algebra on $\bigoplus_{\bfv}\E_\rT(\cM(\bfv,\bfw))$ with general framing $\bfw=(w_0,\cdots,w_{N-1})\in \N^N$. 

We then show that the $\rK$-theoretic vertex functions for general $\cM(\bfv,\bfw)$ are obtained as the highest to highest expectation values of  the corresponding  vertex operators.  Furthermore, by defining the $A^{(1)}_{N-1}$ type elliptic dynamical $R$-matrices $R$ and $R^*$ as the transition matrices of the elliptic stable envelopes with different chambers for $\E_T(\cM(\bfv,\bfw))$, we derive exchange relations among the type I and the type II dual vertex operators.  
The associativity of them yields the dynamical Yang-Baxter equation for $R$ and $R^*$. 
Moreover, the exchange relations among the basic vertex operators allow us to construct a $L$-operator $L^{+(k)}(u)$ on the level-$(1,L)$ representation which satisfies the $RLL=LLR^*$ relation.  
Assuming a universal form of the $L$-operator and defining a comultiplication $\Delta$ in terms of it, we show that our basic vertex operators are intertwining operators of the $U_{t_1,t_2,p}(\gl_{N,tor})$-modules w.r.t. $\Delta$.  

It is also worth to mention that the above Conjecture 1.1 (3) is similar to Theorem 3.1 in \cite{AFO}, which relates the $\rK$-theoretic vertex functions to correlation functions of the deformed $W$-algebras\cite{FrRe}.  
In this sense, the above Example 2) and the result in Sec.\ref{sec:VFs}  below give extensions of this theorem to the affine quiver variety cases.  
In fact, we have a conjecture that the elliptic quantum groups $U_{q,p}(\gh)$ and $U_{t_1,t_2,p}(\gl_{tor})$ give an alternative formulation of deformed $W$-algebras, which are 
deformation of the GKO coset type $W$-algebra associated with $\gh$\cite{GKO,CR,KMQ,LF,Rav}  for $U_{q,p}(\gh)$ and the affine quiver $W$-algebras\cite{KimPes} for $U_{t_1,t_2,p}(\gl_{tor})$, respectively. Note that the coset type $W$-algebra associated with a simply laced $\gh$ is equivalent to the one obtained by a quantum Hamiltonian reduction but this is not the case for the non-simply laced cases, see for example \cite{BS}.  A formulation in \cite{FrRe} gives a deformation of the quantum Hamiltonian reduction type. 
A typical example of the non-simply-laced coset type $W$-algebra is the $W\hspace{-1mm}B_{N}$ algebra introduced by Fateev and Lukyanov\cite{LF}.  Its deformation has been proposed  in \cite{K14,FKO}.  
In wider situation, it has been shown that the elliptic currents of $U_{q,p}(\gh)$ and $U_{t_1,t_2,p}(\gl_{tor})$ give 
screening currents of the corresponding deformed $W$-algebras\cite{K98,JKOS,KK03,FKO,KOgl1},
 the level-1 highest weight $U_{q,p}(\gh)$-module is decomposed  into a direct sum of the deformed $W$-algebra modules\cite{FKO,KKW}, a composition of the vertex operator w.r.t. the Drinfeld comultiplication and its shifted inverse gives a generating function of the deformed $W$-algebra\cite{K14,KOgl1}. 
Moreover we expect that the vertex operators of EQG/EQTA give  deformation of primary fields of the corresponding $W$-algebras\cite{LP,Konno24}. Our construction of the K-theoretic vertex function and Theorem 3.1 in \cite{AFO} support that this expectation is true.

This paper is organized as follows. 
In Sec.2, we review  the EQTA $U_{t_1,t_2,p}(\gl_{N,tor})$ including its  
 $H$-Hopf algebroid structure. 
 Sec.\ref{sec:Reps} is devoted to a summary of representations of 
$U_{t_1,t_2,p}(\gl_{N,tor})$, the level $(1,L)$- and $(0,-1)$-representations.  
  In Sec.\ref{sec:ESE}, we present some known results on the elliptic stable envelopes 
   for $\E_\rT(\cM(\bfv,\bfw))$. 
  In Sec.\ref{sec:VO}, we propose a  construction of  the type I and the type II dual vertex operators of $U_{t_1,t_2,p}(\gl_{N,tor})$ and show their consistency to  the shuffle product formula of the elliptic stable envelopes. 
 Sec.\ref{sec:VFs} is devoted to a derivation of the $\rK$-theoretic vertex functions for $\cM(\bfv,\bfw)$.  
In Sec.\ref{sec:ExchRs}, we derive exchange relations among  the type I and the type II dual  vertex operators. 
Finally, in Sec.\ref{sec:Lop}, we construct a $L$-operator $L^{+(k)}(u)$. Assuming its universal form,  we discuss intertwining properties of the basic vertex operators.

\section{Elliptic Quantum Toroidal Algebras $U_{t_1,t_2,p}(\glntbig)$}\lb{app:1}
We summarize some basic facts on the elliptic quantum toroidal algebra $U_{t_1,t_2,p}(\gl_{N,tor})$ introduced in  \cite{KO24}. 

\subsection{Definition
 }

Let $t_1, t_2, p  \in \C^\times$ be generic complex numbers with $|t_1|, |t_2|, |p|<1$.  We set
\bea
&&\hbar=t_1t_2=q^2,\qquad t_1/t_2=\kappa^{-2}.\lb{hbar}
\ena
These are related to the parameters $\kappa, q_1, q_2, q_3$ used in  \cite{KO24} by
  $t_1=\kappa^{-1}q=q_1^{-1}$, $t_2= \kappa q=q_3^{-1}, \hbar=q_2$. 
We set
\begin{eqnarray}
&&[n]_q=\frac{q^n-q^{-n}}{q-q^{-1}}, \qquad [n]_q !=[n]_q \cdots [1]_q, \qquad \left[ n \atop k \right]_q=\frac{[n]_q !}{[n-k]_q ! [k]_q !}.
\end{eqnarray}

For a natural number $N\geq 3$, we set 
$I=\{0,1, \dots, N-1\}$. Let $\slnh$ be an untwisted affine Lie algebra over $\C$ with  the generalized Cartan matrix $A=(a_{ij})_{i,j \in I}$, $a_{ij}=2\delta_{ij}-\delta_{j i-1}-\delta_{j i+1}$. All indices in $I$ should be considered mod $N$. 
We fix a realization $(\goth{h}, \Pi, \Pi^{\vee})$ of $A$, that is, $\goth{h}$ is an $N+1$-dimensional $\mathbb{Q}$-vector space,
$\Pi=\{\alpha_0,\alpha_1,\dots,\alpha_{N-1}\} \subset \goth{h}^*$ is a set of simple roots, 
and $\Pi^{\vee}=\{ h_0, h_1,\cdots,h_{N-1}
\} \subset \goth{h}$ is a set of simple coroots, 
which are defined such that
 $\langle\alpha_j, h_i\rangle=a_{ij}\ (i,j \in I)$  for a canonical pairing $\langle\, ,\rangle : \h^*\times \h \to \C$.  We also set $\cQ=\sum_{i\in I}\Z \al_i$. Let $c=\sum_{i\in I}h_i$ and $\delta=\sum_{i\in I}\al_i$.  We take $\{h_1,\dots,h_{N-1}, c, d\}$  as the basis of $\h$ and $\{ {\Lambda}_1,\cdots, {\Lambda}_{N-1}, \Lambda_0, \delta \}$ the dual basis of $\h^*$ satisfying
\bea
&&
\langle\delta,d\rangle=1=\langle\Lambda_0,c\rangle,\quad 
\langle{\Lambda}_i,h_j\rangle=\delta_{i,j}\lb{pairinghhs}, 
\ena
with the other pairings being 0. 
The toroidal algebra $\gl_{N,tor}$ has one more central element $c^\perp$. We set $\h_{tor}=\h\oplus\C c^\perp$. Accordingly we add   
an element $\La^\perp_0$ to $\h^*$ and set $\h^*_{tor}=\h^*\oplus\C \La^\perp$. 
We have a paring $\langle \La^\perp_0,c^\perp\rangle=1$.

Let $H_P$ be  a $\C$-vector space spanned by $ P_0, P_1,\cdots,P_{N-1}$
 and $H^{*}_P$ dual space spanned by $Q_0, Q_1,\cdots,Q_{N-1}$ with 
the pairing $\langle Q_i,P_j\rangle =a_{ij}\ (i,j \in I)$.  For $\la=\sum_i c_i\al_i\in \h^*$, we set $Q_\la=\sum_ic_iQ_i$.  
Set $H=\h_{tor}\oplus H_P$ and $H^*=\h^*_{tor}\oplus H_P^*$, and extend the above pairings to $H\times H^*$ with $\langle \h^*_{tor},H_P\rangle =0=\langle H_P^*,\h_{tor}\rangle$.  
We set $\gz_i=\hbar^{P_i+h_i}$ and $\gz^*_i=\hbar^{P_i}$, $i\in I$ and denote by $\cM$ (resp.  $\cM^*$) the field of meromorphic functions of $\gz_i\ (i\in I)$ (resp. $\gz_i^*\  (i\in I)$). 
We regard  $g(\gz_i)\in \cM$  and $g(\gz_j^*)\in \cM^*$ as functions on  ${H^*}$ by 
$g(\gz_i)(\mu)=g(\hbar^{\langle \mu,P_i+h_i\rangle}), g(\gz_j^*)(\mu)=g(\hbar^{\langle \mu,P_j\rangle})$ for $\mu\in H^*$. 

We regard $H\oplus H^*$ as a Heisenberg algebra by imposing the commutation relation 
$[x,y]=\langle y,x\rangle$, $x\in H, y\in H^*$, for example 
\bea
&&[h_i,\al_j]=a_{ij}=[P_i,Q_j],\quad [c,\La_0]=1=[c^\perp, \La^\perp_0],\quad [h_{i},\La_j]=\delta_{ij}. \lb{def:La}
\ena

\begin{dfn}
Let  $m_{ij}=\delta_{j,i-1}-\delta_{j,i+1}\ (i,j \in I)$. 
The elliptic quantum toroidal algebra $U_{t_1,t_2,p}(\glnt)$ is an topological algebra over $\FF=\cM[[p]]\times \cM^*[[p^*]]$
generated by 
\be
\al_{i,l}, \quad x^+_{i,n}, \quad x^-_{i,n}, \quad K^\pm_i,
\quad \gamma^{\pm 1/2} \quad 
 (i\in I, \quad  l \in \mathbb{Z}\backslash\{0\}, \quad n \in \mathbb{Z}). 
\en
Here we set $p^*=p\gamma^{-2}$.  
The defining  relations can be  written conveniently in terms of  the generating functions called the elliptic currents 
given by
\be
&&x^+_i(z)=\sum_{n\in \Z}x^+_{i,n} z^{-n},\qquad x^-_i(z)=\sum_{n\in \Z}x^-_{i,n} z^{-n},\\
&&\phi^+_i(\gamma^{\frac{1}{2}}z)
=K_i^{+ } \exp\left( -(q-q^{-1}) \sum_{n>0} \frac{p^n \al_{i, -n}}{1-p^n}z^{ n}\right)
 \exp\left( (q-q^{-1}) \sum_{n>0} \frac{\al_{i, n}}{1-p^n}z^{ -n}\right),\\
&&\phi^-_i(\gamma^{-\frac{1}{2}}z)
=K_i^{-} \exp\left( -(q-q^{-1}) \sum_{n>0} \frac{ \al_{i, -n}}{1-p^n}z^{ n}\right)
 \exp\left( (q-q^{-1}) \sum_{n>0} \frac{p^n \al_{i, n}}{1-p^n}z^{ -n}\right).
\en
Set also  
\bea
&&K:=\left(\prod_{i\in I}K^+_i(K^-_i)^{-1}\right)^{1/2}. 
\ena
The defining relations are given as follows. 
For any $g(\gz_k,\gz_l^*) \in \FF
$, 
\bea
&& \gamma^{ 1/2}, \quad K\  :\hbox{ central },\\
&&
g(\gz_k,\gz_l^*)x^+_j(z)=x^+_j(z)g(\gz_k,\gz_l^*\hbar^{-a_{lj}}),\qquad
\lb{ge}\\
&&
g(\gz_k,\gz_l^*)x^-_j(z)=x^-_j(z)g(\gz_k\hbar^{-a_{kj}},\gz_l^*),\lb{gf}\\
&&[g(\gz_k,\gz_l^*), \al_{i,m}]=
0,\qquad [g(\gz_k,\gz_l^*),d]=0,\lb{gboson}\\ 
&&g(\gz_k,\gz_l^*)K^\pm_j=K^\pm_jg(\gz_k\hbar^{-a_{kj}},\gz_l^*\hbar^{-a_lj}),
\lb{gKpm}
\\
&& [K^\pm_i,\al_{j,k}]=0, \quad K^\pm_i x^+_j(z)=q^{\pm a_{ij}
  }x^+_j(z)K^\pm_i, \quad K^\pm_i x^-_j(z)=
q^{\mp a_{ij}
}x^-_j(z)K^\pm_i,\\
&&[\al_{i,l},\al_{j,m}]=\delta_{l+m,0}\frac{[a_{ij} l]_q}{l}\frac{\gamma^l-\gamma^{-l}}{q-q^{-1}}\frac{1-p^l}{1-p^{*l}}
\kappa^{-lm_{ij}}\gamma^{-l}, \\
&&
[\al_{i,l},x^+_j(z)]=\frac{[a_{ij} l]_q}{l}\frac{1-p^l}{1-p^{*l}}\gamma^{-l}\kappa^{-lm_{ij}}z^l x^+_j(z),\\
&&
[\al_{i,l},x^-_j(z)]=-\frac{[a_{ij} l]_q}{l}  \kappa^{-lm_{ij}}z^l x^-_j(z),
\\
&&(\kappa^{m_{ij}} z-q^{a_{ij}}w)g^+_{ij}(\kappa^{-m_{ij}}w/z;p^*)
x^+_i(z)x^+_j(w)
=(\kappa^{m_{ij}}q^{a_{ij}}z-w) g^+_{ij}(\kappa^{m_{ij}}z/w;p^*)
x^+_j(w)x^+_i(z),\nn\\
&&\lb{xpxpgen}\\
&&(\kappa^{m_{ij}} z-q^{-a_{ij}}w)g^-_{ij}(\kappa^{-m_{ij}}w/z;p)
 x^-_i(z)x^-_j(w)
=(\kappa^{m_{ij}}q^{-a_{ij}}z-w) g^-_{ij}(\kappa^{m_{ij}}z/w;p)
x^-_j(w)x^-_i(z),\nn\\
&&\lb{xmxmgen}
\\
&&[x^+_i(z),x^-_j(w)]=\frac{\delta_{i,j}}{q-q^{-1}}
\left(\delta\bigl(\gamma {w}/{z}\bigr)\phi^{+}_{i}(\gamma^{\frac{1}{2}} w)
-\delta\bigl(\gamma^{-1} {w}/{z}\bigr)\phi^{-}_{i}(\gamma^{\frac{1}{2}}z)
\right),\lb{xpxm}
\ena
\bea
&&\sum_{\sigma\in \gS_a}\prod_{1\leq k<m\leq a
}g^+_{ii}(z_{\sigma(m)}/z_{\sigma(k)};p^*)
\nonumber\\
&&\qquad\quad\times
\sum_{r=0}^a (-1)^r
\left[\mmatrix{a\cr r\cr}\right]_q \prod_{1\leq s\leq r}
g^+_{ij}(\kappa^{-m_{ij}}w/z_{\sigma(s)};p^*)
\prod_{r+1\leq s\leq a} 
g^+_{ij}(\kappa^{m_{ij}}z_{\sigma(s)}/w;p^*)
\nn \\
&&\qquad\quad\times x^+_{i}(z_{\sigma(1)})\cdots x^+_{i}(z_{\sigma(r)})
x^+_{j}(w) x^+_{i}(z_{\sigma(r+l)})\cdots x^+_{i}(z_{\sigma(a)})=0,\lb{Serrexp}\\
&&\sum_{\sigma\in \gS_a}\prod_{1\leq k<m\leq a
}g^-_{ii}(z_{\sigma(m)}/z_{\sigma(k)};p)
\nonumber\\
&&\qquad\quad\times
\sum_{r=0}^a (-1)^r
\left[\mmatrix{a\cr r\cr}\right]_q \prod_{1\leq s\leq r}
g^-_{ij}(\kappa^{-m_{ij}}w/z_{\sigma(s)};p)
\prod_{r+1\leq s\leq a} 
g^-_{ij}(\kappa^{m_{ij}}z_{\sigma(s)}/w;p)
\nn\\
&&\qquad\quad\times x^-_{i}(z_{\sigma(1)})\cdots x^-_{i}(z_{\sigma(r)})
x^-_{j}(w) x^-_{i}(z_{\sigma(r+l)})\cdots x^-_{i}(z_{\sigma(a)})=0.
\ena
Here we set 
\bea
&&g^{\pm}_{ij}(z;s)=\exp\left(\mp\sum_{m>0}\frac{1}{m}\frac{q^{a_{ij}m}-q^{-a_{ij}m}}{1-s^m}(sz)^m\right)\ \in \ \C[[s]][[z]],\lb{deg:g}
\ena
for $s=p, p^*$. 
\end{dfn}
We treat these relations as formal Laurent series in the argument of the elliptic currents. 
Their coefficients are well defined in the $p$-adic topology\cite{Konno16}.

It is  convenient to set 
\bea
&&\al'_{i,l}=\frac{1-p^{*l}}{1-p^l}\gamma^{l} \al_{i,l}. 
\ena
One then has the following commutation relations. 
\begin{prop}\lb{lem:sec4}
For $m,\,n \in \mathbb{Z}_{\ne 0}$, 
\begin{align}
& [\alpha'_{i,m},\alpha'_{j,n}]=\delta_{m+n,0}\frac{[a_{ij}m]_q}{m} \frac{\gamma^{m}-\gamma^{-m}}{q-q^{-1}}
\frac{1-p^{*m}}{1-p^{m}}\kappa^{-m\, m_{ij}}\gamma^{m},\\
& [\alpha_{i,m},\alpha'_{j,n}]=[\alpha'_{i,m},\alpha_{j,n}]=\delta_{m+n,0}\frac{[a_{ij}m]_q}{m}\frac{\gamma^{m}-\gamma^{-m}}{q-q^{-1}}\kappa^{-m\,m_{ij}}, \\
& [\alpha'_{i,m}, x_j^+(z)]=\frac{[a_{ij}m]_q}{m}\kappa^{-m\,m_{ij}}z^m x_j^+(z), \\
& [\alpha'_{i,m}, x_j^-(z)]=-\frac{[a_{ij}m]_q}{m}\frac{1-p^{*m}}{1-p^m}\gamma^{m}\kappa^{-m\,m_{ij}} z^m x_j^-(z).
\end{align}
\end{prop}

\medskip
For representations, on which $\gamma^{\pm 1/2}$ take complex values 
such as  in Sec.\ref{Sec:Level1Mrep} and Sec.\ref{Sec:Level0-1rep}, 
 we treat $p^*$ as a  complex number satisfying $|p^*|=|p\gamma^{-2}|<1$. 
Then  one has 
\be
 &&g^\pm_{ij}(z;s)=\frac{(sq^{\pm a_{ij}}z;s)_\infty}{(sq^{\mp a_{ij}}z;s)_\infty}
\en
for $|sq^{\pm a_{ij}}z|<1$, $s=p, p^*$, where we set
\be
&&(z;s)_\infty=\prod_{n=0}^\infty(1-zs^n).
\en
In the sense of analytic continuation, one can rewrite, for example \eqref{xpxpgen} and \eqref{xmxmgen}, as
\bea
&&\theta^*(q^{a_{ij}}\kappa^{-m_{ij}}w/z)x^+_i(z)x^+_j(w)
=-\theta^*(q^{a_{ij}}\kappa^{m_{ij}}z/w)x^+_j(w)x^+_i(z),\lb{xpxp}\\
&&\theta(q^{-a_{ij}}\kappa^{-m_{ij}}w/z)x^-_i(z)x^-_j(w)
=-\theta(q^{-a_{ij}}\kappa^{m_{ij}}z/w)x^-_j(w)x^-_i(z).\lb{xmxm}
\ena
Here $\theta(z), \theta^*(z)$ denote Jacobi's odd theta functions given by 
\bea
&&\theta(z)=-z^{-1/2}\theta_p(z),\qquad \theta_p(z)=(z;p)_\infty(p/z;p)_\infty,\lb{def:theta}\\
&&\theta^*(z)=-z^{-1/2}\theta_{p^*}(z).\nn
\ena

\subsection{$\h$-Hopf Algebroid Structure}
The $U_{t_1,t_2,p}(\gl_{N,tor})$ has a $\h$-Hopf algebroid structure as a coalgebra structure w.r.t. the Drinfeld comultiplication. See Sec.3.2 in \cite{KO24}.  
In this paper, we do not use this precise structure except for those explained below. 

\begin{prop}
Let 
$\cU=U_{t_1,t_2,p}(\gl_{N,tor})$  is {a $\h$-algebra} by 
\be
&&\cU=\bigoplus_{\al,\beta\in \h^*}\cU_{\al,\beta}\\
&&(\cU)_{\al\beta}=\left\{x\in \cU \left|\ \mmatrix{\gz_i x\gz_i^{-1}
=\hbar^{\langle\al,P_i+h_i\rangle}x,\quad 
\gz^*_i x\gz_i^{*-1}=\hbar^{\langle Q_\beta,P_i\rangle}x\qquad \forall i\in I\ \cr
 }
 \right.\right\}
\en
and $\mu_l, \mu_r : \FF \to \cU_{0,0}$ defined by 
\be
&&\mu_l(\hf)=f_1(\bgz,p)f_2(\bgz,p),\qquad \mu_r(\hf)=f_1(\bgz^*,p^*)f_2(\bgz^*,p^*)
\en
for $\hf=f_1(\bgz,p)f_2(\bgz^*,p^*)\in \FF$, where $\bgz=\{\gz_i\}_{i\in I}, \ \bgz^*=\{\gz^*_i\}_{i\in I}$.
\end{prop}
Let $\cR_Q=\sum_{i\in I}\Z Q_i$. For $\al=\sum_in_i\al_i$, let $Q_\al=\sum_in_iQ_i\in \cR_Q$ and  
consider its group algebra $\C[\cR_Q]$ with $e^{Q_\al}, e^{Q_\beta}, e^{Q_\al}e^{Q_\beta}=e^{Q_\al+Q_\beta}, e^0=1\in \C[\cR_Q]$. 
We regard $T_\al=e^{-Q_\al}\in \C[\cR_Q]$ as a shift operator 
\be
&&(T_\al \widehat{f})=f_1(\{\gz_i\hbar^{\langle Q_\al,P_i\rangle}\}, p)f_2(\{\gz_j^*\hbar^{\langle Q_\al,P_j\rangle}\},p^*).
\en
Then $\cD=\{\widehat{f} e^{-Q_\al} \ |\  \widehat{f}\in \FF,  e^{-Q_\al}\in \C[\cR_Q]\}$ becomes a $\h$-algebra satisfying
\be
&&\cU\cong \cU\tot \cD\cong  \cD\tot \cU. \lb{Diso} 
\en

The tensor product $\cA {\widetilde{\otimes}}\cB$ of the two $\h$-algebras $\cA$ and $\cB$ is the $\h$-bigraded vector space with 
\be
 (\cA {\widetilde{\otimes}}\cB)_{\al\beta}=\bigoplus_{\gamma\in\h^*} (\cA_{\al\gamma}\otimes_{\FF}\cB_{\gamma\beta}),
\en
where $\otimes_{\FF}$ denotes the usual tensor product 
modulo the following relation.
\bea
\mu_r^\cA(\hf) a\otimes b=a\otimes\mu_l^\cB(\hf) b, \qquad a\in \cA,\  
b\in \cB,\ \hf\in \FF.\lb{reltot}
\ena
The tensor product $\cA {\widetilde{\otimes}}\cB$ is again a $\h$-algebra with the multiplication $(a\otimes b)(c\otimes d)=ac\otimes bd$ and the moment maps 
\be
\mu_l^{\cA {\widetilde{\otimes}}\cB} =\mu_l^\cA\otimes 1,\qquad \mu_r^{\cA {\widetilde{\otimes}}\cB} =1\otimes \mu_r^\cB.
\en

\section{Representations}\lb{sec:Reps}

We summarize the level $(1,L)$- and $(0,-1)$-representations of  $U_{t_1,t_2,p}(\gl_{N,tor})$ 
 constructed  in \cite{KO24}. 

\subsection{Level of representation}
Let us consider a vector space $\hV$ over $\FF$, which is  
${\h}$-diagonalizable, i.e.  
\be
&&\hV=\bigoplus_{\la,\mu\in {\h}^*}\hV_{\la,\mu},\ 
\hV_{\la,\mu}=\left\{ v\in \hV\ \left|\ \mmatrix{
\gz_i\cdot v=\hbar^{\langle\la,P_i+h_i\rangle} v,\ 
\gz_i^*\cdot v=\hbar^{\langle Q_\mu,P_i\rangle} v\quad \forall i\in I\ \cr
 }
 \right.\right\}.
\en
Let us define the $\h$-algebra $\cD_{\h,\hV}$ of the $\C$-linear operators on $\hV$ by
\be
&&\cD_{\h,\hV}=\bigoplus_{\al,\beta\in {\h}^*}(\cD_{\h,\hV})_{\al\beta},\\
&&\hspace*{-10mm}(\cD_{\h,\hV})_{\al\beta}=
\left\{\ X\in \End_{\C}\hV\ \left|\ 
\mmatrix{ \gz_i X \gz_i^{-1}
=\hbar^{\langle\alpha,P_i+h_i\rangle} X ,
\quad \gz^*_i X \gz_i^{*-1}
=\hbar^{\langle Q_\beta,P_i\rangle}X{,}\cr
 X\cdot\hV_{\la,\mu}\subseteq 
 \hV_{\la+\al,\mu+\beta}\cr}  \right.\right\},\\
&&\mu_l^{\cD_{\h,\hV}}(\widehat{f})v=f_1(\{\hbar^{\langle\la,P_i+h_i\rangle}\}, p)f_2(\{\hbar^{\langle\la,P_i+h_i\rangle}\},p)v,
\\ 
&&\mu_r^{\cD_{\h,\hV}}(\widehat{f})v=f_1(\{\hbar^{\langle Q_\mu,P_i\rangle}\},p^*)f_2(\{\hbar^{\langle Q_\mu,P_i\rangle}\},p^*)v,\quad \forall \widehat{f}\in \FF,
\quad\forall v\in \hV_{\la,\mu}.
\en
\begin{dfn}
We define a dynamical representation of $\cU=U_{t_1,t_2,p}(\glnt)$ on $\hV$ to be  
 a $H$-algebra homomorphism ${\pi}: \cU
 \to \cD_{\h,\hV}$. 
\end{dfn}
\begin{dfn}
For $k, l\in \C$, we say that a $\cU$-module has  level $(k,l)$ if the two central elements $\gamma$ and $K$ act 
as $\hbar^{k/2}$ and $\hbar^{l/2}$ on it, respectively.  
\end{dfn}

\begin{dfn}
Let $\cU({\frak H})$, $\cU({\frak N}_+), \cU({\frak N}_-)$ be the subalgebras of $\cU$ generated by 
$\gamma^{\pm 1/2}, 
K^\pm_{i}\\ (i\in I)$, by $\al_{i,n}\ (i\in I, n\in \Z_{>0})$,  $x^+_{i, n}\ (i\in I, n\in \Z_{\geq 0})$  
$x^-_{i, n}\ (i\in I, n\in \Z_{>0})$ and by $\al_{i,-n}\ (i\in I, n\in \Z_{>0}),\ x^+_{i, -n}\ (i\in I, n\in \Z_{> 0}),\ x^-_{i, -n}\ (i\in I, n\in \Z_{\geq 0})$, respectively.   
\end{dfn}

\begin{dfn}
For $k, l\in\C$, $\la, \mu\in \h^*$, 
a dynamical $\cU$-module $V(\la,\mu)$ is called the 
level-$(k,l)$ highest weight module with the highest weight $(\la,\mu)$, if there exists a vector 
$v\in V(\la,\mu)$ such that\footnote{We here changed the definition of the level of representation from the one given in \cite{KO24} so that our level-$(k,l)$ is the level-$(k,-l)$ there. }
\be
&& \gz_i\cdot v=
\hbar^{\langle \la,P_i+h_i\rangle}v,\qquad \gz^*_i\cdot v=\hbar^{\langle Q_\mu,P_i\rangle}v,\\
&&V(\la,\mu)=\cU\cdot v,\qquad \cU({\frak N}_+)\cdot v=0,\\
&&\gamma\cdot v=\hbar^{k/2}v,\qquad  K\cdot v
=\hbar^{l/2}v.
\en
\end{dfn}

\subsection{Level-(1,L) representation}\lb{Sec:Level1Mrep} 
It is convenient  to introduce $\bep_j\ (1\leq j\leq N)$ defined by
${\al}_j=\bep_{j}-\bep_{j+1}$ $(1\leq j\leq N-1)$,  $\al_0=\delta+\bep_N-\bep_1$  and $\sum_{j=1}^{N}\bep_j=0$. One has 
\bea
&&\bep_i=\frac{1}{N}\left(-\sum_{j=1}^{i-1}j\al_{j}+\sum_{j=i}^{N-1}(N-j)\al_j\right),\\
&&\langle \bep_i,h_j\rangle=\delta_{i,j}-\delta_{i,j+1}\qquad 
(1\leq i\leq N,\ 0\leq j\leq N-1 ). 
\ena
One then obtains the following realization of $\Lambda_i$ satisfying \eqref{pairinghhs}.  
\bea
&&\Lambda_i=\Lambda_0+\sum_{j=1}^{i}\bep_{j}\quad (i\in I
).
\ena
Set also  
\bea
&&[c,\delta]=0=[c,\bep_i]=[h_j,\delta]=[h_j,\La_0], \qquad (1\leq i\leq N,\ 1\leq j\leq N-1),\\
&&[h_j,\bep_i]=\delta_{i,j}-\delta_{i,j+1}.
\ena  
One recovers
$[h_i,\al_j]=a_{ij}$ and  $[h_i,\Lambda_j]=\delta_{i,j}$, $(i,j\in I)$.
We also use $h_{\bep_i}\ (1\leq i\leq N)$ and $h_{\La_j}\ (j\in I)$ defined by
\bea
&&h_k=h_{\bep_k}-h_{\bep_{k+1}}\ (1\leq k\leq N-1),\quad h_0=c+h_{\bep_N}-h_{\bep_1},\quad \sum_{i=1}^Nh_{\bep_i}=0,\\
&&h_{\La_j}=h_{\La_0}+\sum_{i=1}^jh_{\bep_i}.
\ena
Setting 
\bea
&&[h_{\La_0},\delta]=1,\quad [h_{\bep_i},\delta]=0=[h_{\La_0},\bep_i],
\ena
one obtains the following commutation relations.
\begin{prop} 
\bea
&&[h_{\La_i},\al_j]=\delta_{i,j},\qquad [h_{\bep_i},\bep_j]=\delta_{i,j}-\frac{1}{N},\lb{ComhLaal}\\
&&[h_{\La_i},\La_j]=\left\{\mmatrix{\frac{1}{N}i(N-j)&\quad (0\leq i\leq j\leq N-1)\cr
\frac{1}{N}j(N-i)&\quad (0\leq j\leq i\leq N-1)\cr}\right. .
\ena 
\end{prop}
We set $\eta_{ij}=[h_{\La_i},\La_j]$. One has $\eta_{ij}=\eta_{ji}$. 

In the same way, 
we introduce  $P_c,\ Q_\delta$,  $P_{\bep_i},\ Q_{\bep_i}\ (1\leq i\leq N)$ and $ Q_{\La_j}\ (j\in I)$  by
\bea
&&P_i=P_{\bep_i}-P_{\bep_{i+1}}\ (1\leq i\leq N-1),\quad P_0=P_c+P_{\bep_N}-P_{\bep_1},\quad \sum_{i=1}^NP_{\bep_i}=0,\\
&&\hspace{-0.3cm}Q_i=Q_{\bep_i}-Q_{\bep_{i+1}}\ (1\leq i\leq N-1),\quad Q_0=Q_\delta+Q_{\bep_N}-Q_{\bep_1},\quad \sum_{i=1}^NQ_{\bep_i}=0,
\ena
satisfying 
\bea
&&[P_{\bep_i},Q_{\bep_j}]=\delta_{i,j}-\frac{1}{N},\\
&&[P_c,Q_\delta]=0=[P_c,Q_{\bep_i}]=[P_{\bep_i},Q_\delta]\qquad \forall 1\leq i,j\leq N. 
\ena
Then one has 
\bea
&&[P_i,Q_{\bep_j}]=\delta_{i,j}-\delta_{i+1,j}
\ena
and recovers  $[P_i,Q_j]=a_{ij}$.  Moreover define $Q_{\La_j}\ (j\in I)$ by
\bea
&& Q_{\La_j}=Q_{\La_0}+\sum_{i=1}^jQ_{\bep_i}
\ena
and require 
\bea
&&[P_c,Q_{\La_0}]=1, \quad [P_{\bep_i},Q_{\La_0}]=0.
\ena
Then one has 
\bea
&&[P_i,Q_{\La_j}]=\delta_{i,j}\qquad i,j\in I.
\ena

Now let us consider 
the group algebra $\C[\cQ]$ with $e^{\al}, e^{\beta}, e^{\al}e^{\beta}=e^{\al+\beta}, e^0=1\in \C[\cQ]$.  
We introduce the following central extension to  $\C[\cQ]$\cite{KO24}\footnote{This makes us avoid the doubling trick by Saito used in Prop.3.2.2 \cite{Saito}. There two elements $H_{i,0}$ and $\partial_{\bar{\al_i}}$ are used to realize the Heisenberg generator corresponding to our ${h}_i$ in $x^\pm_i(z)$ and in $K^\pm_i$, respectively.  Both $H_{i,0}$ and $\partial_{\bar{\al_i}}$ measure the same weight as ${h}_i$, but the former  produces an extra factor $\kappa^{m_{ij}}$
 in the exchange relations among $x^+_i(z)$ and $x^+_j(w)$. 
 }. 
\bea
&&e^{\al_i}e^{\al_j}=(-\kappa)^{-m_{ij}}e^{\al_j}e^{\al_i}.\lb{CEalal}
\ena

\begin{thm}\cite{KO24}\lb{Thm:level1Mrep}
Let $\bxi=(\xi_0,\xi_1,\cdots,\xi_{N-1})\in (\C^\times)^{N}$,  $\bfL=(L_0,L_1,\cdots,L_{N-1}), \bfM=(M_0,M_1,\cdots,M_{N-1})\in \Z^{N}$ and set $L=\sum_{i\in I}L_i$ and $\blLa_L=\sum_{i\in I}L_i\La_i$, 
$\blLa_M=\sum_{i\in I}M_i\La_i$.  
The following gives the level-$(1,L)$ highest weight representation of $U_{t_1,t_2,p}(\glnt)$ with the highest weight $(\blLa_L+\blLa_M,\blLa_M)$. 
\be
&&\gamma^{1/2}=\hbar^{c^\perp/4},\qquad K^\pm_i=q^{\pm h_i}e^{-Q_i},\\
&&\phi^+_i(q^{1/2}z)=K^+_i\exp\left\{-(q-q^{-1})\sum_{m>0}\frac{p^m}{1-p^m}\al_{i,-m}z^m\right\}
\exp\left\{(q-q^{-1})\sum_{m>0}\frac{1}{1-p^m}\al_{i,m}z^{-m}\right\},
\\
&&\phi^-_i(q^{-1/2}z)=K^-_i\exp\left\{-(q-q^{-1})\sum_{m>0}\frac{1}{1-p^m}\al_{i,-m}z^m\right\}
\exp\left\{(q-q^{-1})\sum_{m>0}\frac{p^m}{1-p^m}\al_{i,m}z^{-m}\right\},
\\
&&x^+_i(z)=\exp\left\{\sum_{m>0}\frac{\al_{i,-m}}{[m]_q}z^m\right\}
\exp\left\{-\sum_{m>0}\frac{\al_{i,m}}{[m]_q}z^{-m}\right\}\ e^{\al_i}e^{-Q_i} z^{h_i+1},\\
&&x^-_i(z)=\exp\left\{-\sum_{m>0}\frac{\al'_{i,-m}}{[m]_q}z^m\right\}
\exp\left\{\sum_{m>0}\frac{\al'_{i,m}}{[m]_q}z^{-m}\right\}\ e^{-\al_i}z^{-h_i+1}.
\en
The representation space is given by
\bea
&&\F^{(1,L)}_{\bxi,\blLa_L,\blLa_M}= \FF[\al_{i,-m}\ (i\in I,m\in \Z_{>0})]\prod_{i\in I}\xi_i^{-\sum_{j\in I}\eta_{ij}L_j}\otimes e^{\La^\perp_0}\otimes e^{\blLa_L}\C[Q]\otimes e^{Q_{\blLa_M}}\C[\cR_Q],\nn\\
&&
\ena
where we set $Q_{\blLa_M}=\sum_{i\in I}M_iQ_{\La_i}$. 
The highest weight vector is given by
\bea
&&\ket{0}^{(1,L)}_{\bxi,\blLa_L,\blLa_M}=\prod_{i\in I}\xi_i^{-\sum_{j\in I}\eta_{ij}L_j}1\otimes 
e^{\La_0^\perp}\otimes e^{\blLa_L}\otimes e^{Q_{\blLa_M}}.\lb{def:vac1M} 
\ena
The action is given by
\be
&&
\gamma\cdot u=\hbar^{1/2}u,\qquad 
K\cdot u=\hbar^{L/2}u,  
\\
&&\alpha_{i,-n}\cdot u= \alpha_{i,-n}u,\qquad \\
&&\alpha_{i,n}\cdot u=\sum_{j} \frac{[a_{ij}n]_q[n]_q}{n}\frac{1-p^{n}}{1-p^{*n}}\kappa^{-nm_{ij}}q^{-n}\frac{\partial}{\partial\alpha_{j,-n}}u,\qquad \forall u\in\F^{(1,L)}_{\bxi,\blLa_L,\blLa_M},\ n>0.
\en
In particular one has  $p^*=p\hbar^{-1}$ on $\F^{(1,L)}_{\bxi,\blLa_L,\blLa_M}$. 
\end{thm}

\noindent
{\it Proof.} Noting $K=\hbar^{c/2}$ 
 and using \eqref{def:La} 
and $\hbar^{c/2}e^{L_i\La_i}=\hbar^{L_i/2}e^{L_i\La_i}$. one has $K\cdot u=\hbar^{L/2}u$.  
 The rest of the statement follows from Theorem 5.13 in \cite{KO24}
 \footnote{In \cite{KO24}, 
it is assumed that $e^{Q_i}\cdot 1=1$. Note that $\sum_{i\in I}Q_i$ commutes with $U_{t_1,t_2,p}(\gl_{N,tor})$ and $\cM\times \cM^*$. In addition, we changed the notation : $c$ in \cite{KO24} is $c^\perp$ here.  } . 
\qed

\subsection{Level-(0,-1) representation}\lb{Sec:Level0-1rep}

For $k\in I, u\in \C^\times$, let $V^{(k)}(u)$ be a vector space over $\FF$ with basis $[u]^{(k)}_j, j\in \Z$. 
Define the action of $\phi^\pm_i(z), x^+_i(z), x^-_i(z)\ i\in I$ by
\bea
&&\phi^\pm_i(z)[u]^{(k)}_j=\left\{\mmatrix{\frac{\theta(t_1^{-j}u/\hbar z)}{\theta(t_1^{-j}u/z)}\biggl|_{\pm}e^{-Q_i}[u]^{(k)}_j,\qquad& i+j\equiv k\cr
\frac{\theta(t_1^{-j}t_2u/z)}{\theta(t_1^{-j}u/t_1z)}\biggl|_{\pm}e^{-Q_i}[u]^{(k)}_j,\qquad& i+j+1\equiv k\cr
e^{-Q_i}[u]^{(k)}_j\qquad& \mbox{otherwise}\cr}\right.\nn\\
&&x^+_i(z)[u]^{(k)}_j=\left\{\mmatrix{C_+\delta(t_1^{-j}u/t_1z)e^{-Q_i}[u]^{(k)}_{j+1},\qquad& i+j+1\equiv k\cr
0,\qquad& i+j+1\not\equiv k\cr}\right.\lb{vecrep}\\
&&x^-_i(z)[u]^{(k)}_j=\left\{\mmatrix{C_-\delta(t_1^{-j}u/z)[u]^{(k)}_{j-1},\qquad& i+j\equiv k\cr
0,\qquad& i+j\not\equiv k\cr}\right.,\nn
\ena
where 
\bea
&&C_\pm=\frac{(p\hbar^{\pm 1};p)_\infty}{(p;p)_\infty}.
\lb{def:Cpm}
\ena

\begin{lem}
With the action in \eqref{vecrep}, $V^{(k)}(u)$ is an irreducible, tame, $\Z^N$-graded $\cU$-module of 
level (0,-1).    
\end{lem}

Let $\la=(\la_1\geq \la_2\geq \cdots)$ be a partition, where $\la_i\geq 0\ (i=1,2,\cdots)$, and   $\ell(\la)<+\infty$ be its length. 
For $j\in \Z_{\geq 1}$, let $\la\pm \bold{1}_j=(\la_1,\la_2,\cdots,\la_j\pm1,\cdots)$. 
 We identify $\la$ with a Young diagram $\la=\{(x,y)\ |\ 1\leq x\leq \ell(\la), \ 1\leq y\leq \la_{x}\ \}$. 
 \begin{dfn}\lb{def:colork}
 We call a partition $\la$ equipped with content $c_{(x,y)}=x-y+k$ modulo $N$ for a box $(x,y)\in \la$ a partition with color $k$.  
We denote by $\cP^{(k)}$ a set of all partitions with color $k$. 
\end{dfn}

For $\la\in \cP^{(k)}$, consider  the subspace ${\F}^{(0,-1)(k)}_u$ 
\be
&&{\F}^{(0,-1)(k)}_u\subset V^{(k)}(u)\tot V^{(k)}(u\hbar^{-1})\tot V^{(k)}(u\hbar^{-2})\tot \cdots
\en
 spanned by the vectors 
\be
\ket{\la}^{(k)}_u=[u]^{(k)}_{\la_1-1}\tot [u\hbar^{-1}]^{(k)}_{\la_2-2}\tot [u\hbar^{-2}]^{(k)}_{\la_3-3}\tot \cdots.
\en 
For $j\in I $, let $A^{(k)}_j(\la)$ and $R^{(k)}_j(\la)$ denote the set of all addable and removable boxes in $\la$  
of content $j$, respectively. For $(x,y),  (x',y')\in \la$, let  $(x,y) > (x',y')$ if $c_{(x,y)}> c_{(x',y')}$.  
From Corollary 6.3 in \cite{KO24}, we have the following statement. 
\begin{thm}\lb{lebel0repbox}
For $j\in {I}$, the following gives a level-$(0,-1)$ irreducible lowest weight representation of $\cU$ on 
$\F^{(0,-1)(k)}_u
$.
\bea
&&\phi_j^\pm(z)\ket{\la}^{(k)}_u=\prod_{R\in R^{(k)}_j(\la)} \left. \frac{\theta(u_{R}/z)}{\theta(\hbar u_R/z)}\right|_\pm
\prod_{A\in A^{(k)}_j(\la)} \left. \frac{\theta(\hbar^2 u_A/z)}{\theta(\hbar u_A/z)}\right|_\pm
e^{-Q_{j}}\ket{\la}^{(k)}_u,\lb{phijpm0box}\\
&&x^+_j(z)\ket{\la}^{(k)}_u=C_+\sum_{X\in A^{(k)}_j(\la)}\delta(\hbar u_X/z)A^+_{\la,X}e^{-Q_{j}}
\ket{\la\cup X}^{(k)}_u,\lb{ej0box}\\
&&x^-_j(z)\ket{\la}^{(k)}_u=C_-\sum_{X\in R^{(k)}_j(\la)}\delta(\hbar u_X/z)A^-_{\la,X}
\ket{\la\backslash X}^{(k)}_u,  \lb{fj0box}
\ena
where  $u_{X}=t_1^{-y}t_2^{-x} u$  for $X=(x,y) \in \la$, and  
\bea
A^+_{\la,X}&=&\prod_{R\in R^{(k)}_j(\la) \atop R<X} \frac{\theta(\hbar u_X/u_R)}{\theta(u_X/u_R)}
\prod_{A\in A^{(k)}_j(\la) \atop A<X } \frac{\theta(\hbar^{-1}u_X/u_A)}{\theta(u_X/u_A)},\lb{def:Apbox}\\
&=&\prod_{R\in R^{(k)}_j(\la) \atop R<X} \frac{\theta(\hbar u_X/u_R)}{\theta(u_X/u_R)}
\prod_{A\in A^{(k)}_j(\la\cup X) \atop A<X } \frac{\theta(\hbar^{-1}u_X/u_A)}{\theta(u_X/u_A)},\lb{Apbox}
\ena
\bea
A^-_{\la,X}&=&\prod_{R\in R^{(k)}_j(\la) \atop R> X}  \frac{\theta(\hbar^{-1}u_R/u_X)}{\theta(u_R/u_X)}
\prod_{A\in A^{(k)}_j(\la) \atop A>X} 
\frac{\theta(\hbar u_A/u_X)}{\theta(u_A/u_X)}.\lb{def:Ambox}\\
&=&\prod_{R\in R^{(k)}_j(\la\backslash X) \atop R> X}  \frac{\theta(\hbar^{-1}u_R/u_X)}{\theta(u_R/u_X)}
\prod_{A\in A^{(k)}_j(\la) \atop A>X} 
\frac{\theta(\hbar u_A/u_X)}{\theta(u_A/u_X)}.\lb{def:Ambox}
\ena
In particular,
\bea
&&K\cdot
\ket{\la}^{(k)}_u=\hbar^{(|R(\la)|-|A(\la)|)/2}\ket{\la}^{(k)}_u=\hbar^{-1/2}\ket{\la}^{(k)}_u,\lb{kappa0}
\ena
where we set $|R(\la)|=\sum_{j\in I}|R^{(k)}_j(\la)|$ and $|A(\la)|=\sum_{j\in I}|A^{(k)}_j(\la)|$. In \eqref{phijpm0box}, $|_\pm$ denotes the expansion direction 
\bea
\left. \frac{\theta(z/t_2)(p;p)_\infty^3}{\theta(t_1z)\theta(\hbar^{-1})}\right|_+&=&\sum_{n\in \Z}\frac{1}{1-\hbar^{-1}p^n}(z/t_2)^n,
\ena
for $|p|<| z/t_2|<1$, and 
\bea
\left. \frac{\theta(z/t_2)(p;p)_\infty^3}{\theta(t_1z)\theta(\hbar^{-1})}\right|_-
&=&-\frac{\theta(t_2/z)(p;p)_\infty^3}{\theta(1/t_1z)\theta(\hbar)}\nn\\
&=&\sum_{n\in \Z}\frac{1}{1-\hbar p^n}(t_2/z)^n,
\ena
for $1<|z/t_2|<|p^{-1}|$. 
In \eqref{ej0box} and \eqref{fj0box}, $\ket{\la\pm\bold{1}_i}^{(k)}_u$ vanishes if ${\la\pm \bold{1}_i}$ is  not a partition, respectively. The lowest weight vector is $\ket{\emptyset}^{(k)}_u$ and the lowest weight is 
\be
&&\Phi(z)=\left(\left(\frac{\theta(\hbar u/z)}{\theta(u/z)}\right)^{\delta_{j,k}}\right)_{j\in I}.
\en
\end{thm}
We call this representation the $q$-Fock representation ${\F}^{(0,-1)(k)}_u$ with color $k$ and vacuum weight $u$. This is an elliptic and dynamical analogue of the one in \cite{FJMM2,STU,VV98}. 

From Theorem \ref{lebel0repbox}, one obtains the following.
\begin{cor}\lb{prop:wt}
Setting $K^\pm_i=\hbar^{\pm h_i/2}e^{-Q_i}$, 
one has
\be
&&\hbar^{h_i/2}\ket{\la}^{(k)}_u=\hbar^{\wt(\la)^{(k)}_i/2}\ket{\la}^{(k)}_u
\en
with
\be
&&\wt(\la)^{(k)}_i=|R^{(k)}_i(\la)|-|A^{(k)}_i(\la)|.
\en
Here  $A^{(k)}_i(\la)$ and $R^{(k)}_i(\la)$ denote the set of all addable and removable boxes in $\la\in \cP^{(k)}$ with content $i$, respectively. 
\end{cor}
\noindent
{\it Proof.} \ 
From \eqref{def:theta}, one has in \eqref{phijpm0box}
\be
&&\frac{\theta(u_R/z)}{\theta(\hbar u_R/z)}=\hbar^{1/2}\frac{\theta_p(u_R/z)}{\theta_p(\hbar u_R/z)},\qquad
\frac{\theta(\hbar^2u_A/z)}{\theta(\hbar u_A/z)}=\hbar^{-1/2}\frac{\theta_p(u_R/z)}{\theta_p(\hbar u_R/z)}.
\en
\qed

For $k\in I$, let  ${\cal M}({\bold{v},\bold{w}})$ be the  {quiver variety  of type $A^{(1)}_{N-1}$} with 
 $\bold{v}=(v_0,\cdots,v_{N-1})$, $v_i=\#\{\square \in \la | c_{\square}\equiv i\ \mbox{mod}\ N \}$ and $\bold{w}=(w_0,w_1,\cdots,w_{N-1})$, $w_i=\delta_{i,k}$.  See Sec.\ref{sec:affAquivervar}. 
The trigonometric limit, $p\to 0$, of the case $k=0$ of ${\F}^{(0,-1)(k)}_u$ 
 in Theorem \ref{lebel0repbox} coincides with the representation of the quantum toroidal algebra $U_{t_1,t_2}(\gl_{N,tor})$  
on the equivariant 
$\mathrm{K}$-theory 
$\bigoplus_{\bold{v}}{\mathrm K}_\rT({\cal M}({\bold{v},\bold{w}}))$ 
with $\rT=  (\C^\times)^{|\bold{w}|=1}_u \times \C^\times_{t_1} \times\C^\times_{t_2}$  
studied in  \cite{VV99, Na01} except for a dependence on the dynamical parameters.  
 In particular, the bases $\ket{\la}^{(k)}_u$  can be identified with the fixed point classes labeled by the partition $\la\in \cP^{(k)}$ \cite{Nagao,Negut}. We hence reach the following conjecture. 

\begin{conj} \cite{KO24}\lb{Nvw}
The representation $\F^{(0,-1)(k)}_u$ of $U_{t_1,t_2,p}(\gl_{N,tor})$ in Theorem \ref{lebel0repbox} is equivalent to an expected  level-$(0,-1)$  geometric action of  the same algebra on  $\bigoplus_{\bold{v}}{\mathrm E}_\rT({\cal M}({\bold{v},(\delta_{i,k})}))$ under the identification of $\ket{\la}^{(k)}_u$ with the $\rT$-fixed point class $[\la]$ of ${\cal M}(\bold{v},(\delta_{i,k}))$.   
\end{conj}

This conjecture is one of the motivation of our construction of the vertex operator of $U_{t_1,t_2,p}(\gl_{N,tor})$ in Sec.\ref{sec:VO}. See also Conjecture \ref{level0mw2Mvw}.

\section{Elliptic Stable Envelopes  }\lb{sec:ESE}

We summarize some basic facts and formulas for the elliptic stable envelopes for the equivariant elliptic cohomology of the $A^{(1)}_{N-1}$ quiver variety.

\subsection{The $A^{(1)}_{N-1}$ quiver variety $\cM(\bfv,\bfw)$}\lb{sec:affAquivervar}

\subsubsection{Definition}
Consider the affine Dynkin quiver of type $A^{(1)}_{N-1}$ consists of $N$-vertices labeled by $i\in I=\{0,1,\cdots,N-1\}$ and edges $i\to i+1$ with $N\equiv 0$.  Let $V_i , W_i$ be vector spaces of dimensions $v_1, w_i$ $(i\in I)$, respectively. We set $\bfv=(v_0,v_1,\cdots,v_{N-1})$ and  $\bfw=(w_0,w_1,\cdots,w_{N-1})$. 
Consider 
\be
&&M=\bigoplus_{i=0}^{N-1}\Hom(V_i,V_{i+1})\oplus \bigoplus_{i=0}^{N-1}\Hom(W_i,V_{i}).
\en
The cotangent bundle $T^*M$ has a natural $G_{v}=\prod_{i=0}^{N-1}GL(V_i)$-action which induces a moment map
\be
&&\mu\ :\ T^*M\ \to\ (\Lie G_v)^*.
\en
Let $\theta$ be the cocharacter of $G_v$ given by
\be
&&\theta\ :\ G_v\ \to\ \C^\times,\qquad g=(g_i)_{i\in I}\ \mapsto\ \prod_{i\in I}\det g_i.
\en
\begin{dfn}
The $A^{(1)}_{N-1}$ quiver variety $\cM(\bfv,\bfw)$ is the following symplectic reduction.
\be
&&\cM(\bfv,\bfw)=\mu^{-1}(0)/\hspace{-1mm}/\hspace{-0.8mm}{}_\theta\, G_v=\mu^{-1}(0)^{\theta-ss}/G_v.
\en
\end{dfn} 
Through this paper we fix a positive stability condition $\theta$. 

\subsubsection{Torus fixed points}\lb{sec:TFP}
Let $|\bfv|=\sum_{i=0}^{N-1}v_i$, $|\bfw|=\sum_{i=0}^{N-1}w_i$.  We set  $\rA_{\bfw}=(\C^\times)^{|\bfw|}$ and  $\rT=\rA_{\bfw}\times \C^\times_{t_1}\times \C^\times_{t_2}$ be algebraic tori. There is a natural $\rT$ action on $\cM(\bfv,\bfw)$.  Here we made an identification of $t_1, t_2$ in the elliptic quantum toroidal algebra $U_{t_1,t_2,p}(\glnt)$ with the coordinates of of the torus $\C^\times_{t_1}\times \C^\times_{t_2}$. From \eqref{hbar}, we have
\bea
&&t_1=\kappa^{-1} \hbar^{1/2}, \qquad t_2=\kappa^{-1}\hbar^{1/2}.
\ena
We hence identify $\hbar^{-1}$ and $\kappa^{-1}$ with the $\C^\times_{t_1}\times \C^\times_{t_2}$-weight of the symplectic form on $\cM(\bfv,\bfw)$ and  the coordinate on $\ker(\hbar^{-1}) \subset \C^\times_{t_1}\times \C^\times_{t_2}$, respectively.  The subtorus $\rA=\rA_\bfw\times \C^\times_{\kappa^{-1}}$ preserves the symplectic form. 

For $k\in I
$, let $\la^{(k)}_j=(\la^{(k)}_{j,1}\geq \la^{(k)}_{j,2}\geq \cdots)$ $(j=1,\cdots,w_k) 
$ be partitions with color $k$, 
and consider $|\bfw|$-tuple of partitions $\blla=(\la^{(0)}_1,\cdots,\la^{(0)}_{w_0},\cdots,\la^{(N-1)}_1,\cdots,\la^{(N-1)}_{w_{N-1}})$. 
For $\square=(x,y)\in \la^{(k)}_j$, 
the content of $\square$ is given by  $c_\square=x-y+k$ as before. 
\begin{prop}
There are finitely many $\rT$-fixed points on $\cM(\bfv,\bfw)$ labeled by 
$|\bfw|$-tuple of partitions $\blla$ satisfying
\bea
&&v_i=\#\{\ \square\in \blla\ |\ c_\square\equiv i\ \mbox{mod}\ N \}\qquad i\in I. \lb{fixedpt}
\ena
Hence $|\blla|=\sum_{k=0}^{N-1}\sum_{j=1}^{w_k}|\la^{(k)}_j|=|\bfv|$. 
 \end{prop}
 Let $\Lambda_{\bfv,\bfw}$ be the set of all $|\bfw|$-tuple of partitions $\blla$ satisfying \eqref{fixedpt}. 
We identify $\cM(\bfv,\bfw)^\rT$ with $\Lambda_{\bfv,\bfw}$.

We introduce an ordering on the boxes of $\blla\in \Lambda_{\bfv,\bfw}$ as follows.
For a box $\square=(x,y)\in \la^{(k)}_j$, let us set 
\be
&&h_\square:=x+y-2,\\
&&\trho_\square:=c_\square-\ep h_\square
\en
with $0<\ep\ll1$. We often denote boxes in $\blla$ by $a, b, c,  \cdots$, and  define a  canonical ordering on them
 by 
\be
&&a\prec b\ \Leftrightarrow \ \trho_a<\trho_b \ \  \mbox{for}\  \ a,b \in \la^{(k)}_i,\\
&&a\prec b\ \ \mbox{for}\  a\in  \la^{(k)}_i,\ b\in  \la^{(k)}_j,\quad i<j,\\
&&a\prec b\ \ \mbox{for}\  a\in  \la^{(k)}_i,\ b\in  \la^{(l)}_j,\quad k<l,\ \forall i, j.
\en 
Set also $I_\blla^{(i)}=\{\ a\in \blla\ |\ c_a\equiv i\ \}, i\in I$ so that $|I_\blla^{(i)}|=v_i$ and $\blla=\cup_{i\in I}I_\blla^{(i)}$. 
Let $a^{(i)}_j$ denote the $j$-th box in $I_\blla^{(i)}$ in the order $\prec$, i.e. 
\be
&&  I_\blla^{(i)}=\{\ a^{(i)}_1\prec a^{(i)}_2\prec \cdots \prec a^{(i)}_{v_i} \ \}.
\en

\medskip
\noindent
{\it Example 1.}\ $N=3$, $\la=(6,5,4,1)$ with the vacuum color $0$. Then $\bfv=(v_0=6, v_1=5, v_2=5)$ and 
$\bfw=(w_0=1,w_1=0,w_2=0)$. 
\be
\la&=&I^{(0)}_\la\cup I^{(1)}_\la\cup I^{(2)}_\la,\\
I^{(0)}_\la&=&\{(2,5)\prec(1,4) \prec (3,3)\prec(2,2)\prec (1,1)\prec (4,1)\},\\ 
 I^{(1)}_\la&=&\{(1,6)\prec (2,4)\prec(1,3)\prec (3,2)\prec(2,1)\},\\
 I^{(2)}_\la&=&  \{(1,5)\prec(3,4)\prec(2,3)\prec(1,2)\prec(3,1)\}.
\en
\qed 

\subsubsection{Tautological vector bundles}\lb{tautVB}
Let  $\cV_i$ be the tautological vector bundles of rank $v_i$ $(i\in I)$ on $\cM(\bfv,\bfw)$ given by
\be
&&\cV_i=\mu^{-1}(0)^{\theta-ss}\times_{G_v} V_i. 
\en
The bundles $\cV_i,\ i\in I$ have a natural $\rT$-equivariant structure.
We denote the corresponding tautological $\rK$-theory class by the same symbol $\cV_i$, which is expressed by the Chern roots  $x^{(i)}_j$ as 
\bea
&&\cV_i=\sum_{j=1}^{v_i}x^{(i)}_j.\lb{cVi}
\ena 
We have an obvious  bijection between $I^{(i)}_\blla=\{a^{(i)}_j\}$ and the set of  $\{x^{(i)}_j\}$, $j=1,\cdots,v_i$ and define $x_{a^{(i)}_j}:=x^{(i)}_{j}$, $j=1,\cdots, v_i $.

\medskip
\noindent
{\it Example 2.}\  Let $\la$ be the same as in Example 1. Then one has
\be
\cV_0&:&  x^{(0)}_1=x_{(2,5)},\ x^{(0)}_2=x_{(1,4)},\ x^{(0)}_3=x_{(3,3)},\ x^{(0)}_4=x_{(2,2)},\  x^{(0)}_5=x_{(1,1)},\ x^{(0)}_6=x_{(4,1)},\\
\cV_1&:&  x^{(1)}_1=x_{(1,6)},\ x^{(1)}_2=x_{(2,4)},\ x^{(1)}_3=x_{(1,3)},\ x^{(1)}_4=x_{(3,2)},\  x^{(1)}_5=x_{(2,1)},\\
\cV_2&:&  x^{(2)}_1=x_{(1,5)},\ x^{(2)}_2=x_{(3,4)},\ x^{(2)}_3=x_{(2,3)},\ x^{(2)}_4=x_{(1,2)},\  x^{(2)}_5=x_{(3,1)}.
\en 
\qed

We also have the tautological vector bundles $\cW_i$ of rank $w_i$, $i\in I$ by using the vector spaces $W_i$. 
Using  the framing weights  $u^{(i)}_j$,\ $j=1,\cdots, w_i$ at the $i$-th vertex, i.e. the coordinates on the framing torus $(\C^\times)^{w_i}\subset \rA_{\bfw}$, we have 
\bea
&&\cW_i=\sum_{j=1}^{w_i}u^{(i)}_j.\lb{cWi}
\ena
For a box $a=(x,y)\in \la^{(k)}_j$ with $c_a\equiv i$, the restriction of $\cV_i$ on the fixed point $\blla$ is given by a 
replacement of the Chern root $x_a$ with the $\rT$-weight
\be
&&\varphi^{\blla}_a:=t_1^{-y+1}t_2^{-x+1}.
\en
Namely,
\be
&&\cV_i\Bigr|_{\blla}=\sum_{j=1}^{v_i}\varphi^{\blla}_{a^{(i)}_j}=\sum_{a\in \blla\atop c_a\equiv i}\varphi^{\blla}_a.
\en

\subsection{Elliptic stable envelopes for $\E_\rT(\cM(\bfv,\bfw))$}
Let $E=\C^\times/p^\Z$ be an elliptic curve with $0<|p|<1$  and consider the equivariant elliptic cohomology $\E_\rT(X)$ for $X=\cM(\bfv,\bfw)$\cite{GKV,AO}.  The elliptic stable envelope $\Stab_{\gC,T^{1/2}}(\blla)$ for  $\E_T(X)$ was constructed in \cite{Dinkins21} following the definition given in \cite{AO} and the preceding work \cite{Smirnov18} for the Hilbert scheme on points in $\C^2$. They are certain sections of line bundle over $\E_\rT(X)$ determined depending on a choice of fixed point $\blla\in X^\rT$, polarization $T^{1/2}\in \rK_\rT(X)$  satisfying 
\be
&&
TX=T^{1/2}+\hbar^{-1}(T^{1/2})^*, 
\en
and chamber $\gC\subset {\rm Lie}_\R(\rA)\backslash\{w^{\perp}\}
$, where $\{w\}$ be a set of $\rA$-weights appearing in the normal bundle to $X^\rA$ in $X$.
We chose the chamber $\gC$ as
\bea
&&\hspace{-1cm}|u^{(0)}_1|> 
\cdots> |u^{(0)}_{w_0}|> |u^{(1)}_1|> 
\cdots> |u^{(1)}_{w_1}|> \cdots\cdots > |u^{(N-1)}_{1}|> 
\cdots > |u^{(N-1)}_{w_{N-1}}|\lb{chamber}
\ena
and $|t_1/t_2|=|\kappa^{-1}|<1$. Note that the latter choice corresponds to the choice $t_+=t_1, t_-=t_2$ in \cite{Dinkins21}.

For $\blla\in \Lambda_{\bfv,\bfw}$, denote by $\bfx=\{ x_a | a\in \blla\ \}=\cup_{i\in I}\{ x^{(i)}_j, j=1,\cdots,v_i  \}$ the set of all Chern roots and by $\bfu=\{ u^{(k)}_j\ |\ k\in I, j=1,\cdots,w_k \}$ the set of all framing weights. Define
\bea
&&S_{\blla}(\bfx, \bfu; p, t_1, t_2)=\frac{S_1(\bfx; p, t_1, t_2)S_2(\bfx,\bfu; p, t_1, t_2)}{S_3(\bfx; p, t_1, t_2)}\lb{Sla}
\ena
with
\bea
&&S_1(\bfx; p, t_1, t_2)=
\prod_{a,b\in\blla \atop {c_a+1\equiv c_b\atop 
\rho_a+1<\rho_b}}\theta(t_1x_a/x_b)\prod_{a,b\in\blla \atop {c_a+1\equiv c_b\atop 
\rho_a+1>\rho_b}}\theta(t_2x_b/x_a),\lb{S1}\\
&&S_2(\bfx,\bfu; p, t_1, t_2)=\prod_{k=0}^{N-1}\prod_{j=1}^{w_k}\left(
\prod_{a\in\blla \atop {c_a\equiv k\atop 
\rho_a\leq \rho_{r_{k,j}}}}\theta(x_a/u^{(k)}_j)
\prod_{a\in\blla \atop {c_a\equiv k\atop 
\rho_a> \rho_{r_{k,j}}}}\theta(t_1t_2u^{(k)}_j/x_a)
\right),\lb{S2}\\
&&S_3(\bfx; p, t_1, t_2)=\prod_{a,b\in\blla \atop {c_a\equiv c_b\atop 
\rho_a<\rho_b}}\theta(x_a/x_b)\theta(t_1t_2x_a/x_b).\lb{S3}
\ena
Here $r_{k,j}$ denotes a box at the corner $(1,1)$ of a partition  $\la^{(k)}_j$.  
 Our $S_\blla$
 corresponds to the choice $t_+=t_1, t_-=t_2$, $t(i,i+1)=i, h(i,i+1)=i+1$, $i\in I$
   in \cite{Dinkins21}. 

Let $t^{(i)}_j$ be a $\la$-tree in a partition $\la^{(i)}_j$, see Definition 1 in Sec.4.2 of \cite{Smirnov18} and Sec.2.3, 2.4 of \cite{Dinkins21}. Namely, $t^{(i)}_j$ is a rooted tree in a Young diagram $\la^{(i)}_j$ with vertices corresponding to boxes of $\la^{(i)}_j$, edges connecting only adjacent boxes and the root at the box $(1,1) \in\la^{(i)}_j$. Let us set  (Definition 12 in \cite{Dinkins21}):
\bea
&&W^{\la^{(i)}_j}_{t^{(i)}_j}(\bfx,\bfu,\bfz; p,t_1,t_2)\nn\\
&&=(-1)^{\kappa(t^{(i)}_j)}\phi\left(\frac{x_{r_{i,j}}}{u^{(i)}_j
},\prod_{a\in [r_{i,j},t^{(i)}_j]}z_{c_a}(t_1t_2)^{d^\blla_a}\right)
\prod_{e\in t^{(i)}_j}
\phi\left(\frac{x_{h(e)}\varphi^\blla_{t(e)}}{x_{t(e)}\varphi^\blla_{h(e)}},\prod_{a\in [h(e),t^{(i)}_j]}z_{c_a}(t_1t_2)^{d^\blla_a}\right),\nn\\
&&\lb{Wlat}
\ena
where $\bfz=\{z_i\}_{ i\in I}$ denotes a set of the K\"ahler parameters  and 
\bea
&&\phi(x,y)=\frac{\theta(xy)\theta(\hbar)}{\theta(x)\theta(y)}.
\ena
Here  the product runs over the edges of the tree $t^{(i)}_j$ and $h(e) \in \lambda^{(i)}_j$, $t(e)\in \lambda^{(i)}_j$ denote the head and tail box of the edge $e$. 
For a box  $a\in \la^{(i)}_j$, the symbol $[a,t^{(i)}_j]$ denotes the set of boxes in $\la^{(i)}_j$ appearing as vertices in the subtree of $t^{(i)}_j$ rooted at $a$. 
The integer $\kappa({{t^{(i)}_j}})$ denotes 
\be
&&\#\mbox{(arrows in $t^{(i)}_j$ directed decreasing $y$)}+\#\mbox{(arrows in $t^{(i)}_j$ directed decreasing $x$)} 
\en
for  $\la^{(i)}_j=\{(x,y)\ |\ 1\leq x\leq \ell(\la^{(i)}_j), 1\leq y\leq \la^{(i)}_{j,x} \}$. 
The integer $d^\blla_a$ is defined by 
\be
&&\det {\rm ind}_\blla=\prod_{a\in \blla} x_a^{d^{\blla}_a}
\en
with the index class ${\rm ind}_\blla=T^{1/2}|_{X^\rA,>0}\in \rK_\rT(X^\rA)$. 
We refer to Sec. 2.8, 4.5 in \cite{Smirnov18} and Sec.4.1 in \cite{Dinkins21} for definitions of these integers. 
Let 
\be
&&\Gamma_\blla=\prod_{i=0}^{N-1}\prod_{j=1}^{v_i}\Gamma_{\la^{(i)}_j}
\en
with $\Gamma_{\la^{(i)}_j}
$ being the set of trees  in $\la^{(i)}_j$ without \reflectbox{$\mathsf{L}$}-shaped subgraphs, see section 4.6 in \cite{Smirnov18}. We set $\bft=(t^{(i)}_j)$, $i\in I, j\in \{1,\cdots,w_i\}$ and 
\be
&&W^{\bla}_{\bft}(\bfx,\bfu,\bfz;p,t_1,t_2)=\prod_{i\in I}\prod_{j=1}^{w_i}
W^{\la^{(i)}_j}_{t^{(i)}_j}(\bfx,\bfu,\bfz; p,t_1,t_2).
\en

\begin{thm}\lb{EllStab}\cite{Dinkins21}
The elliptic stable envelope of $\blla\in X^\rT$ is given by
\be
&&\Stab_{\gC,T^{1/2}}(\blla;\bfz)
=\Sym_0\Sym_1\cdots\Sym_{N-1}\left(S_{\blla}(\bfx, \bfu; p, t_1, t_2)\sum_{\bft\in \Gamma_\blla}W^{\blla}_{\bft}(\bfx,\bfu,\bfz;p,t_1,t_2)\right),
\en
where $\Sym_i$ denotes the symmetrization over the Chern roots $\{x^{(i)}_j,\ j=1,\cdots,v_i \}=\{ x_a,\ a\in \blla,  c_a\equiv i\ \mbox{mod}\ N\}$ of $\cV_i$.  
\end{thm}
Later we identify the dynamical parameter $\gz_i$ in $U_{t_1,t_2,p}(\gl_{N,tor})$ with the K\"ahler parameter in $\E_\rT(X)$\cite{AO} in the same way as in \cite{Konno18}. See Sec.\ref{sec:shuffleStab}.

 In the next section, we make a connection to the representations of $U_{t_1,t_2,p}(\gl_{N,tor})$. 
 For this purpose, we introduce $\widehat{\Stab}$ and $\widetilde{\Stab}$ defined as follows. 
Let us set
\bea
K_I(\bfx,\bfu)&=&\frac{
\mbox{
$\ds{\prod_{a,b\in\blla \atop {c_a+1\equiv c_b\atop \rho_a+1<\rho_b}}{\theta(t_2x_b/x_a)}
\prod_{a,b\in\blla \atop {c_a+1\equiv c_b\atop \rho_a+1>\rho_b}}{\theta(t_1x_a/x_b)}
\prod_{i=0}^{N-1}\prod_{j=1}^{w_i}\prod_{a\in\blla \atop c_a\equiv i}{\theta(x_a/u^{(i)}_j)}}$
}
}{
\mbox{
$\ds{\prod_{a,b\in\blla \atop {c_a\equiv c_b\atop \rho_a\not=\rho_b}}\theta(x_a/x_b)}$
}
}.\lb{Thompol}
\ena
Define $\widehat{S}_\blla(\bfx,\bfu;p,t_1,t_2)$
 by
\bea
{S}_\blla(\bfx,\bfu;p,t_1,t_2)&=&(-1)^{\vep(\bla)}K_I(\bfx,\bfu)
\widehat{S}_\blla(\bfx,\bfu;p,t_1,t_2),\lb{def:hatSla}
\ena
where
\bea
&&\vep(\bla)=\sum_{a,b\in \bla\atop c_a+1\equiv c_b}1+\sum_{i\in I}\sum_{j=1}^{w_i}\sum_{a\in \bla\atop{c_a\equiv i \atop \rho_a>\rho_{r_{i,j}}}}1.
\ena
One finds
\bea
&&\hspace{-1.5cm}\widehat{S}_\blla(\bfx,\bfu;p,t_1,t_2)=
\prod_{a,b\in\blla \atop {c_a+1\equiv c_b\atop \rho_a+1<\rho_b}}\left(-\frac{\theta(t_1x_a/x_b)}{\theta(t_2x_b/x_a)}\right)
\prod_{a,b\in\blla \atop {c_a+1\equiv c_b\atop \rho_a+1>\rho_b}}\left(-\frac{\theta(t_2x_b/x_a)}{\theta(t_1x_a/x_b)}\right)\nn\\
&&\hspace{2cm}\times
\prod_{i=0}^{N-1}\prod_{j=1}^{w_i}\hspace{-0.3cm}\prod_{a\in\blla \atop {c_a\equiv i \atop \rho_a>\rho_{r_{i,j}} }}
\left(-\frac{\theta(\hbar u^{(i)}_j/x_a)}{\theta(x_a/u^{(i)}_j)}\right)
\prod_{a,b\in\blla \atop {c_a\equiv c_b\atop \rho_a<\rho_b}}\frac{\theta(x_b/x_a)}{\theta(\hbar x_a/x_b)}.\lb{hatSla}
\ena
We then define the hatted $\Stab$ as follows.
\bea
&&\hspace{-1cm}\widehat{\Stab}_{\gC,T^{1/2}}(\blla;\bfz)\nn\\
&&\hspace{-1cm}:=\Sym_0\Sym_1\cdots\Sym_{N-1}\left(\widehat{S}_{\blla}(\bfx, \bfu; p, t_1, t_2)\sum_{\bft\in \Gamma_\blla}W^{\blla}_{\bft}(\bfx,\bfu,\bfz;p,t_1,t_2)\right).
\lb{def:hatStab}
\ena

Similarly, let us set
\bea
K_{II}(\bfx,\bfu)&=&\frac{
\mbox{
$\ds{\prod_{a,b\in\blla \atop {c_a+1\equiv c_b\atop \rho_a+1<\rho_b}}{\theta(t_2x_b/x_a)}
\prod_{a,b\in\blla \atop {c_a+1\equiv c_b\atop \rho_a+1>\rho_b}}{\theta(t_1x_a/x_b)}
\prod_{i=0}^{N-1}\prod_{j=1}^{w_i}\prod_{a\in\blla \atop c_a\equiv i}{\theta(\hbar u^{(i)}_j/x_a)}}$
}
}{
\mbox{
$\ds{\prod_{a,b\in\blla \atop {c_a\equiv c_b\atop \rho_a\not=\rho_b}}\theta(\hbar x_a/x_b)}$
}
}\lb{ThompolTypeII}
\ena
and define $\widetilde{S}_\blla(\bfx,\bfu;p,t_1,t_2)$
 by
\bea
{S}_\blla(\bfx,\bfu;p,t_1,t_2)&=&(-1)^{\vep^*(\bla)}K_{II}(\bfx,\bfu)
\widetilde{S}_\blla(\bfx,\bfu;p,t_1,t_2),\lb{def:hatSla}
\ena
where
\bea
&&\vep^*(\bla)=\sum_{a,b\in \bla\atop c_a+1\equiv c_b}1+\sum_{i\in I}\sum_{j=1}^{w_i}\sum_{a\in \bla\atop{c_a\equiv i \atop \rho_a\leq\rho_{r_{i,j}}}}1.
\ena
One finds
\bea
&&\hspace{-1.5cm}\widetilde{S}_\blla(\bfx,\bfu;p,t_1,t_2)=
\prod_{a,b\in\blla \atop {c_a+1\equiv c_b\atop \rho_a+1<\rho_b}}\left(-\frac{\theta(t_1x_a/x_b)}{\theta(t_2x_b/x_a)}\right)
\prod_{a,b\in\blla \atop {c_a+1\equiv c_b\atop \rho_a+1>\rho_b}}\left(-\frac{\theta(t_2x_b/x_a)}{\theta(t_1x_a/x_b)}\right)\nn\\
&&\hspace{2cm}\times
\prod_{i=0}^{N-1}\prod_{j=1}^{w_i}\hspace{-0.3cm}\prod_{a\in\blla \atop {c_a\equiv i \atop \rho_a\leq\rho_{r_{i,j}} }}
\left(-\frac{\theta(x_a/ u^{(i)}_j)}{\theta(\hbar u^{(i)}_j/x_a)}\right)
\prod_{a,b\in\blla \atop {c_a\equiv c_b\atop \rho_a<\rho_b}}
\frac{\theta(\hbar x_b/x_a)}{\theta(x_a/x_b)}.\lb{hatSla}
\ena
We then define $\widetilde{\Stab}$ as follows.
\bea
&&\hspace{-1cm}\widetilde{\Stab}_{\gC,T^{1/2}}(\blla;\bfz)\nn\\
&&\hspace{-1cm}:=\Sym_0\Sym_1\cdots\Sym_{N-1}\left(\widetilde{S}_{\blla}(\bfx, \bfu; p, t_1, t_2)\sum_{\bft\in \Gamma_\blla}W^{\blla}_{\bft}(\bfx,\bfu,\bfz;p,t_1,t_2)\right).
\lb{def:hatStabTypeII}
\ena

\subsection{Shuffle product of  $\Stab_{\gC,T^{1/2}}(\bla;\bfz)$}\lb{sec:shuffleStab}

The  shuffle product formula for $\Stab_{\gC,T^{1/2}}(\bla;\bfz)$ is obtained in a similar way to the Jordan quiver case\cite{Botta,KS}.  
Let  $\cM(\bfv',\bfw')$ and $\cM(\bfv'',\bfw'')$ be the $A^{(1)}_{N-1}$ quiver varieties. 
Let $\Stab_{\gC',T^{1/2'}}(\bla';\bfz')$ be 
 the elliptic stable envelopes  for $\E_\rT(\cM(\bfv',\bfw'))$  with $\bla'\in \Lambda_{\bfv',\bfw'}$, the Chern roots $\bfx'=\{x_a'\}_{a\in \bla'}$ and  the framing weights $\bu'=\{u^{'(i)}_j |\  i\in I, j=1,\cdots,w_i'\}$. Similarly let $\Stab_{\gC'',T^{1/2''}}(\bla'';\bfz'')$  for  $\E_\rT(\cM(\bfv'',\bfw''))$ with $\bla''\in \Lambda_{\bfv'',\bfw''}$,   the Chern roots $\bfx''=\{x_a''\}_{a\in \bla''}$ and
  the framing weights  $\bu''=\{u^{''(i)}_j |\  i\in I, j=1,\cdots,w_i''\}$. The  chambers 
  $\gC'$ and $\gC''$ are taken as \eqref{chamber}, respectively,  with a common stability condition $|t_1/t_2|<1$.
The polarizations $T^{1/2'}$ and $T^{1/2''}$ are taken in such a way that $\Stab_{\gC',T^{1/2'}}(\bla';\bfz')$ and $\Stab_{\gC'',T^{1/2''}}(\bla'';\bfz'')$ 
are given by the formula in Theorem \ref{EllStab} for $\bla'$ and $\bla''$, respectively.  
We set $\bfx=\bfx'\cup \bfx''$ and $\bu=\bu'\cup\bu''$. 

\begin{thm}\lb{shuffleStabgen}
The following shuffle product formula gives  the elliptic stable envelope $\Stab_{\gC,T^{1/2}}((\bla',\bla'');\bfz)$ for  $\E_T(\cM(\bfv,\bfw))$  with $\bfv=\bfv'+\bfv''$, $\bfw=\bfw'+\bfw''$,  
the chamber $\gC$ given by $|u^{'(0)}_1|> \cdots> |u^{'(N-1)}_{w_{N-1}'}|> |u^{''(0)}_1|>\cdots > |u^{''(N-1)}_{w_{N-1}''}|$ and the stability condition $|t_1/t_2|<1$. The polarization $T^{1/2}$ 
 is chosen in such a way that $\Stab_{\gC,T^{1/2}}((\bla',\bla'');\bfz)$ is given by the formula in Theorem \ref{EllStab} for $\bla=(\bla',\bla'')$.  
\be
&&\hspace{-1cm}{\Stab}_{\gC,T^{1/2}}((\bla',\bla'');\bfz)\nn\\
&&\hspace{-1cm}=\Sym_0\Sym_1\cdots \Sym_{N-1} \Biggl[
\prod_{a\in \bla', b\in \bla'' \atop {c_a+1\equiv c_b 
}}{\theta(t_1 x'_a/x''_b)}\prod_{a\in \bla', b\in \bla'' \atop {c_b+1\equiv c_a 
}}{\theta(t_2 x'_a/x''_b)}
 \prod_{a\in \bla', b\in \bla'' \atop {c_a\equiv c_b 
 }}\frac{1}{\theta( x'_a/x''_b)\theta(\hbar x'_a/x''_b)}
 \nn\\
&&\qquad\hspace{2.5cm}\times \prod_{i=0}^{N-1}\left(\prod_{j=1}^{w_i'}\prod_{b\in \bla'' \atop {c_b\equiv i 
}}{\theta(\hbar u^{'(i)}_j/x''_b)}
\prod_{j=1}^{w_i''}\prod_{a\in \bla' \atop {c_a\equiv i 
}}{\theta( x'_a/u^{''(i)}_j)}\right)
\nn\\
&&\qquad\hspace{2.5cm}\times  {\Stab}_{\gC',T^{1/2'}}(\bla';\{ z_i'\hbar^{w_i''-v_i''+v_{i+1}''}\}) {\Stab}_{\gC'',T^{1/2''}}(\bla'';\{ z_i''\hbar^{v_i'-v_{i-1}'}\})\Biggr].
\en
Here ${\rm Sym}_i$ denotes the symmetrization over the Chern roots $\{x_a, a\in (\bla',\bla''), c_a\equiv i \ \mbox{mod}\ N \}$. 
\end{thm}
\noindent
{\it Proof.}\ The statement follows from
\be
&&\prod_{a\in \bla', b\in \bla'' \atop {c_a+1\equiv c_b 
}}{\theta(t_1 x'_a/x''_b)}\prod_{a\in \bla', b\in \bla'' \atop {c_b+1\equiv c_a 
}}{\theta(t_2 x'_a/x''_b)}
 \prod_{a\in \bla', b\in \bla'' \atop {c_a\equiv c_b 
 }}\frac{1}{\theta( x'_a/x''_b)\theta(\hbar x'_a/x''_b)}
 \nn\\
&&\times \prod_{i=0}^{N-1}\left(\prod_{j=1}^{w_i'}\prod_{b\in \bla'' \atop {c_b\equiv i 
}}{\theta(\hbar u^{'(i)}_j/x''_b)}
\prod_{j=1}^{w_i''}\prod_{a\in \bla' \atop {c_a\equiv i 
}}{\theta( x'_a/u^{''(i)}_j)}\right)
\times {S}_{\bla'}(\bfx',\bu';p,t_1,t_2) {S}_{\bla''}(\bfx'',\bu'';p,t_1,t_2)\nn\\
&&={S}_{(\bla',\bla'')}(\bfx,\bu;p,t_1,t_2)
\en
and 
\be
&&\sum_{\bft\in \Gamma_{\bla'}}W^{\bla'}(\bfx',\bu',\{ z_i'\hbar^{w_i''-v_i''+v_{i+1}''}\};p,t_1,t_2)
\sum_{\bft\in \Gamma_{\bla''}}W^{\la''}(\bfx'',\bu'',\{ z_i''\hbar^{v_i'-v_{i-1}'}\} ;p,t_1,t_2)\nn\\
&&\qquad =\sum_{\bft\in \Gamma_{(\bla',\bla'')}}W^{(\bla',\bla'')}(\bfx,\bu,\bfz ;p,t_1,t_2).
\en
In particular the shifts in the K\"ahler parameters are determined following the general argument by Botta\cite{Botta}.
Namely we have the $\rK$-theory class
\bea
\Delta&=&t_1\sum_{k=0}^{N-1}\Hom(\cV''_k,\cV_{k+1}')+\sum_{k=0}^{N-1}\Hom(\cV_k'',\cW_k')
-\sum_{k=0}^{N-1}\Hom(\cV_k'',\cV_k')\nn\\
&=&t_1\sum_{k=0}^{N-1}\sum_{i=1}^{v''_k}\sum_{j=1}^{v_{k+1}'}\frac{x^{(k+1)'}_j}{x^{(k)''}_i}
+\sum_{k=0}^{N-1}\sum_{i=1}^{v_k''}\sum_{j=1}^{w_{k}'}\frac{u^{'(k)}_j}{x^{(k)''}_i}
-\sum_{k=0}^{N-1}\sum_{i=1}^{v_k''}\sum_{j=1}^{v_{k}'}\frac{x^{(k)'}_j}{x^{(k)''}_i}.
\ena
Hence 
\bea
&&\delta=\det\Delta=\mbox{const.}\prod_{k=0}^{N-1}\left(
\prod_{v_k'+1}^{v_k}(x^{(k)''}_i)^{-w_k'+v_k'-v_{k+1}'}\prod_{j=1}^{v_k'}(x^{(k)'}_i)^{-v_k''+v_{k-1}''}
\right). 
\ena
\qed

Noting 
\be
K_I(\bfx,\bfu)&=&\frac{
\mbox{
$\ds{\prod_{a\in\bla', b\in \bla'' \atop {c_a+1\equiv c_b}}{\theta(t_2x_b''/x_a')}
\prod_{a\in\bla', b\in \bla'' \atop {c_b+1\equiv c_a}}{\theta(t_1x_b''/x_a')}
\prod_{i=0}^{N-1}\left(\prod_{j=1}^{w_i''}\prod_{a\in\bla' \atop c_a\equiv i}{\theta(x_a'/u^{''(i)}_j)}
\prod_{j=1}^{w_i'}\prod_{b\in\bla' \atop c_b\equiv i}{\theta(x_b''/u^{'(i)}_j)}\right)
}$
}
}{
\mbox{
$\ds{\prod_{a\in\bla', b\in \bla'' \atop {c_a\equiv c_b}}\theta(x_a'/x_b'')\theta(x_b''/x_a')}$
}
}\nn\\
&&\hspace{0.5cm}\times 
K_I(\bfx',\bfu')K_I(\bfx'',\bfu''),\\[2mm]
\vep((\bla',\bla''))&=&\vep(\bla')+\vep(\bla'')+\sum_{a\in \bla', b\in \bla''\atop c_a+1\equiv c_b}1+\sum_{a\in \bla'', b\in \bla'\atop c_a+1\equiv c_b}1+\sum_{i\in I}w'_iv''_i
,
\en
one obtains the shuffle product formula for $\widehat{\Stab}_{\gC,T^{1/2}}((\bla',\bla'');\bfz)$ from \eqref{def:hatSla}. 
\begin{cor}
\bea
&&\hspace{-1.5cm}\widehat{\Stab}_{\gC,T^{1/2}}((\bla',\bla'');\bfz)\nn\\
&&\hspace{-1.5cm}=\Sym_0\Sym_1\cdots \Sym_{N-1} \Biggl(\prod_{a\in \bla', b\in \bla'' \atop {c_b+1\equiv c_a 
}}\left(-\frac{\theta(t_2 x'_a/x''_b)}{\theta(t_1x''_b/x'_a)}\right)
\prod_{a\in \bla', b\in \bla'' \atop {c_a+1\equiv c_b 
}}\left(-\frac{\theta(t_1 x'_a/x''_b)}{\theta(t_2x''_b/x'_a)}\right)
\nn\\
&&\qquad\qquad\qquad\qquad\times
\prod_{i=0}^{N-1}\prod_{j=1}^{w_i'}\prod_{b\in \bla'' \atop {c_b\equiv i 
}}\left(-\frac{\theta(\hbar u^{'(i)}_j/x''_b)}{\theta(x''_b/u^{'(i)}_j)} \right)
\times \prod_{a\in \bla', b\in \bla'' \atop {c_a\equiv c_b 
}}\frac{\theta( x''_b/x'_a)}{\theta(\hbar x'_a/x''_b)}\nn\\
&&\qquad\qquad\times  \widehat{\Stab}_{\gC',T^{1/2'}}(\bla';\{ z_i'\hbar^{w_i''-v_i''+v_{i+1}''}\})
\widehat{\Stab}_{\gC'',T^{1/2''}}(\bla'';\{ z_i''\hbar^{v_i'-v_{i-1}'}\})\Biggr).
\lb{shufflehStabgen}
\ena
\end{cor}

We relate the dynamical parameters $\gz_i$ in $U_{t_1,t_2,p}(\gl_{N,tor})$ to the K\"ahler parameters $z_i$ as follows.  
\bea
&&z'_i=z_{0,i}\hbar^{v_i'-v_{i-1}'}\ \mbox{for} \ \cM(\bfv',\bfw'),\nn\\
&&z''_i=z_{0,i}\hbar^{v_i''-v_{i-1}''}\ \mbox{for} \ \cM(\bfv',\bfw'),\lb{relzz0}\\ 
&&\gz_i=z_i=z_{0,i}\hbar^{v_i-v_{i-1}}\ \mbox{for} \ \cM(\bfv,\bfw),\nn
\ena
where   $z_{0,i}$, $i\in I$ denote certain standard K\"ahler parameters.  One then can rewrite the shifts in \eqref{shufflehStabgen} in terms of $\gz_i$ as 
\bea
&&z'_i\hbar^{w_i''-v_i''+v_{i+1}''}=\gz_i\hbar^{w_i''-2v_i''+v_{i-1}''+v_{i+1}''},\qquad
z''_i\hbar^{v_i'-v_{i-1}'}=\gz_i.\lb{dynamicalshift}
\ena
In the next section, we show that the shuffle product formula  \eqref{shufflehStabgen} with \eqref{dynamicalshift} is equivalent to a composition of the vertex operators. See Sec.\ref{shufflebasic}. 

Similarly, we have the following shuffle product formula for $\widetilde{\Stab}_{\gC,T^{1/2}}(\bla;\bgz)$. 
\begin{cor}\lb{cor:shuffletStabgen}
\bea
&&\hspace{-1.5cm}\widetilde{\Stab}_{\gC,T^{1/2}}((\bla',\bla'');\bgz)\nn\\
&&\hspace{-1.5cm}=\Sym_0\Sym_1\cdots \Sym_{N-1} \Biggl(\prod_{a\in \bla', b\in \bla'' \atop {c_b+1\equiv c_a 
}}\left(-\frac{\theta(t_2 x'_a/x''_b)}{\theta(t_1x''_b/x'_a)}\right)
\prod_{a\in \bla', b\in \bla'' \atop {c_a+1\equiv c_b 
}}\left(-\frac{\theta(t_1 x'_a/x''_b)}{\theta(t_2x''_b/x'_a)}\right)
\nn\\
&&\qquad\qquad\qquad\qquad\times
\prod_{i=0}^{N-1}\prod_{j=1}^{w_i''}\prod_{a\in \bla' \atop {c_a\equiv i 
}}\left(-\frac{\theta(x'_a/u^{''(i)}_j)}{\theta(\hbar u^{''(i)}_j/x'_a)} \right)
\times \prod_{a\in \bla', b\in \bla'' \atop {c_a\equiv c_b 
}}\frac{\theta(\hbar x''_b/x'_a)}{\theta( x'_a/x''_b)}\nn\\
&&\qquad\qquad\times 
\widetilde{\Stab}_{\gC'',T^{1/2''}}(\bla'';\bgz
)
 \widetilde{\Stab}_{\gC',T^{1/2'}}(\bla';
 \{ \gz_i\hbar^{w_i''-2v_i''+v_i''+v_{i+1}''}\}
 )
\Biggr).
\lb{shuffletStabgen}
\ena
\end{cor}

\medskip

\noindent
{\it Remark.} \ 
In this occasion,  we remark the shuffle product formula for 
 the elliptic stable envelopes for the cotangent bundle of the partial flag varieties $T^*\F_\la$ obtained in  \cite{Konno17, 
 Konno18}.  Here $T^*\F_\la$ is isomorphic to the quiver variety of $A_{N-1}$ type with $\bfv=(\la^{(1)},\cdots,\la^{(N-1)})$, $0\leq \la^{(1)}\leq \cdots \leq \la^{(N-1)}\leq n$ and $\bfw=(n\delta_{i,N-1})$ in the notation of \cite{Konno17}.  In Sec.\ref{sec:intro}, we denote $\F_\la$ by $fl(\bfv)$. 
It was derived by considering a composition of the vertex operators $\Phi_{\mu_j}(u_j)$ $(1\leq j\leq n, \mu_j\in \{1,\cdots,N\})$ of the elliptic quantum group $U_{q,p}(\slnh)$.  
   Eq.(5.11) in \cite{Konno17} indicates 
the dynamical shift
\be
&&\Pi_{i,i+1}\ \mapsto\ \Pi_{i,i+1}q^{-2\sum_{j=1}^{n}\langle\bep_{\mu_j'},h_{i}\rangle}
\en
in the first function $F$ and no dynamical shifts in the second function $G$. 
Here the sum $\sum_{j=1}^{n}\langle\bep_{\mu_j'},h_{i}\rangle$ is the weight of the composition of the vertex operators
\be
\phi_{\mu'_1,\cdots,\mu'_n}(\bfu)&=&\Phi_{\mu'_1}(u_1)\cdots\Phi_{\mu'_n}(u_n)\\
&=&\int d\bfx\ \Phi_N(u_1)\cdots \Phi_{N}(u_n) \prod_{1\leq j\leq N-1}^{\curvearrowleft}:F_j(x^{(j)}_1)\cdots F_j(x^{(j)}_{v_j}): \nn\\
&&\times \prod_{j=1}^{N-1}\prod_{1\leq a<b\leq v_j}
<F_j(x^{(j)}_a)F_j(x^{(j)}_b)>^{\Sym}\times \widehat{\Stab}_{\gC,T^{1/2}}(I';\boldsymbol{\Pi}). 
\en
Here $I'=(I_1',\cdots,I_{N}')$ denotes a partition of $[1,n]$ determined by an index set $(\mu'_1,\cdots,\mu_n')$. $F_k(x)\ (1\leq k\leq N-1)$ and $\Phi_N(u)$  denote the elliptic currents and the top component of the vertex operator, respectively. See Appendix C in \cite{Konno17} for details. One can identify the elliptic weight function $\widetilde{W}_I(\bfx,\bfu,\boldsymbol{\Pi})$ with the elliptic stable envelope $ \widehat{\Stab}_{\gC,T^{1/2}}(I;\boldsymbol{\Pi})$ for $\E_\rT(T^*\F_\la)$. 
See 
Sec.5 in \cite{Konno18}. Since $\Phi_N(u)$ and $F_j(x)$ carry the weight $-\bep_N$ and $-\al_j$, respectively, one finds
\be
&&-\sum_{j=1}^{n}\langle\bep_{\mu_j'},h_{i}\rangle
 =\langle n\bep_N+\sum_{j=1}^{N-1}v_j\al_j,h_i\rangle=n\delta_{i,N-1}-2v_i+v_{i-1}+v_{i+1}
\en  
$i=1,\cdots,N-1$ with $v_0=v_N=0$. Here $v_j=\# F_j$'s in $\phi_{\mu'_1,\cdots,\mu'_n}(\bfu)$.  Therefore the dynamical shifts of $\Pi_{i,i+1}$ in the function $F$ and $G$ are  the same  as in the RHS of \eqref{dynamicalshift} in terms of $\gz_i$. 
Therefore identifying $\Pi_{i,i+1}$ with $\gz_i$ and using the same relations as in \eqref{relzz0}, one finds that the formula (5.11) in \cite{Konno17} is equivalent to the shuffle product formula for the elliptic stable envelopes for $T^*\F_\la$, which can be derived by the  
general argument by Botta\cite{Botta}.

\section{Vertex Operators for $\cM(\bfv,\bfw)$ }\lb{sec:VO}

We construct the type I and the type II dual vertex operators of $U_{t_1,t_2,p}(\gl_{N,tor})$ by combining 
the level-$(1,L)$ representation and the elliptic stable envelope for $\E_\rT(\cM(\bfv,\bfw))$.  Our construction is  step-by-step starting from the basic ones  $\Phi^{(k)}_\la(u)$ and $\Psi^{*(k)}_{\la}(u)$ for the case $|\bfw|=1$, $\bfw=(\delta_{i,k})_{i\in I}$ and $|\bfv|=|\la|$  for each $k\in I$ (Sec.\ref{sec:VObasic}), then composing $w_k$ of them we obtain the  ones $\Phi^{(k)}_{\blla^{(k)}}(\bfu^{(k)})$ and $\Psi^{*(k)}_{\blla^{(k)}}(\bfu^{(k)})$ for the case $|\bfw|=|w_k|$, $\bfw=(w_k\delta_{i,k})_{i\in I}$ and $|\bfv|=\sum_{j=1}^{w_k}|\la^{(k)}_j|$ (Sec.\ref{sec:VOkgen}).  Finally  composing $\Phi^{(k)}_{\blla^{(k)}}(\bfu^{(k)})\ (\mbox{resp.}\ \Psi^{*(k)}_{\blla^{(k)}}(\bfu^{(k)}))$ for all $k\in I$,  we obtain the most general ones $\Phi_{\blla}(\bfu)$ and $\Psi^*_{\blla}(\bfu)$ for  $\cM(\bfv,\bfw)$ with $\bfv, \bfw\in \N^N$ (Sec.\ref{sec:VOgen}). 
As the case of $\gl_{1,tor}$\cite{KS},  these vertex operators turn out to be  intertwining operators of the  
$U_{t_1,t_2,p}(\glnt)$-modules w.r.t. the standard comultiplication. See Sec.\ref{sec:intrel} for  $\Phi^{(k)}_\la(u)$ and $\Psi^{*(k)}_{\la}(u)$.

\subsection{Operator product expansion (OPE) of the elliptic currents}

In the level $(1,L)$ representation we have $\gamma=\hbar^{1/2}$,  $p^*=p\hbar^{-1}$ and 
\bea
x^+_i(x_a)x^+_i(x_b)
&=&x_a^2\frac{(\hbar^{-1} x_b/x_a;p^*)_\infty(x_b/x_a;p^*)_\infty}{(p^*x_b/x_a;p^*)_\infty(p^*\hbar x_b/x_a;p^*)_\infty}:x^+_i(x_a)x^+_i(x_b):\nn\\
&=&{<x^+_i(x_a)x^+_i(x_b)>^{Sym}}\times\frac{\theta^*(\hbar^{-1} x_b/x_a)}{\theta^*(x_a/x_b)}:x^+_i(x_a)x^+_i(x_b):,\\
x^-_i(x_a)x^-_i(x_b)
&=&x_a^2\frac{(\hbar x_b/x_a;p)_\infty(x_b/x_a;p)_\infty}{(px_b/x_a;p)_\infty(p\hbar^{-1}x_b/x_a;p)_\infty}:x^-_i(x_a)x^-_i(x_b):\nn\\
&=&{<x^-_i(x_a)x^-_i(x_b)>^{Sym}}\times\frac{\theta(x_b/x_a)}{\theta(\hbar x_a/x_b)}:x^-_i(x_a)x^-_i(x_b):
\lb{OPExmxm}
\ena
for $|x_b/x_a|<1$ with 
\bea
&&<x^+_i(x_a)x^+_i(x_b)>^{Sym}=\hbar^{-1/2}x_ax_b\frac{(x_a/x_b,x_b/x_a;p^*)_\infty}{( p^*\hbar x_a/x_b, p^*\hbar x_b/x_a;p^*)_\infty},\\
&&<x^-_i(x_a)x^-_i(x_b)>^{Sym}=\hbar^{-1/2}x_ax_b\frac{( \hbar x_a/x_b, \hbar x_b/x_a;p)_\infty}{(px_a/x_b, px_b/x_a;p)_\infty}.
\ena
Here we set 
\be
&&(z_1,z_2;p)_\infty=(z_1;p)_\infty(z_2;p)_\infty,\qquad \mbox{etc.}
\en
Hence in the sense of analytic continuation we have
\bea
&&x^+_i(x_a)x^+_i(x_b)=-\frac{\theta(\hbar x_a/x_b)}{\theta(\hbar x_b/x_a)}x^+_i(x_b)x^+_i(x_a),\\
&&x^-_i(x_a)x^-_i(x_b)=-\frac{\theta(\hbar x_b/x_a)}{\theta(\hbar x_a/x_b)}x^-_i(x_b)x^-_i(x_a).
\ena
Similarly, 
\bea
x^+_{i+1}(x_a)x^+_i(x_b)
&=&x_a^{-1}\frac{(p^*t_1x_b/x_a;p^*)_\infty}{(t_2^{-1}x_b/x_a;p^*)_\infty}:x^+_{i+1}(x_a)x^+_i(x_b):,\\
x^-_{i+1}(x_a)x^-_i(x_b)
&=&x_a^{-1}\frac{(pt_2^{-1}x_b/x_a;p)_\infty}{(t_1x_b/x_a;p)_\infty}:x^-_{i+1}(x_a)x^-_i(x_b):
\ena
for $|x_b/x_a|<1$ and
\bea
x^+_i(x_b)x^+_{i+1}(x_a)
&=&x_b^{-1}\frac{(p^*t_2x_a/x_b;p^*)_\infty}{(t_1^{-1}x_a/x_b;p^*)_\infty}:x^+_{i+1}(x_a)x^+_i(x_b):,\\
x^-_i(x_b)x^-_{i+1}(x_a)
&=&x_b^{-1}\frac{(pt_1^{-1}x_a/x_b;p)_\infty}{(t_2x_a/x_b;p)_\infty}:x^-_{i+1}(x_a)x^-_i(x_b):
\ena
for $|x_a/x_b|<1$. 
Hence noting \eqref{CEalal} we have 
\bea
&&x^+_{i+1}(x_a)x^+_i(x_b)=-\frac{\theta^*(t_1^{-1}x_a/x_b)}{\theta^*(t_2^{-1}x_b/x_a)}x^+_i(x_b)x^+_{i+1}(x_a),\\
&&x^-_{i+1}(x_a)x^-_i(x_b)=-\frac{\theta(t_2x_a/x_b)}{\theta(t_1x_b/x_a)}x^-_i(x_b)x^-_{i+1}(x_a).
\ena
In addition, we have
\bea
&&x^+_j(z)x^-_i(w)=:x^-_i(w)x^+_j(z):\left\{\mmatrix{\frac{1}{(1-qw/z)(1-q^{-1}w/z)}&j\equiv i\cr
1-\kappa w/z&j\equiv i-1\cr
1-\kappa^{-1}w/z&j\equiv i+1\cr
1&\mbox{otherwise}\cr}\right.,\lb{xpxm}
\ena
for $|w/z|<1$ and
\bea
&&x^-_i(w)x^+_j(z)=:x^-_i(w)x^+_j(z):\left\{\mmatrix{\frac{1}{(1-qz/w)(1-q^{-1}z/w)}&j\equiv i\cr
1-\kappa^{-1}z/w&j\equiv i-1\cr
1-\kappa z/w&j\equiv i+1\cr
1&\mbox{otherwise}\cr}\right.\lb{xmxp}
\ena
for $|z/w|<1$. 
Hence from \eqref{CEalal}, we have
\bea
&&x^+_i(z)x^-_j(w)=x^-_j(w)x^+_i(z),\qquad i,j\in I
\ena

Let  
$A_{k,m}$ and 
$A'_{k,m}=\frac{1-p^{*m}}{1-p^m}\gamma^{m}A_{k,m}\ (i\in I, m\in \Z\backslash \{0\})$ be the Heisenberg generators satisfying
\bea
&&[A_{k,m},\al_{i,n}]=\delta_{k,i}\frac{[m]_q}{m}\frac{\gamma^m-\gamma^{-m}}{q-q^{-1}}\frac{1-p^{m}}{1-p^{*m}}\gamma^{-m}\delta_{m+n,0},\\
&&[A'_{k,m},{\al}'_{i,n}]=\delta_{k,i}\frac{[m]_q}{m}\frac{\gamma^m-\gamma^{-m}}{q-q^{-1}}\frac{1-p^{*m}}{1-p^m}\gamma^{m}\delta_{m+n,0},\\
&&[A_{i,m},A_{j,n}]=-\frac{1}{(q-q^{-1})^2}\frac{\gamma^{m}-\gamma^{-m}}{m}\frac{1-p^m}{1-p^{*m}}\gamma^{-m}
\left(\frac{t_1^{Nm}t_1^{(i-j)m}}{1-t_1^{Nm}}+\frac{t_2^{-(i-j)m}}{1-t_2^{Nm}}\right)\delta_{m+n,0}\nn\\
&&\qquad\qquad\quad=-[A_{j,n},A_{i,m}]\quad\qquad (i\leq j),\lb{ComAA}\\
&&[A'_{i,m},A'_{j,n}]=-\frac{1}{(q-q^{-1})^2}\frac{\gamma^{m}-\gamma^{-m}}{m}\frac{1-p^{*m}}{1-p^{m}}\gamma^{m}
\left(\frac{t_1^{Nm}t_1^{(i-j)m}}{1-t_1^{Nm}}+\frac{t_2^{-(i-j)m}}{1-t_2^{Nm}}\right)\delta_{m+n,0}\nn\\
&&\qquad\qquad\quad=-[A'_{j,n},A'_{i,m}]\quad\qquad (i\leq j),\lb{ComApAp}
\ena
Set
\be
\Psi^{*(k)}_\emptyset(u)
&=&e^{-\Lambda_k}e^{Q_{\Lambda_k}}(-\hbar^{1/2}u)^{-h_{\Lambda_k}}\\
&&\qquad\times\exp\left\{-\sum_{m>0}\frac{1}{[m]_q}A_{k,-m} (\hbar^{1/2}u)^m
\right\}\exp\left\{\sum_{m>0}\frac{1}{[m]_q}A_{k,m}(\hbar^{1/2}u)^{-m}
\right\},\\
\Phi^{(k)}_\emptyset(u)
&=&e^{\Lambda_k}(-\hbar^{1/2}u)^{h_{\Lambda_k}}\\
&&\qquad \times\exp\left\{\sum_{m>0}\frac{1}{[m]_q}A'_{k,-m} (\hbar^{1/2}u)^m
\right\}\exp\left\{-\sum_{m>0}\frac{1}{[m]_q}A'_{k,m}(\hbar^{1/2}u)^{-m}
\right\}.
\en
It is easy to verify the following. 
\bea
&&\Psi^{*(k)}_\emptyset(u)x^+_j(z)=\left((-\hbar^{1/2}u)^{-1}\frac{(p^*z/u;p^*)_\infty}{(\hbar^{-1}z/u;p^*)_\infty}\right)^{\delta_{k,j}}:\Psi^{*(k)}_\emptyset(u)x^+_j(z):,\lb{Psisexp2}\\
&&\Psi^{*(k)}_\emptyset(u)x^-_j(z)=(z-\hbar^{1/2}u)^{\delta_{k,j}}:\Psi^{*(k)}_\emptyset(u)x^-_j(z):,\lb{Psisexm}\\
&&\Phi^{(k)}_\emptyset(u)x^+_j(z)=(z-\hbar^{1/2}u)^{\delta_{k,j}}:\Phi^{(k)}_\emptyset(u)x^+_j(z):,\lb{Phiexm}\\
&&\Phi^{(k)}_\emptyset(u)x^-_j(z)=\left((-\hbar^{1/2}u)^{-1}\frac{(p\hbar^{-1}z/u;p)_\infty}{(z/u;p)_\infty}\right)^{\delta_{k,j}}:\Phi^{(k)}_\emptyset(u)x^+_j(z):,\lb{Phiexp2}
\ena
for $|z/u|<1$ and
\bea
&&x^+_j(z)\Psi^{*(k)}_\emptyset(u)=\left(z^{-1}\frac{(p^*\hbar u/z;p^*)_\infty}{(u/z;p^*)_\infty}\right)^{\delta_{k,j}}:\Psi^{*(k)}_\emptyset(u)x^+_j(z):,
\lb{xpPsise2}\\
&&x^-_j(z)\Psi^{*(k)}_\emptyset(u)=(z-\hbar^{1/2}u)^{\delta_{k,j}}:\Phi^{(k)}_\emptyset(u)x^-_j(z):,
\lb{xmPsise}\\
&&x^+_j(z)\Phi^{(k)}_\emptyset(u)=(z-\hbar^{1/2}u)^{\delta_{k,j}}:\Phi^{(k)}_\emptyset(u)x^+_j(z):,
\lb{xpPhie}\\
&&x^-_j(z)\Phi^{(k)}_\emptyset(u)=\left(z^{-1}\frac{(pu/z;p)_\infty}{(\hbar u/z;p)_\infty}\right)^{\delta_{k,j}}:\Phi^{(k)}_\emptyset(u)x^+_j(z):
\lb{xpPhie2}
\ena
for $|u/z|<1$. 
Hence we have
\bea
&&\Psi^{*(k)}_\emptyset(u)x^+_j(z)=\left(-\frac{\theta^*( z/u)}{\theta^*(\hbar u/z)}\right)^{\delta_{k,j}}x^+_j(z)\Psi^{*(k)}_\emptyset(u),\lb{commPsisexp}\\
&&\Psi^{*(k)}_\emptyset(u)x^-_j(z)= x^-_j(z)\Psi^{*(k)}_\emptyset(u),\lb{commPsisexm}\\
&&\Phi^{(k)}_\emptyset(u)x^+_j(z)= x^+_j(z)\Phi^{(k)}_\emptyset(u),\lb{commPhiexp}\\
&&\Phi^{(k)}_\emptyset(u)x^-_j(z)=\left(-\frac{\theta(\hbar u/z)}{\theta(z/u)}\right)^{\delta_{k,j}}x^-_j(z)\Phi^{(k)}_\emptyset(u).\lb{commPhiexm}
\ena

\subsection{The basic vertex operators : the case $\cM(\bfv,\bfw=(\delta_{i,k}))$}\lb{sec:VObasic}

Fix $k\in I$, let $\la$ be a partition with  color $k$. See Definition \ref{def:colork}. Let  
  $X={\cal M}(\bold{v},\bold{w})$ be  the  quiver variety of type $A^{(1)}_{N-1}$ with 
$\bfv=(v_i)_{i\in I}$, $v_i=\#\{\square \in \la\ |\ c_{\square}\equiv i\ \mbox{mod}\ N \}$ and 
 $\bold{w}=(w_i)_{i\in I}$, $w_i=\delta_{i,k}$. 
We denote by $\Lambda_{\bfv,\bfw}$ the corresponding set of fixed points. 

Consider the elliptic stable envelopes for $\E_\rT(X)$
\bea
&&\hspace{-1.5cm}\widehat{\Stab}_{\gC,T^{1/2}}(\la;\bgz)=\Sym_0\cdots\Sym_{N-1}\left(\widehat{S}_\la(\bfx,u;p,t_1,t_2)\sum_{\bft\in \Gamma_\la}W^\la_{\bft}(\bfx,u,\bgz;p,t_1,t_2)\right)\lb{basicStab}
,\\
&&\hspace{-1.5cm}\widetilde{\Stab}^*_{\gC^*,T^{1/2}_{opp}}(\la;\bgz^{*-1})=\Sym_0\cdots\Sym_{N-1}\left(\widetilde{S}_\la(\bfx,u;p^*,t_1,t_2)\sum_{\bft\in \Gamma_\la}W^\la_{\bft}(\bfx,u,\bgz^{*-1};p^*,t_1,t_2)\right)\nn\\
&&
\ena
with  chambers $\gC : | t_1/t_2|=|\kappa^{-1}|<1$ and
\be
&&\widehat{S}_\la(\bfx,u;p,t_1,t_2)\nn\\
&&=
\prod_{a,b\in\la \atop {c_a+1\equiv c_b\atop 
\rho_a+1<\rho_b}}\left(-\frac{\theta(t_1x_a/x_b)}{\theta(t_2x_b/x_a)}\right)
\prod_{a,b\in\la \atop {c_a+1\equiv c_b\atop 
\rho_a+1>\rho_b}}\left(-\frac{\theta(t_2x_b/x_a)}{\theta(t_1x_a/x_b)}\right)
\prod_{a\in\la \atop {c_a\equiv k\atop 
\rho_a> \rho_{r}}}\left(-\frac{\theta(\hbar u/x_a)}{\theta(x_a/u)}\right)
\prod_{a,b\in\la \atop {c_a\equiv c_b\atop 
\rho_a<\rho_b}}\frac{\theta(x_b/x_a)}{\theta(\hbar x_a/x_b)},\\
&&\widetilde{S}_\la(\bfx,u;p^*,t_1,t_2)\nn\\
&&=
\prod_{a,b\in\la \atop {c_a+1\equiv c_b\atop 
\rho_a+1<\rho_b}}\left(-\frac{\theta^*(t_1x_a/x_b)}{\theta^*(t_2x_b/x_a)}\right)
\prod_{a,b\in\la \atop {c_a+1\equiv c_b\atop 
\rho_a+1>\rho_b}}\left(-\frac{\theta^*(t_2x_b/x_a)}{\theta^*(t_1x_a/x_b)}\right)
\prod_{a\in\la \atop {c_a\equiv k\atop 
\rho_a\leq \rho_{r}}}\left(-\frac{\theta^*(x_a/u)}{\theta^*(\hbar u/x_a)}\right)
\prod_{a,b\in\la \atop {c_a\equiv c_b\atop 
\rho_a>\rho_b}}\frac{\theta^*(\hbar x_a/x_b)}{\theta^*(x_b/x_a)}.
\en
where $r$ denotes the root box in $\la$.

Let $\ds{\prod_{a\in \la} x^-_{c_a}(x_a)}$ denotes a product in the increasing order w.r.t. $\prec$.  Namely 
\be
&&\prod_{a\in \la}x^-_{c_a}(x_a)=x^-_{c_{a_1}}(x^-_{a_1})\cdots x^-_{c_{a_{|\la|}}}(x^-_{a_{|\la|}})
\en 
with $\rho_{a_i}<\rho_{a_j}$ for $1\leq i<j\leq |\la|$.  See Sec.\ref{sec:TFP}. Similarly, let $\ds{\prod^{\curvearrowleft}_{a\in \la\atop \rho_a> \rho_r}x^+_{c_a}(x_a)}$  denotes a product in the decreasing order w.r.t. $\prec$ : 
\be
\prod^{\curvearrowleft}_{a\in \la}x^+_{c_a}(x_a)=
x^+_{c_{a_{|\la|}}}(x^+_{a_{|\la|}})\cdots x^+_{c_{a_{1}}}(x^+_{a_{1}})
\en
\begin{dfn}\lb{TypeITypeIIVObasic}
For $k\in I$, we define the type I and the type II dual vertex operators $\Phi^{(k)}(u)$ and $\Psi^{*(k)}(u)$, respectively, as follows. 
\bea
&&\Phi^{(k)}(u)\ :\ \F^{(1,L)}_{\bxi,\blLa_L,\blLa_M} \to \F^{(0,-1) (k)}_u\tot \F^{(1,L+1)}_{\left(\xi_0,\cdots,-\frac{\xi_k}{h^{1/2}u},\cdots,\xi_{N-1}\right),\blLa_L+\La_k,\blLa_M},\\
&&\Psi^{*(k)}(u)\ :\   \F^{(1,L)}_{\bxi,\blLa_L,\blLa_M}
\tot\F^{(0,-1) (k)}_u\to \F^{(1,L-1)}_{\left(\xi_0,\cdots,-{h^{1/2}u}{\xi_k},\cdots,\xi_{N-1}\right),\blLa_L-\La_k,\blLa_M+\La_k}
\ena
whose components are defined by
\bea
&&\Phi^{(k)}(u)=\sum_{\la\in \cP^{(k)}}\ket{\la}^{(k)}_u\tot \Phi^{(k)}_\la(u),\\
&&\Psi^{*(k)}_\la(u)\eta =\Psi^{*(k)}(u)(\eta\tot \ket{\la}^{(k)}_u),\qquad \forall \eta\in \F^{(1,L)}_{\bxi,\blLa_L,\blLa_M}.
\ena
\bea
&&\Phi^{(k)}_{\la}(u)=\int_\cC\prod_{a\in \la}\underline{dx_a}\prod_{i=0}^{N-1}
\left(:\prod_{a\in \la\atop {c_a\equiv i}}x^-_{i}(x_a) :\right) \ 
 \Phi^{(k)}_\emptyset(u)\ 
\prod_{a,b\in \la \atop {c_a\equiv c_b\atop \rho_a<\rho_b}}{<x^-_{c_a}(x_a)x^-_{c_b}(x_b)>^{Sym}}\nn\\
&&\hspace{3cm}\times \widehat{\Stab}_{\gC,T^{1/2}}(\la;\bgz),\lb{def:TypeIbasic}
\\
&&\Psi^{*(k)}_{\la}(u)=\int_{\cC^*}\prod_{a\in \la}\underline{dx_a}\ 
\widetilde{\Stab}^*_{\gC,T^{1/2}_{opp}}(\la;\bgz^{*-1})\prod_{i=0}^{N-1}
\left(:\prod_{a\in \la\atop {c_a\equiv i}}x^+_{i}(x_a) :\right) \ 
 \Psi^{*(k)}_\emptyset(u)\ \nn\\
&&\hspace{3cm}\times 
\prod_{a,b\in \la \atop {c_a\equiv c_b\atop \rho_a>\rho_b}}{<x^+_{c_a}(x_a)x^+_{c_b}(x_b)>^{Sym}}.
\lb{TypeIIbasic}
\ena
\end{dfn}
Note that  $\Phi^{(k)}_{\la}(u)$ and $\Psi^{*(k)}_{\la}(u)$ carry the $h_i$-weights 
$\langle h_i,\La_k-\sum_{j\in I}v_j\al_j\rangle$ and  $-\langle h_i,\La_k-\sum_{j\in I}v_j\al_j\rangle$, respectively. 
Hence by the identity 
\bea
&&
|R^{(k)}_i(\la)|-|A^{(k)}_i(\la)| 
=-\delta_{i,k}+\sum_{j\in I}a_{ij}v_j\lb{IdGeoRep}
\ena
for $\la\in \cP^{(k)}$ and $v_i=\#\{ \square \in \la\ |\ c_{\square}\equiv i\ \mbox{mod}\ N\}$,  the $h_i$-weight $\wt(\la)^{(k)}_i$ of $\ket{\la}^{(k)}_u$ (Corollary \ref{prop:wt}) compensates the one  of  the vertex operator $\Phi^{(k)}_{\la}(u)$.  Similarly, the weight  of $\Psi^{*(k)}_{\la}(u)$ coincides with the one of $\ket{\la}^{(k)}_u$.  These balancing conditions  of the weights ensure the existence of the vertex operators $\Phi^{(k)}_{\la}(u)$ and $\Psi^{*(k)}_{\la}(u)$.  

\vspace{2mm}
\noindent
{\it Remark.}\ Note also that $\Phi_\la(u)$ and $\Psi^*_\la(u)$ depend on $\la\in \La_{\bfv,\bfw}$ only through 
$\widehat{\Stab}_{\gC,T^{1/2}}(\la;\bgz)$ and $\widetilde{\Stab}^*_{\gC,T^{1/2}_{opp}}(\la;\bgz^{*-1})$, respectively.  This is due to the fact that the integrands of them 
are symmetric in $\{x_a, a\in \la, c_a\equiv i \}$ for each $i\in I$ so that they are 
 independent from the choice of $\la\in \La_{\bfv,\bfw}$ except for the elliptic stable envelopes.  
   The same is true for the vertex operators for general $\cM(\bfv,\bfw)$ we will discuss in the following subsections. 
\qed

The following proposition indicates a representation theoretical origin of some factors in $\widehat{S}_\la(\bfx,u;p,t_1,t_2)$  and $\widetilde{S}_\la(\bfx,u;p^*,t_1,t_2)$.
\begin{prop}
\bea
\Phi^{(k)}_{\la}(u)&=&\int\prod_{a\in \la}\underline{dx_a}\prod_{a\in \la\atop \rho_a\leq \rho_r}x^-_{c_a}(x_a)\cdot \Phi^{(k)}_\emptyset(u)\cdot \prod_{a\in \la\atop \rho_a>\rho_r}x^-_{c_a}(x_a) \ 
\nn\\
&&\hspace{2cm}\times
\prod_{a,b\in \la \atop {c_a+1\equiv c_b\atop h_a<h_b}}\left(\frac{\theta(t_2x_b/x_a)}{\theta(t_1x_a/x_b)}\right)^2
\sum_{t\in \Gamma_\la}W^\la_{t}(\bfx,u,\bgz;p,t_1,t_2),\\
\Psi^{*(k)}_{\la}(u)&=&\int\prod_{a\in \la}\underline{dx_a}\ \sum_{t\in \Gamma_\la}W^\la_{t}(\bfx,u,\bgz;p^*,t_1,t_2)
\prod_{a,b\in \la \atop {c_a+1\equiv c_b\atop h_a<h_b}}\left(\frac{\theta^*(t_2x_b/x_a)}{\theta^*(t_1x_a/x_b)}\right)^2
\nn\\
&&\hspace{2cm}\times
\prod^{\curvearrowleft}_{a\in \la\atop \rho_a> \rho_r}x^+_{c_a}(x_a)\cdot \Psi^{*(k)}_\emptyset(u)\cdot \prod^{\curvearrowleft}_{a\in \la\atop \rho_a\leq \rho_r}x^+_{c_a}(x_a).
\ena
\end{prop}
\noindent
{\it Proof.}\ 
One has
\bea
&&\{a\in \blla\ |\ c_a+1\equiv c_b, \rho_a<\rho_b\}\nn\\
&&=\{a\in \blla\ |\  c_a+1\equiv c_b, \rho_a+1<\rho_b\}
\cup\{a\in \blla\ |\  c_a+1= c_b, h_a<h_b\},\lb{rhoa<rhob}\\
&&\{a\in \blla\ |\ c_a+1\equiv c_b, \rho_a+1>\rho_b\}\nn\\
&&=\{a\in \blla\ |\  c_a+1\equiv c_b, \rho_a>\rho_b\}
\cup\{a\in \blla\ |\  c_a+1= c_b, h_a<h_b\}.\lb{rhoa>rhob}
\ena
The statement follows from \eqref{commPhiexm}, \eqref{rhoa>rhob}, \eqref{rhoa<rhob} and 
\be
&&\prod_{a\in \la}x^-_{c_a}(x_a)=\prod_{i=0}^{N-1}\prod_{a\in \la\atop c_a\equiv i}x^-_i(x_a)
\prod_{a,b\in \la \atop {c_a+1\equiv c_b\atop \rho_a<\rho_b}}
\left(-\frac{\theta(t_1x_a/x_b)}{\theta(t_2x_b/x_a)}\right)
\prod_{a,b\in \la \atop {c_a+1\equiv c_b\atop \rho_a>\rho_b}}
\left(-\frac{\theta(t_2x_b/x_a)}{\theta(t_1x_a/x_b)}\right),\\
&&
\prod_{i=0}^{N-1}\prod_{a\in \la\atop c_a\equiv i}x^-_i(x_a)=
\prod_{i=0}^{N-1}\left(:\prod_{a\in \la\atop c_a\equiv i}x^-_i(x_a):\right)
\prod_{a,b\in \la \atop {c_a\equiv c_b\atop \rho_a<\rho_b}}
<x^-_{c_a}(x_a)x^-_{c_b}(x_b)>^{Sym}
\frac{\theta(x_b/x_a)}{\theta(\hbar x_a/x_b)}
\en
for the type I, and 
from \eqref{commPhiexm}, \eqref{rhoa>rhob}, \eqref{rhoa<rhob} 
\be
&&\prod^{\curvearrowleft}_{a\in \la}x^+_{c_a}(x_a)=\prod_{i=0}^{N-1}\prod^{\curvearrowleft}_{a\in \la\atop c_a\equiv i}x^+_i(x_a)
\prod_{a,b\in \la \atop {c_a+1\equiv c_b\atop \rho_a<\rho_b}}
\left(-\frac{\theta^*(t_1x_a/x_b)}{\theta^*(t_2x_b/x_a)}\right)
\prod_{a,b\in \la \atop {c_a+1\equiv c_b\atop \rho_a>\rho_b}}
\left(-\frac{\theta^*(t_2x_b/x_a)}{\theta^*(t_1x_a/x_b)}\right),\\
&&
\prod_{i=0}^{N-1}\prod^{\curvearrowleft}_{a\in \la\atop c_a\equiv i}x^+_i(x_a)=
\prod_{i=0}^{N-1}\left(:\prod_{a\in \la\atop c_a\equiv i}x^+_i(x_a):\right)
\prod_{a,b\in \la \atop {c_a\equiv c_b\atop \rho_a>\rho_b}}
<x^+_{c_a}(x_a)x^+_{c_b}(x_b)>^{Sym}
\frac{\theta^*(\hbar x_a/x_b)}{\theta^*(x_b/x_a)}.
\en
for the type II dual. 
\qed

\subsection{Shuffle product as a composition of vertex operators}\lb{shufflebasic}

We show that compositions of the vertex operators are consistent to the shuffle product formula for $\widehat{\Stab}_{\gC,{T^{1/2}}}(\la;\bgz)$. 
Let us consider the composition 
\be
&& (\id\tot\Phi^{(k)}(u_{1}))\circ\Phi^{(k)}(u_{2})\\
&&\qquad  :\ \F^{(1,L)}_{\bxi,\blLa_L,\blLa_M}\ \to   \ \F^{(0,-1)(k)}_{u_{2}}\tot \F^{(0,-1)(k)}_{u_1}\tot\F^{(1,L+2)}_{\left(\xi_0,\cdots,\frac{\xi_k}{(-\hbar^{1/2})^2u_1u_2},\cdots,\xi_{N-1}\right), \blLa_L+2\La_k, \blLa_M}.
\en
Its components $\Phi^{(k)}_{\bla}(u_1,u_2)$ are defined by 
\bea
&&(\id\tot\Phi^{(k)}(u_{1}))\circ\Phi^{(k)}(u_{2})=\sum_{n\in \Z_{\geq 0}}\sum_{\bla=(\la',\la'')\in \cP^{(k)}\times \cP^{(k)}\atop |\bla|=n}
\ket{\la''}^{(k)}_{u_2}\tot\ket{\la'}^{(k)}_{u_1}\tot\ \Phi^{(k)}_{\bla}(u_1,u_2).   \nn
\lb{TypeIMn2total}
\ena
One has
\be
&&\Phi^{(k)}_{\bla}(u_1,u_2)=\Phi^{(k)}_{\la'}(u_1)\Phi^{(k)}_{\la''}(u_2). 
\en

Let 
$\la', \la'' \in \cP^{(k)}
$ and $\cM(\bfv',\bfw'), \cM(\bfv'',\bfw'')$ be the corresponding $A^{(1)}_{N-1}$ quiver varieties with $|\bfv'|=|\la'|$,  
$|\bfv''|=|\la''|$ and $\bfw'=\bfw''=(\delta_{i,k})$, respectively. Let $u_1$ and $\{ x_a\}_{a\in \la'}$ ( resp. $u_2$ and $\{ x_b\}_{b\in \la''}$ ) be the framing weights and  the Chern roots associated with    
$\cM(\bfv',\bfw')$  ( resp. $\cM(\bfv'',\bfw'')$), respectively. 
We assume $|u_1|> |u_2|$ and  $\rho_a<\rho_b$ for $a\in \la', b\in \la''$. 
Let $ \widehat{\Stab}_{\gC',{T^{1/2}}}(\la';\bgz)$ and $ \widehat{\Stab}_{\gC'',T^{1/2}}(\la'';\bgz)$ be the corresponding elliptic stable envelopes given in \eqref{basicStab} with $\gC'=\gC'' : | t_1/t_2|<1$. Let $\Phi^{(k)}_{\la'}(u_1)$ and $\Phi^{(k)}_{\la''}(u_2)$ be the basic vertex operators given in the same form as \eqref{def:TypeIbasic} and consider the composition 
\be
&&\hspace{-1.8cm}\Phi^{(k)}_{\la'}(u_1)\Phi^{(k)}_{\la''}(u_2)
=\int\prod_{a\in \la'}\underline{dx'_a}\int\prod_{b\in \la''}\underline{dx''_b}\ 
\prod_{i=0}^{N-1}\left(:\prod_{a\in \la'\atop {c_a\equiv i}}x^-_{i}(x'_a) :\right) \ 
 \Phi^{(k)}_\emptyset(u_1)\  \widehat{\Stab}_{\gC',T^{1/2'}}(\la';\bgz)\nn\\
 &&\qquad\qquad\qquad \times \prod_{j=0}^{N-1}
\left(:\prod_{b\in \la''\atop {c_b\equiv j}}x^-_{j}(x''_b) :\right) \ 
 \Phi^{(k)}_\emptyset(u_2)\  \widehat{\Stab}_{\gC'',T^{1/2''}}(\la'';\bgz)\nn\\
&&\qquad\qquad\qquad\times\prod_{a,b\in \la' \atop {c_a\equiv c_b\atop \rho_a<\rho_b}}{<x^-_{c_a}(x'_a)x^-_{c_b}(x'_b)>^{Sym}}\ \prod_{c,d\in \la'' \atop {c_c\equiv c_d\atop \rho_c<\rho_d}}{<x^-_{c_c}(x''_c)x^-_{c_d}(x''_d)>^{Sym}}.
\en
We arrange the order of the elements in the integrand as follows\cite{Konno17,KS}. 
\begin{itemize}
\item[1.] Move $\widehat{\Stab}_{\gC',T^{1/2}}(\la';\bgz)$ to the right of $\Phi^{(k)}_\emptyset(u_2)$. Then $\gz_i=\hbar^{P_i+h_i}$ gets a shift by $\hbar^{\delta_{i,k}-2v''_i+v''_{i-1}+v''_{i+1}
}$. This is due to $\hbar^{h_i}e^{\Lambda_k}=\hbar^{\delta_{i,k}}e^{\Lambda_k}\hbar^{h_i}$ and 
\be
&&\hbar^{h_i}\left(:\prod_{b\in \la''\atop {c_b\equiv j}}x^-_{j}(x''_b) :\right)=\hbar^{-\sum_{j=0}^{N-1}a_{ij}v_j''}\left(:\prod_{b\in \la''\atop {c_b\equiv j}}x^-_{j}(x''_b) :\right)\hbar^{h_i}. 
\en
\item[2.] Move $\ds{\prod_{i=0}^{N-1}\left(:\prod_{b\in \la''\atop c_b\equiv i}x^-_i(x''_b):\right)}$ to the left of $\Phi^{(k)}_\emptyset(u_1)$ by using the formula \eqref{commPhiexm}. One gets factors $\ds{
\prod_{b\in \la'' \atop c_b\equiv k}\left(-\frac{\theta(\hbar u_1/x''_b)}{\theta(x''_b/u_1)} \right) }$.
\item[3.] Make the normal ordering 
\be
&&\hspace{-1.5cm} \prod_{i=0}^{N-1}\left(:\prod_{a\in \la' \atop c_a\equiv i}x^-_i(x'_a):\right)\prod_{j=0}^{N-1}\left(:\prod_{b\in \la'' \atop c_b\equiv j}x^-_j(x''_b):\right)\nn\\
&&\hspace{-1.5cm}=\prod_{i=0}^{N-1}\left(:\prod_{a\in \la' \atop c_a\equiv i}x^-_i(x'_a)::\prod_{b\in \la'' \atop c_b\equiv i}x^-_i(x''_b):\right)\times \prod_{a\in \la', b\in \la'' \atop {c_b+1\equiv c_a 
}}\left(-\frac{\theta(t_2 x'_a/x''_b)}{\theta(t_1x''_b/x'_a)}\right)
\prod_{a\in \la', b\in \la'' \atop {c_a+1\equiv c_b 
}}\left(-\frac{\theta(t_1 x'_a/x''_b)}{\theta(t_2x''_b/x'_a)}\right).
\en
\item[4.] Using \eqref{OPExmxm}, make
\be
&&\prod_{i=0}^{N-1}\left(:\prod_{a\in \la' \atop c_a\equiv i}x^-_i(x'_a)::\prod_{b\in \la'' \atop c_a\equiv i}x^-_i(x''_b):\right)\nn\\
&&=\prod_{a\in \la', b\in \la'' \atop c_a\equiv c_b}<x^-_{c_a}(x'_a)x^-_{c_b}(x''_b)>^{Sym}\frac{\theta(x_b''/x_a')}{\theta(\hbar x_a'/x_b'')}
\times \prod_{i=0}^{N-1}\left(:\prod_{a\in (\la',\la'')\atop c_a\equiv i}x^-_i(x_a):\right).
\en
\item[5.] Symmetrize the integrand over $\{x_a\}_{a\in (\la',\la'')}$. 
\end{itemize}
One then obtains 
\be
\Phi^{(k)}_{\la'}(u_1)\Phi^{(k)}_{\la''}(u_2)
&=&\int\prod_{a\in (\la',\la'')}\underline{dx_a}\prod_{i=0}^{N-1}
\left(:\prod_{a\in (\la',\la'') \atop {c_a\equiv i}}x^-_{i}(x_a) 
:\right) \  \Phi^{(k)}_\emptyset(u_1)\  \Phi^{(k)}_\emptyset(u_2)
\nn\\
&& \times\prod_{a,b\in (\la',\la'') \atop {c_a\equiv c_b\atop \rho_a<\rho_b}}{<x^-_{c_a}(x_a)x^-_{c_b}(x_b)>^{Sym}} \ \ \widehat{\Stab}_{\gC,T^{1/2}}((\la',\la'');\bgz).
\en
Here  we set $\{x_a\}_{a\in (\la',\la'')}=\{ x'_a \}_{a\in \la'} \cup \{ x''_b \}_{b\in \la''} $ and 
\bea
&&\hspace{-1.5cm}\widehat{\Stab}_{\gC,T^{1/2}}((\la',\la'');\bgz)\nn\\
&&
\hspace{-1.5cm}
=\Sym_0\Sym_1\cdots \Sym_{N-1} \left(\prod_{a\in \la', b\in \la'' \atop {c_b+1\equiv c_a 
}}\left(-\frac{\theta(t_2 x'_a/x''_b)}{\theta(t_1x''_b/x'_a)}\right)
\prod_{a\in \la', b\in \la'' \atop {c_a+1\equiv c_b 
}}\left(-\frac{\theta(t_1 x'_a/x''_b)}{\theta(t_2x''_b/x'_a)}\right)
\prod_{b\in \la'' \atop {c_b\equiv i 
}}\left(-\frac{\theta(\hbar u_1/x''_b)}{\theta(x''_b/u_1)} \right)
 \right.\nn\\
&&
\hspace{-1.5cm}
\left.\qquad\quad\times \prod_{a\in \la', b\in \la'' \atop {c_a\equiv c_b 
}}\frac{\theta( x''_b/x'_a)}{\theta(\hbar x'_a/x''_b)}
\times  \widehat{\Stab}_{\gC',T^{1/2}}(\la';\{\gz_i\hbar^{\delta_{i,k}-2v''_i+v''_{i-1}+v''_{i+1}}\})
\widehat{\Stab}_{\gC'',T^{1/2}}(\la'';\bgz)\right).
\lb{shufflebasicStab}
\ena
This coincides with the shuffle product formula \eqref{shufflehStabgen} with $\bla'=\la', \bla''=\la''$ and $\bfw'=\bfw''=(\delta_{i,k})$. Hence  
 $\widehat{\Stab}_{\gC,T^{1/2}}((\la',\la'');\bgz)$ gives the elliptic stable envelope for $\E_\rT(\cM(\bfv'+\bfv'',(2\delta_{i,k})))$  with  $\gC: |u_1|>|u_2|, |t_1/t_2|<1$ and  $T^{1/2}$.
We therefore regard $\Phi^{(k)}_{\bla}(u_1,u_2)=\Phi^{(k)}_{\la'}(u_1)\Phi^{(k)}_{\la''}(u_2)$ as the type I  vertex operator for $\cM(\bfv'+\bfv'',(2\delta_{i,k}))$. 

Similarly, consider the composition of the type II dual vertex operators.  
\be
&&\Psi^{*(k)}(u_2)\circ (\Psi^{*(k)}(u_1)\tot \id) \\
&&\qquad  :\ \F^{(1,L)}_{\bxi,\blLa_L,\blLa_M} \tot \F^{(0,-1)(k)}_{u_1}\tot \F^{(0,-1)(k)}_{u_2}\to \F^{(1,L-2)}_{\left(\xi_0,\cdots,(-\hbar^{1/2})^2u_1u_2\xi_k,\cdots,\xi_{N-1},\right),
\blLa_L-2\La_k,\blLa_M+2\La_k},\\
&& \Psi^{*(k)}(u_2)\circ (\Psi^{*(k)}(u_1)\tot \id)(\eta\tot \ket{\la'}^{(k)}_{u_1}\tot \ket{\la''}^{(k)}_{u_2})
=\Psi^{*(k)}_{\la''}(u_2)\Psi^{*(k)}_{\la'}(u_1),\qquad \forall \eta\in \F^{(1,L)}_{\bxi,\blLa_L,\blLa_M}. 
\en
One finds 
\be
\Psi^{*(k)}_{\la''}(u_2)\Psi^{*(k)}_{\la'}(u_1)
&=&\int\prod_{a\in (\la',\la'')}\underline{dx_a}\ 
\widetilde{\Stab}^*_{\gC,T^{1/2}_{opp}}((\la',\la'');\bgz^{*-1})
\prod_{i=0}^{N-1}
\left(:\prod_{a\in (\la',\la'') \atop {c_a\equiv i}}x^+_{i}(x_a) 
:\right) \  
\nn\\
&& \times  \Psi^{*(k)}_\emptyset(u_2)\Psi^{*(k)}_\emptyset(u_1)\ 
\prod_{a,b\in (\la',\la'') \atop {c_a\equiv c_b\atop \rho_a>\rho_b}}{<x^+_{c_a}(x_a)x^+_{c_b}(x_b)>^{Sym}}.
\en
Here  we set $\{x_a\}_{a\in (\la',\la'')}=\{ x'_a \}_{a\in \la'} \cup \{ x''_b \}_{b\in \la''} $ and 
\bea
&&\hspace{-1.5cm}\widetilde{\Stab}_{\gC,T^{1/2}}((\la',\la'');\bgz^{*-1})\nn\\
&&
\hspace{-1.5cm}
=\Sym_0\Sym_1\cdots \Sym_{N-1} \left(\prod_{a\in \la', b\in \la'' \atop {c_b+1\equiv c_a 
}}\left(-\frac{\theta^*(t_2 x'_a/x''_b)}{\theta^*(t_1x''_b/x'_a)}\right)
\prod_{a\in \la', b\in \la'' \atop {c_a+1\equiv c_b 
}}\left(-\frac{\theta^*(t_1 x'_a/x''_b)}{\theta^*(t_2x''_b/x'_a)}\right)
\prod_{a\in \la' \atop {c_a\equiv k
}}\left(-\frac{\theta^*(x'_a/u_2)}{ \theta^*(\hbar u_2/x'_a)}\right)
 \right.\nn\\
&&
\hspace{-1.5cm}
\left.\qquad\quad\times \prod_{a\in \la', b\in \la'' \atop {c_a\equiv c_b 
}}\frac{\theta^*(\hbar x''_b/x'_a)}{\theta^*( x'_a/x''_b)}
\times  \widetilde{\Stab}^*_{\gC',T^{1/2}_{opp}}(\la';\{\gz^{*-1}_i\hbar^{\delta_{i,k}-2v''_i+v''_{i-1}+v''_{i+1}}\})
\widetilde{\Stab}^*_{\gC'',T^{1/2}_{opp}}(\la'';\bgz^{*-1})\right).
\lb{shufflebasicStabs}
\ena
Hence we regard $\Psi^{*(k)}_{\la''}(u_2)\Psi^{*(k)}_{\la'}(u_1)$ as the type II dual vertex operator for $\cM(\bfv'+\bfv'',(2\delta_{i,k}))$.

\subsection{Vertex operators for $\cM(\bfv,\bfw=(w_k\delta_{i,k}))$}\lb{sec:VOkgen}

By repeating the composition of $\Phi^{(k)}(u)$ (resp. $\Psi^{*(k)}(u)$) $w_k$-times, one obtains the type I  vertex operator $\Phi^{(k)}(\bfu)$ (resp. the type II dual vertex operator $\Psi^{*(k)}(\bfu)$)  for $\bigoplus_\bfv \cM(\bfv,\bfw=({w_k}\delta_{i,k})))$. 
Here we set $\bfu=(u_1,\cdots,u_{w_k})$.  Namely, we have    
\be
&&\hspace{-1cm}\Phi^{(k)}(\bfu)=(\id\tot\cdots\tot\id\tot\Phi^{(k)}(u_1))\circ  \cdots\circ (\id\tot\Phi^{(k)}(u_{{w_k}-1}))\circ\Phi^{(k)}(u_{{w_k}})\\
&&\hspace{-0.1cm} :\ \F^{(1,L)}_{\bxi,\blLa_L,\blLa_M}\ \to   \ \F^{(0,-1)(k)}_{u_{{w_k}}}\tot \cdots\tot \F^{(0,-1)(k)}_{u_1}\tot\F^{(1,L+{w_k})}_{\left(\xi_0,\cdots,\frac{\xi_k}{(-\hbar^{1/2})^{{w_k}}u_1\cdots u_{{w_k}}},\cdots,\xi_{N-1}\right), \blLa_L+{w_k}\La_k,\blLa_M},\\[1mm]
&&\hspace{-1cm}\Psi^{*(k)}(\bfu)=\Psi^{*(k)}(u_{{w_k}})\circ(\Psi^{*(k)}(u_{{w_k}-1})\tot\id)\circ \cdots\circ(\Psi^{*(k)}(u_1)\tot\id\tot\cdots\tot\id)  \\
&&\hspace{-0.7cm} :\ \F^{(1,L)}_{\bxi,\blLa_L,\blLa_M}\tot \F^{(0,-1)(k)}_{u_{1}}\tot \cdots\tot \F^{(0,-1)(k)}_{u_{w_k}} \to   \ \F^{(1,L-{w_k})}_{\left(\xi_0,\cdots,{(-\hbar^{1/2})^{{w_k}}u_1\cdots u_{{w_k}}}{\xi_k},\cdots,\xi_{N-1}\right), \blLa_L-{w_k}\La_k,\blLa_M+w_k\La_k}. 
\en
Its components $\Phi^{(k)}_{\bla^{(k)}}(\bfu)$ (resp. $\Psi^{*(k)}_{\bla^{(k)}}(\bfu)$) labeled by  a $w_k$-tuple of partitions $\bla^{(k)}=(\la^{(k)}_1,\cdots,\la^{(k)}_{w_k})
\in (\cP^{(k)})^{{w_k}}$  with color $k$  are defined by 
\bea
&&\Phi^{(k)}(\bfu)=
\sum_{\bla^{(k)}
}\ket{\bla^{(k)}}_{\bfu}
\tot\ \Phi^{(k)}_{\bla^{(k)}}(\bfu) \lb{TypeIMnrtotal}
\ena
and 
\bea
&&\Psi^{*(k)}_{\bla^{(k)}}(\bfu)=\Psi^{*(k)}(\bfu)(\eta\tot \ket{\bla^{(k)}}_{\overline{\bfu}}),\qquad \forall \eta\in \F^{(1,L)}_{\bxi,\blLa_L,\blLa_M}, 
\ena
respectively. Here we set
\be
&&\ket{\bla^{(k)}}_{\bfu}=\ket{\la^{(k)}_{w_k}}_{u_{w_k}}\tot\cdots\tot\ket{\la^{(k)}_1}_{u_1},\\
&&\ket{\bla^{(k)}}_{\overline{\bfu}}=\ket{\la^{(k)}_1}_{u_1}\tot\cdots\tot\ket{\la^{(k)}_{w_k}}_{u_{w_k}}.
\en
Note that the $h_i$-weights of $\ket{\bla^{(k)}}_{\bfu}$ and $\ket{\bla^{(k)}}_{\overline{\bfu}}$ are the same and given by
\be
&&\wt(\bla^{(k)})_i:=\sum_{j=1}^{w_k}\wt(\la^{(k)}_j)_i=|R^{(k)}_i(\bla^{(k)})|-|A^{(k)}_i(\bla^{(k)})|
=-\langle h_i, {w_k}\La_k-\sum_{j\in I}v_j\al_j\rangle 
\en
for
 $\bfv=\{v_i\}_{i\in I}$ with $v_i=\#\{ a\in \blla^{(k)}\ |\ c_a\equiv i\ \mbox{mod}\ N \}$. 
 
One finds 
\bea
\Phi^{(k)}_{\bla^{(k)}}(\bfu)&=&\Phi^{(k)}_{\la^{(k)}_1}(u_1)\cdots \Phi^{(k)}_{\la^{(k)}_{w_k}}(u_{w_k})\nn\\
&=&\int_{\cC^{w_k}} \prod_{a\in \bla^{(k)}}dx_a \prod_{i=0}^{N-1}\left(:\prod_{a\in \bla^{(k)} \atop c_a\equiv i} x^-_i(x_a): \right)\ \Phi^{(k)}_{\emptyset}(u_1)\cdots \Phi^{(k)}_{\emptyset}(u_{w_k})\ \nn\\
&&\times
\prod_{a,b\in \bla^{(k)} \atop {c_a\equiv c_b \atop \rho_a<\rho_b}}<x^-_{c_a}(x_a)x^-_{c_b}(x_b)>^{sym}\ \widehat{\Stab}_{\gC,T^{1/2}}(\bla^{(k)};\bgz).
\lb{typeIVOMnr}
\ena
Here $\widehat{\Stab}_{\gC,T^{1/2}}(\bla^{(k)};\gz)$ is  the elliptic stable envelope 
for $\E_\rT(\cM(\bfv,\bfw=({w_k}\delta_{i,k})))$ obtained by repeating the shuffle product \eqref{shufflebasicStab}. 
We hence regard $\Phi^{(k)}_{\bla}(\bfu)$  as the type I  vertex operator for $\cM(\bfv,({w_k}\delta_{i,k}))$. 
Explicitly,  $\widehat{\Stab}_{\gC,T^{1/2}}(\blla^{(k)};\bgz)$ is given by
\be
&&\widehat{\Stab}_{\gC,T^{1/2}}(\blla^{(k)};\bgz)=
\Sym_0\cdots\Sym_{N-1}\left(\widehat{S}_{\blla^{(k)}}(\bfx,\bfu;p,t_1,t_2)\sum_{\bft\in \Gamma_{\bla^{(k)}}}W^{\bla^{(k)}}_{\bft}(\bfx,\bfu^{(k)},\bgz ;p,t_1,t_2)\right)
\en
with the chamber $\gC : |u_1|> \cdots > |u_{{w_k}}|$, $|t_1/t_2|<1$ 
and 
\be
\widehat{S}_{\blla^{(k)}}(\bfx,\bfu;p,t_1,t_2)&=&\frac{\mbox{$\ds{\prod_{a,b\in\blla^{(k)} \atop {c_a+1\equiv c_b\atop 
\rho_a+1<\rho_b}}\frac{\theta(t_1x_a/x_b)}{\theta(t_2x_b/x_a)}\prod_{a,b\in\blla^{(k)} \atop {c_a+1\equiv c_b\atop 
\rho_a+1>\rho_b}}\frac{\theta(t_2x_b/x_a)}{\theta(t_1x_a/x_b)}
\ \prod_{j=1}^{{w_k}}\hspace{-2mm}\prod_{a\in\blla^{(k)} \atop {c_a\equiv k\atop 
\rho_a> \rho_{r_{k,j}}}}\frac{\theta(\hbar u_j/x_a)}{\theta(x_a/u_j)}}$}
}{\mbox{$\ds{\prod_{a,b\in\blla^{(k)} \atop {c_a\equiv c_b\atop 
\rho_a<\rho_b}}\theta(x_a/x_b)\theta(\hbar x_a/x_b)}$}},
\en
where $\bfx=\{x_a\}_{a\in \bla^{(k)}}$, $\bfu=(u_1,\cdots,u_{{w_k}})$, and  $r_{k,j}$ denotes the root box $(1,1)$ of the partition $\la^{(k)}_j$. 

Similarly, one obtains the  type II dual vertex operator for $\cM(\bfv,({w_k}\delta_{i,k}))$ 
\bea
\Psi^{*(k)}_{\bla^{(k)}}(\bfu)&=&\Psi^{*(k)}_{\la^{(k)}_{w_k}}(u_{w_k})\cdots \Psi^{*(k)}_{\la^{(k)}_1}(u_1)\nn\\
&=&\int_{\cC^{w_k}} \prod_{a\in \bla^{(k)}}dx_a \ \widetilde{\Stab}^*_{\gC,T^{1/2}}(\bla^{(k)};\bgz^{*-1})
\prod_{i=0}^{N-1}\left(:\prod_{a\in \bla^{(k)} \atop c_a\equiv i} x^+_i(x_a): \right)\ \nn\\
&&\times 
\Psi^{*(k)}_{\emptyset}(u_r)\cdots \Psi^{*(k)}_{\emptyset}(u_1)\ \prod_{a,b\in \bla^{(k)} \atop {c_a\equiv c_b \atop \rho_a>\rho_b}}<x^+_{c_a}(x_a)x^+_{c_b}(x_b)>^{sym}.  
\lb{typeIIdualVOMnr}
\ena

\subsection{Vertex operator for general $\cM(\bfv,\bfw)$, 
$\bfv, \bfw\in \N^N$}\lb{sec:VOgen}

 Let us consider  the general $A^{(1)}_{N-1}$ quiver variety $\cM(\bfv,\bfw)$ with $\bfv=(v_0,v_1,\cdots,v_{N-1})$, $\bfw=(w_0,w_1,\cdots,w_{N-1}) \in \N^N$. 
Let $\blla=(\blla^{(0)}, \blla^{(1)},\cdots,\blla^{(N-1)})\in \Lambda_{\bfv,\bfw}$, i.e.
\be
&&\blla^{(k)}=(\la^{(k)}_1,\cdots,\la^{(k)}_{w_k})\ :\ w_k\mbox{-tuple of partitions with  color }\  k\in I,\\
&&v_i=\#\{\ \square \in \blla\ |\ c_\square\equiv i\  \mbox{mod}\ N \}\qquad i\in I.
\en
The Chern roots for the tautological vector bundle $\cV_i$ of rank $v_i$ are denoted  by $\{x_a\}_{a\in I^{(i)}_\blla}=\{x^{(i)}_j, j=1,\cdots,v_i\}
$. We set $\bfx=\cup_{i\in I} \{x_a\}_{a\in I^{(i)}_\blla}$.  
We denote the framing weights by $\bfu=(\bfu^{(0)},\bfu^{(1)},\cdots,\bfu^{(N-1)})$ with $\bu^{(k)}=(u^{(k)}_1,\cdots,u^{(k)}_{w_k})$, $k\in I$. 

Now by composing $\Phi^{(k)}(\bfu^{(k)})$ for all $k\in I$ constructed in the last subsection,  
one obtains the most general vertex operator $\Phi(\bfu)$ for $\bigoplus_\bfv \cM(\bfv,\bfw)$ by  
\be
&&\hspace{-1cm}\Phi(\bfu)=(\id_{w_{N-1}}\tot\cdots\tot\id_{w_1}\tot\Phi^{(0)}(\bfu^{(0)})\circ  \cdots (\id_{w_{N-1}}\tot\Phi^{(N-2)}(\bfu^{(N-2)}))\circ\Phi^{(N-1)}(\bfu^{(N-1)})\\
&&\hspace{1cm} :\ \F^{(1,L)}_{\bxi,\blLa_L,\blLa_M}\ \to   \ \F^{(0,-1)(N-1)}_{\bfu^{(N-1)}}\tot \cdots\tot \F^{(0,-1)(0)}_{\bfu^{(0)}}\tot\F^{(1,L+|\bfw|)}_{\bxi_{(\bfw)}, \blLa_L+\sum_{i\in I}w_i\La_i,\blLa_M}.
\en
where we set $\id_{w_k}=\underbrace{\id\tot\cdots\tot\id}_{w_k}$,\quad  $\F^{(0,-1)(k)}_{\bfu^{(k)}}=\F^{(0,-1)(k)}_{u^{(k)}_1}\tot
\cdots\tot \F^{(0,-1)(k)}_{u^{(k)}_{w_k}}$  and \\
$\bxi_{(\bfw)}=(\xi_{(w_0)},\cdots,\xi_{(w_{N-1})})$ with 
\be
&&\xi_{(w_k)}=\frac{\xi_k}{(-\hbar^{1/2})^{w_k}u^{(k)}_1\cdots u^{(k)}_{w_k}}.
\en
Its components $\Phi_\bla(\bfu)$ labeled by $\bla=(\bla^{(0)},\cdots,\bla^{(N-1)})
\in  (\cP^{(0)})^{w_0}\times \cdots\times (\cP^{(N-1)})^{w_{N-1}}$ 
are defined by 
\bea
&&\Phi(\bfu)=
\sum_{\bla} \ket{\bla}_{\bfu}\tot\ \Phi_{\bla}(\bfu),\lb{TypeIMvwgen}
\ena
where we set
\be
&&\ket{\bla}_{\bfu}=\ket{\bla^{(N-1)}}_{\bfu^{(N-1)}}\tot\cdots\tot\ket{\bla^{(0)}}_{\bfu^{(0)}}
\en
The $h_i$-weight of $\ket{\bla}_{\bfu}$ is given by
\bea
\wt(\bla)_i:&=&\sum_{k\in I}\wt(\bla^{(k)})_i=\sum_{k\in I}\left(|R^{(k)}_i(\bla^{(k)})|-|A^{(k)}_i(\bla^{(k)})|\right)\\
&=&-\langle h_i, \sum_{j\in I}(w_j\La_j-v_j\al_j)\rangle \lb{wtla=wtPhi}
\ena
for
 $\bfv=\{v_i\}_{i\in I}$ with $v_i=\#\{ a\in \bla \ |\ c_a\equiv i\ \mbox{mod}\ N \}$. 

One finds
\bea
\Phi_{\blla}(\bfu)
&=&\Phi^{(0)}_{\blla^{(0)}}(\bfu^{(0)})\Phi^{(1)}_{\blla^{(1)}}(\bfu^{(1)})\cdots \Phi^{(N-1)}_{\blla^{(N-1)}}(\bfu^{(N-1)})\nn\\
&=&\int\prod_{a\in \blla}\underline{dx_a}
\prod_{i=0}^{N-1}
\left(:\prod_{a\in \blla \atop c_a\equiv i}x^-_{i}(x_a) :\right) \  
\Phi_{\blempty}(\bfu)
\nn\\
&&\qquad\times\prod_{a,b\in \blla\atop {c_a\equiv c_b\atop \rho_a<\rho_b}}{<x^-_{c_a}(x_a)x^-_{c_b}(x_b)>^{Sym}}\ \widehat{\Stab}_{\gC,T^{1/2}}(\blla;\bgz)\lb{VOgencyclic}
\ena
with
\bea
\Phi_{\blempty}(\bfu)&=&\Phi^{(0)}_\emptyset(\bfu^{(0)})\Phi^{(1)}_\emptyset(\bfu^{(1)})\cdots  \Phi^{(N-1)}_\emptyset(\bfu^{(N-1)}).
\ena
Here $\widehat{\Stab}_{\gC,T^{1/2}}(\blla;\bgz)$ is the elliptic stable envelope for $\E_\rT(\cM(\bfv,\bfw))$ \eqref{def:hatStab} obtained by using the shuffle product formula \eqref{shufflehStabgen} repeatedly. 

Similarly, the type II dual vertex operator for $\cM(\bfv,\bfw)$ is given by the composition 
\bea
&&\hspace{-1cm}\Psi^{*}(\bfu)=\Psi^{*(N-1)}(\bfu^{(N-1)})\circ(\Psi^{*(N-2)}(\bfu^{(N-2)})\tot\id_{w_{N-1}})\circ \cdots\circ(\Psi^{*(0)}(\bfu^{(0)})\tot\id_{w_1}\tot\cdots\tot\id_{w_{N-1}}) \nn \\
&&\hspace{1cm} :\ \F^{(1,L)}_{\bxi,\blLa_L,\blLa_M}\tot \F^{(0,-1)(0)}_{\overline{\bfu^{(0)}}}\tot \cdots\tot \F^{(0,-1)(N-1)}_{\overline{\bfu^{(N-1)}}} \to   \ \F^{(1,L-|\bfw|)}_{\bxi^{(\bfw)}, \blLa_L-\sum_{i\in I}w_i\La_i,\blLa_M+\sum_{i\in I}w_i\La_i},\nn\\
&&
\ena
where we set $\ds{\F^{(0,-1)(k)}_{\overline{\bfu^{(k)}}}=\F^{(0,-1)(k)}_{u^{(k)}_{w_k}}\tot
\cdots\tot \F^{(0,-1)(k)}_{u^{(k)}_{1}}}$ and $\bxi^{(\bfw)}=(\xi^{(w_0)},\cdots,\xi^{(w_{N-1})})$ with 
\bea
&&\xi^{(w_k)}=(-\hbar^{1/2})^{w_k}u^{(k)}_1\cdots u^{(k)}_{w_k}{\xi_k}.
\ena
Its components are defined by
\be
&&\Psi^{*}_{\blla}(\bfu)\eta
=\Psi^{*}(\bfu)\left(\eta\tot \ket{\bla}_{\overline{\bfu}}
\right)\qquad \qquad \forall \eta \in \F^{(1,L)}_{\bxi,\blLa_L,\blLa_M}
\en
with 
\be
&&\ket{\bla}_{\overline{\bfu}}=
\ket{\bla^{(0)}}_{\overline{\bfu^{(0)}}}\tot\cdots\tot\ket{\bla^{(N-1)}}_{\overline{\bfu^{(N-1)}}}.
\en
One finds 
\bea
\Psi^{*}_{\blla}(\bfu)
&=& \Psi^{*(N-1)}_{\blla^{(N-1)}}(\bfu^{(N-1)})
\cdots \Psi^{*(1)}_{\blla^{(1)}}(\bfu^{(1)})\Psi^{*(0)}_{\blla^{(0)}}(\bfu^{(0)})
\nn\\
&=&\int\prod_{a\in \blla}\underline{dx_a}\ 
\widetilde{\Stab}^*_{\gC,T^{1/2}_{opp}}(\blla;\bgz^{*-1})
\prod_{i=0}^{N-1}
\left(:\prod_{a\in \blla \atop c_a\equiv i}x^+_{i}(x_a) :\right) \  
\Psi^*_{\blempty}(\bfu)
\nn\\
&&\qquad\times\prod_{a,b\in \blla\atop {c_a\equiv c_b\atop \rho_a>\rho_b}}{<x^+_{c_a}(x_a)x^+_{c_b}(x_b)>^{Sym}},\ \lb{TypeIIVOgencyclic}\\
\Psi^*_{\blempty}(\bfu)&=&\Psi^{*(N-1)}_\emptyset(\bfu^{(N-1)})\cdots \Psi^{*(1)}_\emptyset(\bfu^{(1)}) \Psi^{*(0)}_\emptyset(\bfu^{(0)}).
\ena

The results in this section and Conjecture \ref{Nvw} lead us to the following conjecture, whose statement becomes closer to Example 1) in Sec.\ref{sec:intro}. 
\begin{conj}\lb{level0mw2Mvw}
The tensor product representation $\F^{(0,-1)(N-1)}_{\bfu^{(N-1)}}\tot \cdots \tot \F^{(0,-1)(0)}_{\bfu^{(0)}}$ of\\  $U_{t_1,t_2,p}(\gl_{N,tor})$ is equivalent to an expected  level-$(0,-|\bfw|)$  geometric action of the same algebra on  $\bigoplus_{\bold{v}}{\mathrm E}_\rT({\cal M}({\bold{v},\bfw}))$ under the identification of the bases $\ket{\bla}_\bfu$ with the $\rT$-fixed point classes $[\bla]$ of ${\cal M}(\bold{v},\bfw)$.   
\end{conj}

In \cite{KS}, we studied a similar construction of the vertex operators of the elliptic quantum toroidal algebra  $U_{t_1,t_2,p}(\gl_{1,tor})$ by using the elliptic stable envelope for the equivariant elliptic cohomology $\E_\rT(\cM(n,r))$ of the instanton moduli space $\cM(n,r)$. See Example 2) in Sec.\ref{sec:intro} in this paper and Sec.4 of \cite{KS}. We hence have the following similar conjecture.
\begin{conj}
The tensor product $\F^{(0,-1)}_{u_r}\tot \cdots \tot \F^{(0,-1)}_{u_1}$ of  the level-$(0,-1)$ representation $\F^{(0,-1)}_{u}$ of $U_{t_1,t_2,p}(\gl_{1,tor})$ is equivalent to an expected  level-$(0,-r)$  geometric action of the same algebra  on  $\bigoplus_{n}{\mathrm E}_\rT({\cal M}(n,r))$ under the identification of the bases $\ket{\la^{(r)}}_{u_r}\tot \cdots \tot\ket{\la^{(1)}}_{u_1}$ with the $\rT$-fixed point classes $[(\la^{(r)}, \cdots, \la^{(1)})]$ of ${\cal M}(n,r)$.   
\end{conj}

\section{The $\rK$-Theoretic Vertex Functions for $\cM(\bfv,\bfw)$}\lb{sec:VFs}

We show that the highest to highest expectation value of the vertex operator $\Phi_\bla(\bfu)$ of 
$U_{t_1,t_2,p}(\gl_{N,tor})$ for   $X=\cM(\bfv,\bfw)$ constructed in the last subsection gives the $\rK$-theoretic  vertex function for the quasi-maps $\PP^1\dashrightarrow X$\cite{AO,AFO,Ok}.

\begin{lem}\lb{OPEs}
The following OPEs are hold. 
\be
(1)&&\prod_{k=0}^{N-1}\prod_{j=1}^{w_k}\Phi^{(k)}_{\emptyset}(u^{(k)}_j)=\mu(\bfu)
:\prod_{k=0}^{N-1}\prod_{j=1}^{w_k}\Phi^{(k)}_{\emptyset}(u^{(k)}_j):,\\
&&\mu(\bfu)=\prod_{0\leq k\leq l\leq N-1}\prod_{i=1}^{w_k}\prod_{j=1}^{w_l}
\left(-\hbar^{1/2}u^{(k)}_i\right)^{\eta_{kl}}
\frac{(t_1^Nt_1^{k-l}u^{(l)}_j/u^{(k)}_i;p,t_1^N)_\infty}{(t_1^Nt_1^{k-l}u^{(l)}_j/u^{(k)}_i;\hbar,t_1^N)_\infty}
\frac{(t_2^{l-k}u^{(l)}_j/u^{(k)}_i;p,t_2^N)_\infty}{(t_2^{l-k}u^{(l)}_j/u^{(k)}_i;\hbar,t_2^N)_\infty},
\en
for $|u^{(l)}_j/u^{(k)}_i|<1$, 
\be
(2)&&\prod_{i=0}^{N-1}:\prod_{a\in \blla \atop c_a\equiv i}x^-_{i}(x_a) :
=\prod_{a,b\in \blla\atop c_a+1\equiv c_b}x_a^{-1}
\frac{(p\hbar^{-1}t_2x_b/x_a;p)_\infty}{(t_2x_b/x_a;p)_\infty}:\prod_{i=0}^{N-1}\prod_{a\in \blla \atop c_a\equiv i}x^-_{i}(x_a) :,\qquad |x_b/x_a|<1,\\
(3)&&:\prod_{i=0}^{N-1}\prod_{a\in \blla \atop c_a\equiv i}x^-_{i}(x_a) :
 :\prod_{k=0}^{N-1}\prod_{j=1}^{w_k}\Phi^{(k)}_{\emptyset}(u^{(k)}_j):\nn\\
 &&=\prod_{k=0}^{N-1}\prod_{j=1}^{w_k}\prod_{a\in \blla\atop c_a\equiv k} x_a\frac{(p u^{(k)}_j/x_a;p)_\infty}{(\hbar u^{(k)}_j/x_a;p)_\infty}
:\prod_{i=0}^{N-1}\prod_{a\in \blla \atop c_a\equiv i}x^-_{i}(x_a) 
\prod_{k=0}^{N-1}\prod_{j=1}^{w_k}\Phi^{(k)}_{\emptyset}(u^{(k)}_j):, \qquad |u^{(k)}_j/x_a|<1.
\en
\end{lem}

For a torus fixed point $\blmu \in \Lambda_{\bfv,\bfw}$ of $\cM(\bfv,\bfw)$, 
we set $\varphi^\blmu_a=t_1^{-(y-1)}t_2^{-(x-1)}$ for a box $a=(x,y)\in \blmu$ as before. Noting the remark below \eqref{IdGeoRep}, we take the Jackson integral 
\be
&&\int_0^{\varphi^\blmu_{a}}d_px_a f(x_a)=(1-p)\varphi^\blmu_{a}\sum_{d_a\in \N} f(\varphi^\blmu_{a} p^{d_a})p^{d_a}
\en
 in the vertex operators \eqref{VOgencyclic} for $\cM(\bfv,\bfw)$. We denote thus obtained vertex operator by $\Phi^{(\bmu)}_\bla(\bfu)$. 
For  the highest weight vector $\ket{0}^{(1,L)}_{\bxi,\blLa_L,\blLa_M}$ in $\F^{(1,L)}_{\bxi,\blLa_L,\blLa_M}$,  
 let  ${}^{\qquad(1,L)}_{\ \ \bxi,\blLa_L,\blLa_M}\bra{0}$ be the dual vector satisfying 
\be
&&{}^{\qquad(1,L)}_{\ \ \bxi,\blLa_L,\blLa_M}\bra{0}\ket{0}^{(1,L)}_{\bxi,\blLa_L,\blLa_M}=1.
\en

We then obtain the following expectation value. 
\begin{prop}\lb{VEVPhi}
\bea
&&\hspace{-0.5cm}{}_{\quad \bxi_{(\bfw)}, \blLa'_L,\blLa_M}^{(1,L+|\bfw|)}\bra{0}\Phi^{(\bmu)}_{\blla}(\bfu)\ket{0}^{(1,L)}_{\bxi,\blLa_L,\blLa_M}
\nn\\&&
=\mu(\bfu)\prod_{a\in \blla}\left(\int_0^{\varphi^{\blmu}_{a}}\hspace{-0.1cm}{d_px_a}\right)\ 
\prod_{k=0}^{N-1}\prod_{j=1}^{w_k}\prod_{a\in \blla\atop c_a\equiv k} x_a\frac{(p u^{(k)}_j/x_a;p)_\infty}{(\hbar u^{(k)}_j/x_a;p)_\infty}
 \prod_{a,b\in \blla\atop c_a+1\equiv c_b}\hspace{-0.2cm}x_a^{-1}
\frac{(p\hbar^{-1}t_2x_b/x_a;p)_\infty}{(t_2x_b/x_a;p)_\infty}\nn\\
&&\hspace{1cm}\qquad\times
\prod_{a,b\in \blla\atop {
c_a\equiv c_b \atop \rho_a\not= \rho_b}}<x^-_{c_a}(x_a)x^-_{c_b}(x_b)>^{Sym}
\ \widehat{\Stab}_{\gC,T^{1/2}}(\blla;\bgz).\lb{HtHEPhi}
\ena
Here we set
\be
&&\blLa'_L=\blLa_L+\sum_{i\in I}\left(w_i\La_i-v_i\al_i\right).
\en
\end{prop}
We show that the normalized expectation value ${}_{\quad \bxi_{(\bfw)}, \blLa'_L,\blLa_M}^{\qquad(1,L+|\bfw|)}\bra{0}\Phi^{(\bmu)}_{\blla}(\bfu)\ket{0}^{(1,L)}_{\bxi,\blLa_L,\blLa_M}/\cN_{\blmu}$ 
gives a vertex function for $X=\cM(\bfv,\bfw)$\cite{AFO,AOMosc,Ok}. Here the factor  $\cN_{\blmu}$ is the restriction of the integrand of 
${}_{\quad \bxi_{(\bfw)}, \blLa'_L,\blLa_M}^{\qquad(1,L+|\bfw|)}\bra{0}\Phi^{(\bmu)}_{\blla}(\bfu)\ket{0}^{(1,L)}_{\bxi,\blLa_L,\blLa_M}$ without $\widehat{\Stab}_{\gC,T^{1/2}}(\blla;\bgz)$ to the fixed point $\bmu$\cite{KPSZ,KS,Konno24} given by 
\be
\cN_\blmu
&=&\mu(\bfu)\left[\prod_{k=0}^{N-1}\prod_{j=1}^{w_k}\prod_{a\in \blla\atop c_a\equiv k} x_a\frac{(p u^{(k)}_j/x_a;p)_\infty}{(\hbar u^{(k)}_j/x_a;p)_\infty}
 \prod_{a,b\in \blla\atop c_a+1\equiv c_b}x_a^{-1}
\frac{(p\hbar^{-1}t_2 x_b/x_a;p)_\infty}{(t_2 x_b/x_a;p)_\infty}\right.\nn\\
&&\qquad\qquad\qquad\qquad\left.\times
\prod_{a,b\in \blla\atop {
c_a\equiv c_b \atop \rho_a\not= \rho_b}}\frac{x_ax_b}{\hbar^{1/2}}
\frac{(\hbar  x_a/x_b;p)_\infty}{(p x_a/x_b;p)_\infty} \right]
_{x_a=\varphi^\blmu_{a}\atop a\in \bla}.
\en
To evaluate the Jackson integrals, we use the following quasi-periodicity of $\widehat{\Stab}_{\gC,T^{1/2}}(\blla;\bgz)$ obtained by the same argument as Proposition 5.1 of \cite{KS}. 
\begin{prop}\lb{QPStab}
For the elliptic stable envelope $\widehat{\Stab}_{\gC,T^{1/2}}(\blla;\bgz)$ for $\cM(\bfv,\bfw)$, one has
\be
&& \widehat{\Stab}_{\gC,T^{1/2}}(\blla;\bgz)\bigr|_{x_a=\varphi^\blmu_{a}p^{d_a}\atop  a\in \bla}=\prod_{i\in I}\gz_i^{-\sum_{a\in \bla \atop c_a\equiv i}d_a}\times \widehat{\Stab}_{\gC,T^{1/2}}(\bla;\bgz)\bigr|_{\blmu},  
\en
where $\widehat{\Stab}_{\gC,T^{1/2}}(\bla;\bgz)\bigr|_{\blmu}$ stands for the elliptic stable envelope restricted to the fixed point $\blmu$ i.e. 
\be
&& \widehat{\Stab}_{\gC,T^{1/2}}(\bla;\bgz)\bigr|_{\blmu}:=\widehat{\Stab}_{\gC,T^{1/2}}(\bla;\bgz)\bigr|_{x_a=\varphi^\blmu_{a}
\atop  a\in \bla },
\en
and is independent from ${\bf d}=(d_a)_{a\in \bla}$. 
\end{prop}

\medskip
\noindent
{\it Remark.}\ Let us consider the following function 
\be
&&{\bf e}({\bf x}, { \bgz})=\exp\left\{-\frac{1}{\log p}\sum_{i\in I}\sum_{a\in \bla\atop c_a\equiv i} \log x_a \log \gz_i\right\}=\prod_{i\in I}\prod_{j=1}^{v_i}(x^{(i)}_j)^{-\frac{P_i+h_i}{r}},
\en
where we set 
 $p=\hbar^r$, $r\in \C^\times$ with $|p|<1$. 
This has the same transformation property as  $\widehat{\Stab}_{\gC,T^{1/2}}(\bla;\bgz)$ in Proposition \ref{QPStab}.

We hence obtain the normalized expectation value
\bea
V_{\blla,\blmu}(\bfu)&=&\frac{1}{\cN_\blmu}\times{}_{\quad \bxi_{(\bfw)}, \blLa'_L,\blLa_M}^{(1,L+|\bfw|)}\bra{0}\Phi^{(\bmu)}_{\blla}(\bfu)\ket{0}^{(1,L)}_{\bxi,\blLa_L,\blLa_M}
\nn\\
&=&\sum_{\bfd\in \N^{|\bfv|}} \prod_{k=0}^{N-1}
\left(\hbar^{w_k}p^{2-2v_k+v_{k+1}-2w_k}\gz_k\right)^{-\sum_{a\in \blmu \atop c_a\equiv k}d_a}
\prod_{k=0}^{N-1}
\prod_{j=1}^{w_k}
\prod_{a\in \blla\atop c_a\equiv k} 
\frac{( \varphi^\blmu_{a}/u^{(k)}_j;p)_{-d_a}}{(p\hbar^{-1} \varphi^\blmu_{a}/u^{(k)}_j;p)_{-d_a}}
\nn\\
&&\qquad\qquad\times 
\prod_{a,b\in \blla\atop c_a+1\equiv c_b}
\frac{(t_2\varphi^\blmu_{b}/\varphi^\blmu_{a};p)_{d_b-d_a}}
{(p\hbar^{-1}t_2\varphi^\blmu_{b}/\varphi^\blmu_{a};p)_{d_b-d_a}}
\prod_{a,b\in \blla\atop {
c_a\equiv c_b \atop \rho_a\not= \rho_b}}
\frac{(p \varphi^\blmu_{a}/\varphi^\blmu_{b};p)_{d_a-d_b}}
{(\hbar \varphi^\blmu_{a}/\varphi^\blmu_{b};p)_{d_a-d_b}}\nn\\
&&\qquad\qquad\qquad\qquad\times \widehat{\Stab}_{\gC,T^{1/2}}(\blla;\bgz)\Bigr|_{x_a=\varphi^\bmu_{a}\atop a\in \bla}\nn\\
&=&\sum_{\bfd\in \N^{|\bfv|}} 
\prod_{k=0}^{N-1}
\left(\hbar^{w_k}p^{2-2v_k+v_{k+1}-2w_k}\gz_k\right)^{-\sum_{a\in \blmu \atop c_a\equiv k}d_a}
\prod_{k=0}^{N-1}\prod_{j=1}^{w_k}
\prod_{l=1}^{v_k}
\frac{(\varphi^\blmu_{a^{(k)}_l}/u^{(k)}_j;p)_{d_{k,l}}}{(p\hbar^{-1}
\varphi^\blmu_{a^{(k)}_l}/ u^{(k)}_j;p)_{d_{k,l}}}
\nn\\
&&\hspace{-1cm}\times 
\prod_{k=0}^{N-1}\prod_{j=1}^{v_k}\prod_{l=1}^{v_{k+1}}
\frac{(t_2 \varphi^\blmu_{a^{(k+1)}_l}
/\varphi^\blmu_{a^{(k)}_j};p)_{d_{k+1,l}-d_{k,j}}}
{(p\hbar^{-1}t_2\varphi^\blmu_{a^{(k+1)}_l}
/\varphi^\blmu_{a^{(k)}_j};p)_{d_{k+1.l}-d_{k,j}}}
\prod_{k=0}^{N-1}\prod_{j,l=1\atop j\not= l}^{v_k}
\frac{(p \varphi^\blmu_{a^{(k)}_l}
/\varphi^\blmu_{a^{(k)}_j};p)_{d_{k,l}-d_{k,j}}}
{(\hbar  \varphi^\blmu_{a^{(k)}_l}
/\varphi^\blmu_{a^{(k)}_j};p)_{d_{k,l}-d_{k,j}}}\nn\\
&&\qquad\qquad\qquad\qquad\times \widehat{\Stab}_{\gC,T^{1/2}}(\blla;\bgz)\Bigr|_{x_a=\varphi^\bmu_{a}\atop a\in \bla}.\lb{VF}
\ena
Here $a^{(k)}_j\in \blmu$ are given in Sec.\ref{tautVB} and we set $d_{k,j}=d_{a^{(k)}_j}$. We also   used the formula for $d>0$
\be
&&(z;p)_{-d}=(-z)^{-d}p^{\frac{d(d-1)}{2}}\frac{1}{(pz^{-1};p)_d}.
\en 
Hence by the shift $u^{(k)}_j \ \mapsto\ \hbar^{-1}u^{(k)}_j $ for all $k, j$, the function $V_{\bla,\blmu}(\hbar^{-1}\bfu)$ gives a natural higher rank extension of the $\rK$-theoretic vertex functions obtained in \cite{KS,SD20,DinkinsThesis}. 

In order to confirm this, let us consider a map $\Phi : \rK_T(X)\to \rK_T(X)\otimes \Z[[p]]$ with $|p|<1$ defined by  
\be
&&\Phi(\sum_{i}a_i-\sum_{j}b_j)=\frac{\prod_i(a_i;p)_\infty}{\prod_j(b_j;p)_\infty},\quad \sum_{i}a_i-\sum_{j}b_j\in \rK_T(X). 
\en
One then finds that the polarization 
\be
T^{1/2}&=&\sum_{i=0}^{N-1}\cW_i\otimes \cV_i^*+ t_1^{-1}\sum_{i=0}^{N-1} \cV_{i+1}\otimes \cV_i^*-
\sum_{i=0}^{N-1} \cV_{i}\otimes \cV_i^*\\
&=&\sum_{i=0}^{N-1}\sum_{j=1}^{w_i}\sum_{k=1}^{v_i}\frac{u^{(i)}_j}{x^{(i)}_k}+ t_1^{-1}\sum_{i=0}^{N-1} \sum_{j=1}^{v_{i+1}}\sum_{k=1}^{v_{i}}\frac{x^{(i+1)}_j}{x^{(i)}_k}-
\sum_{i=0}^{N-1}  \sum_{j,k=1\atop j\not=k}^{v_i}\frac{x^{(i)}_j}{x^{(i)}_k} \quad \in  \rK_T(X)
\en
gives 
\be
\Phi((p-\hbar)T^{1/2})&=&\prod_{i=0}^{N-1}\prod_{j=1}^{w_i}\prod_{k=1}^{v_i}\frac{(pu^{(i)}_j/x^{(i)}_k;p)_\infty}{(\hbar u^{(i)}_j/x^{(i)}_k;p)_\infty} \prod_{i=0}^{N-1} \prod_{j=1}^{v_{i+1}}\prod_{k=1}^{v_{i}}\frac{(pt_1^{-1}x^{(i+1)}_j/x^{(i)}_k;p)_\infty}{(\hbar t_1^{-1}x^{(i+1)}_j/x^{(i)}_k;p)_\infty}\\
&&\times\prod_{i=0}^{N-1}  \prod_{j,k=1\atop j\not=k}^{v_i}\frac{(\hbar x^{(i)}_j/x^{(i)}_k;p)_{\infty}}{(p x^{(i)}_j/x^{(i)}_k;p)_{\infty}}. 
\en
Remembering $x_{a^{(i)}_j}=x^{(i)}_j$, this coincides with the integrand of \eqref{HtHEPhi} up to a factor which drops after normalization. 
Combining this with the above Remark, \eqref{HtHEPhi} coincides with the integral formula for a vertex function for $X$ derived in \cite{AFO} Theorem 3.1. Here our $\widehat{\Stab}_{\gC,T^{1/2}}(\blla;\bgz)$ corresponds to $\mathscr{F}'(s)$ in \cite{AFO} and plays a role of a pole subtraction matrix\cite{AO}.  

We emphasize that our representation theoretical derivation of the $\rK$-theoretic vertex functions 
 given here as well as  in \cite{KS,Konno24} is relying on the {\it elliptic} nature of EQG/EQTA,  $U_{q,p}(\slnh)$, $U_{t_1,t_2,p}(\gl_{1,tor})$ and $U_{t_1,t_2,p}(\gl_{N,tor})$. 
Namely the parameter $p$ appearing in the $\rK$-theoretic counting of the quasimap is nothing but  the elliptic nome of  EQG/EQTA. It is also interesting to note that in \cite{AFO} it is identified with a deformation parameter $q$ of the deformed $W$-algebra $\cW_{q,t}(\g)$\cite{FrRe}. This is consistent to 
the relationship between EQG/EQTA and deformed $W$-algebras including the affine quiver ones 
under the identification of the parameters $(p,p^*)$ with $(q,t)$. See introduction section and \cite{KOgl1,Konno24}. In the current case of $U_{t_1,t_2,p}(\gl_{N,tor})$, the corresponding deformed $W$-algebra should be the $A^{(1)}_{N-1}$ type affine quiver $W$-algebra\cite{Kim25}. We expect  the same relationship as $U_{t_1,t_2,p}(\gl_{1,tor})$ realizes the Jordan quiver $W$-algebra\cite{KOgl1}.

\subsection{Bethe ansatz equatioin}

Let us consider the limit $p\to1$ of the integral in Proposition \ref{VEVPhi}. We have the following conjugate modulus transformation.
\begin{prop}
\bea
\widehat{\Stab}_{\gC,T^{1/2}}(\bla;\bgz)&=&\exp\left\{-\frac{1}{\log p}\sum_{i\in I}\sum_{a\in \bla\atop c_a\equiv i}\log x_a\log \gz_i\right\}\exp\left\{\frac{1}{\log p}\sum_{i\in I}\sum_{j=1}^{w_i}\log u^{(i)}_j
\sum_{a\in \la^{(i)}_j}\log\gz_{c_a}
\right\}\nn\\
&&\times \widehat{\Stab}_{\gC,T^{1/2}}(\bla;\bgz)\biggl|_{\theta\mapsto \vartheta_1}.\lb{conjMT}
\ena
Here $\theta \mapsto \vartheta_1$ denotes the replacement of all theta-function-factors $\theta(X)$ appearing in 
the expression of $\widehat{\Stab}_{\gC,T^{1/2}}(\bla;\bgz)$ by 
\be
&&\vartheta_1\left(\left.\frac{\log X}{\log p}\right| \tau
\right).
\en
Here $\vartheta_1(v |\tau)$ is given by
\be
&&\vartheta_1(v |\tau)=\sqrt{-1}\sum_{n\in \Z}(-1)^ne^{\pi \sqrt{-1}\tau(n-1/2)^2}e^{2\pi \sqrt{-1}v(n-1/2)}
\en
with $\tau=-2\pi i/\log p$. 
\end{prop}

\noindent
{\it Proof.}\ Use the formula
\be
&&\theta(X)=e^{-\frac{\pi i}{4}}\tau^{1/2}p^{-\frac{1}{8}}\hbar^{-\frac{(\log X)^2}{2\log p}}\vartheta_1\left(\left.\frac{\log X}{\log p}\right|\tau
\right).
\en
\qed

\begin{prop}
In the limit $p\to 1$, the integrand is dominated by the saddle point satisfying the following equations.
\be
&&\hspace{-1cm}\prod_{j=1}^{w_k}\frac{1-u^{(k)}_j/x^{(k)}_i}{1-\hbar u^{(k)}_j/x^{(k)}_i}
\prod_{l=1}^{v_{k+1}}\frac{1-t_1^{-1}x^{(k+1)}_l/x^{(k)}_i}{1-t_2 x^{(k+1)}_l/x^{(k)}_i}
\prod_{m=1}^{v_{k-1}}\frac{1-t_2 x^{(k)}_i/x^{(k-1)}_m}{1-t_1^{-1}x^{(k)}_i/x^{(k-1)}_m}
\prod_{n=1\atop \not=i}^{v_{k}} \frac{1-\hbar x^{(k)}_n/x^{(k)}_i}{1-\hbar^{-1}  x^{(k)}_n/x^{(k)}_i}
=\gz_k\hbar^{v_k-1}
\en
for $k\in I$ and $i\in\{1,\cdots,v_k\}$. Here  $x^{(k)}_i$ is defined in Sec.\ref{tautVB}.  
\end{prop}

\noindent
{\it Proof.}\ In the limit $p\to 1$, the integrand decays exponentially as 
\be
&&\exp\left\{-\frac{1}{1-p}W(\bfx,\bgz)\right\}.
\en
Note that $W(\bfx,\bgz)$ contains a contribution from the first factor in the RHS of \eqref{conjMT} due to 
$1/\log p\ \sim \ -1/(1-p)$ as $p\to 1$.  Then the required equations are the conditions 
\be
&&\frac{\partial W(\bfx,\bgz)}{\partial\ x^{(k)}_i}=0. 
\en  
\qed

\vspace{2mm}
\noindent
{\it Remark.}\ In our previous paper \cite{KS}, we derived the $\rK$-theoretic vertex function 
for the ADHM instanton moduli space $\cM(n,r)$. In the same way as above, one can derive 
the saddle point condition in the limit $p\to 1$ as
\bea
&&\hspace{-1cm}
\prod_{j=1}^{r}\frac{1-u_j/x_a}{1-\hbar u_j/x_a}
\prod_{b\in \bla \atop \not=a} \frac{(x_a-\hbar^{-1} x_b)(x_a-t_1  x_b)(x_a-t_2 x_b)}
{(x_a-\hbar x_b)(x_a-t_1^{-1} x_b)(x_a-t_2^{-1}  x_b)}
=\gz
\ena  
for $a\in \bla$. Here  $\bla$ is a $r$-tuple of partitions $(\la^{(1)},\cdots,\la^{(r)})$ satisfying $|\bla|=n$. 
This coincides with  (6.1) in \cite{FJMM}.

\section{Exchange Relations of the Vertex Operators} \lb{sec:ExchRs}

The exchange relations among the vertex operators $\Phi_\blla(\bfu)$ and $\Psi^{*}_\blla(\bfu)$ 
 can be derived  in the same way as the $\gl_{1,tor}$ case\cite{KS}. 

Let $\blla'\in \Lambda_{\bfv',\bfw'}$, $\blla''\in \Lambda_{\bfv'',\bfw''}$ and 
\be
&&\bfu_1=(u^{(k)'}_i)_{k\in I, 1\leq i\leq w'_k},\qquad \bfu_2=(u^{(l)''}_j)_{l\in I, 1\leq j\leq w''_l}. 
\en
Consider the vertex operators $\Phi_{\bla'}(\bfu_1)$ and  $\Phi_{\bla''}(\bfu_2)$ for $\cM(\bfv',\bfw')$ and $\cM(\bfv'',\bfw'')$, respectively. 
The composition $\Phi_{\bla'}(\bfu_1)\Phi_{\bla''}(\bfu_2)$  gives the vertex operator for $\cM(\bfv,\bfw)$ with $\bfv=\bfv'+\bfv''$ and  $\bfw=\bfw'+\bfw''$. 
\bea
&&\hspace{-0.7cm}\Phi^{(\blmu')}_{\bla'}(\bfu_1)\Phi^{(\blmu'')}_{\bla''}(\bfu_2)=\prod_{a\in \bla'}\left(\int_0^{\varphi^{\bmu'}_{\iota(a)}}d_px'_a\right)
\prod_{b\in \bla''}\left(\int_0^{\varphi^{\bmu''}_{\iota(b)}}d_px''_b\right)
\ \prod_{i=0}^{N-1}\left(:\prod_{a\in (\bla',\bla'')\atop c_a\equiv i}x_i^-(x_a): \right)\nn\\
&&\qquad\qquad\qquad\times \Phi_{\blempty}(\bfu_1)\Phi_{\blempty}(\bfu_2)
\prod_{a, b\in (\bla',\bla'')\atop {c_a\equiv c_b\atop \rho_a<\rho_b}}<x^-_{c_a}(x_a)x^-_{c_b}(x_b)>^{sym}\ \widehat{\Stab}_{\gC,T^{1/2}}((\bla',\bla'');\bgz),\nn\\
&&\lb{Phi1Phi2}
\ena
where the chamber $\gC$ is given by $|u^{(0)'}_1| >\cdots>|u^{(N-1)'}_{w'_{N-1}}|>|u^{(0)''}_1| >\cdots>|u^{(N-1)''}_{w''_{N-1}}|$ with the stability condition $|t_1/t_2|<1$.
Similarly we have
\bea
&&\hspace{-0.7cm}\Phi^{(\blmu'')}_{\bla''}(\bfu_2)\Phi^{(\blmu')}_{\bla'}(\bfu_1)=
\prod_{a\in \bla'}\left(\int_0^{\varphi^{\bmu'}_{\iota(a)}}d_px'_a\right)
\prod_{b\in \bla''}\left(\int_0^{\varphi^{\bmu''}_{\iota(b)}}d_px''_b\right)
\ \prod_{i=0}^{N-1}\left(:\prod_{a\in (\bla',\bla'')\atop c_a\equiv i}x_i^-(x_a): \right)\nn\\
&&\qquad\qquad\qquad \times
\Phi_{\blempty}(\bfu_2)\Phi_{\blempty}(\bfu_1)
  \prod_{a, b\in (\bla',\bla'')\atop {c_a\equiv c_b\atop \rho_a<\rho_b}}<x^-_{c_a}(x_a)x^-_{c_b}(x_b)>^{sym}
\ \widehat{\Stab}_{\gC,T^{1/2}}((\bla'',\bla');\bgz),\nn\\
&&\lb{Phi2Phi1}
\ena
where the chamber $\overline{\gC}$ is given by $|u^{(0)''}_1| >\cdots>|u^{(N-1)''}_{w''_{N-1}}|>|u^{(0)'}_1| >\cdots>|u^{(N-1)'}_{w'_{N-1}}|$ with the stability condition $|t_1/t_2|<1$.

Let  
\be
&&\hspace{-0.5cm}\F^{(0,-1)}_{\bfu_1}:=\F^{(0,-1)(N-1)}_{\bfu^{(N-1)'}}\tot\cdots \tot \F^{(0,-1)(0)}_{\bfu^{(0)'}},\qquad 
\F^{(0,-1)}_{\bfu_2}:=\F^{(0,-1)(N-1)}_{\bfu^{(N-1)''}}\tot\cdots \tot \F^{(0,-1)(0)}_{\bfu^{(0)''}}.
\en
We define the elliptic dynamical  $R$-matrix $R_{T^{1/2}}(\bfu_1,\bfu_2;\bgz)\in\End_\FF(\F^{(0,-1)}_{\bfu_1}\tot\F^{(0,-1)}_{\bfu_2})$ 
as the following transition matrix of the elliptic stable envelopes.  For 
$\bal=(\bal',\bal''),\ \bbe=(\bbe',\bbe''), \bga=(\bga',\bga'') 
\in\Lambda_{\bfv',\bfw'}\times \Lambda_{\bfv'',\bfw''}$, 
\bea
&&R_{T^{1/2}}(\bfu_1,\bfu_2;\bgz)\ket{{\bbe}'}_{\bfu_1}\tot \ket{\bbe''}_{\bfu_2}=\sum_{\bal
}R_{T^{1/2}}(\bfu_1,\bfu_2;\bgz)^{\bbe}_{\bal}\ket{\bal'}_{\bfu_1}\tot \ket{\bal''}_{\bfu_2},\\
&&{R}_{T^{1/2}}(\bfu_1,\bfu_2;\bgz)_\bal^\bbe
=\mu(\bfu_1/\bfu_2)\bar{R}_{T^{1/2}}(\bfu_1,\bfu_2;\bgz)_\bal^\bbe
,\lb{Rmat}\\
&&\widehat{\Stab}_{\overline{\gC},T^{1/2}}(\bar{\bal};\bgz)\bigr|_{\bar{\bga}}=\sum_{\bbe
}\widehat{\Stab}_{\gC,T^{1/2}}({\bbe};\bgz)\bigr|_{{\bga}}\ \bar{R}_{T^{1/2}}(\bfu_1,\bfu_2;\bgz
)_\bal^\bbe.\lb{bRmat}
\ena
Here $\bar{\bal}=(\bal'',\bal')$ etc. and $\mu(\bfu)$ is a scalar function defined by 
\bea
&&\mu(\bfu_1/\bfu_2)\Phi_{\blempty}(\bfu_1)\Phi_{\blempty}(\bfu_2)=\Phi_{\blempty}(\bfu_2)\Phi_{\blempty}(\bfu_1).\lb{comPhiempty2}
\ena
It is explicitly given by
\be
&&\mu(\bfu_1/\bfu_2)=\prod_{k,l\in I}\prod_{i=1}^{w'_k}\prod_{j=1}^{w_l''}\mu({u^{(l)''}_j}/{u^{(k)'}_i})_{kl},\\
&&\hspace{-1cm}\mu(v/u)_{kl}=\left(v/u\right)^{-\eta_{kl}}
\frac{\Gamma(t_2^{-(k-l)}v/u;t_1^N,t^N_2,\hbar)\Gamma(t_1^Nt_1^{k-l}v/u;t_1^N.t_2^N,\hbar)}{\Gamma(t_1^Nt_2^{-(k-l)}v/u;t^N_1,t_2^N,\hbar)\Gamma(t_1^Nt_2^Nt_1^{k-l}v/u;t_1^N,t^N_2,\hbar)}\Biggl/(\hbar\mapsto p)\qquad (k\leq l)\\
&&\hspace{-1cm}\mu(v/u)_{kl}=\frac{1}{\mu(u/v)_{lk}}\qquad (k> l).
\en
Here $\Gamma(z;t_1^N,t_2^N,\hbar)$ etc. denotes the triple Gamma function  defined by
\be
&&\Gamma(z;t_1,t_2,\hbar)=(z;t_1,t_2,\hbar)_\infty(t_1t_2\hbar/z;t_1,t_2,\hbar)_\infty,\qquad \\
&&(z;t_1,t_2,\hbar)_\infty=\prod_{m_1,m_2,m_3=0}^\infty(1-zt_1^{m_1}t_2^{m_2}\hbar^{m_3}).
\en

By definition, we have
\bea
&&[\hbar^{h^{(1)}_j+h^{(2)}_j}, R(\bfu_1, \bfu_2;\bgz)]=0\qquad \forall j\in I.
\ena
Here $h^{(1)}_j=h_j\tot 1$ and $h^{(2)}_j=1\tot h_j$ and
\be
&&\hbar^{h_j}\ket{\bla'}_{\bfu_1}=\hbar^{\wt(\bla')_j}\ket{\bla'}_{\bfu_1}\qquad \mbox{etc.} 
\en
We assume 
\bea
&&{R}_{T^{1/2}}(\bfu_1,\bfu_2;\bgz\hbar^{\bfh^{(1)}+\bfh^{(2)}})={R}_{T^{1/2}}(\bfu_1,\bfu_2;\bgz).
\lb{gzshiftinv}
\ena
Here we set $\bgz\hbar^{\bfh}=\{ \gz_j \hbar^{h_j} \}_{j\in I}$.
Note also that $\widehat{\Stab}_{{\gC},T^{1/2}}({\bullet};\gz)$ is depend on $\bfu_1, \bfu_2$ only through the chamber 
$\gC
$ essentially.  Hence for any $a\in \C^\times$
\bea
    &&{R}_{T^{1/2}}(a\bfu_1,a\bfu_2;\gz)={R}_{T^{1/2}}(\bfu_1,\bfu_2;\gz).\lb{Rauau}
\ena

\begin{prop}\lb{comPhi2}
The type I vertex operators satisfy the following exchange relation.
\be
&&\Phi^{(\blmu'')}_{\bom''}(\bfu_2)\Phi^{(\blmu')}_{\bom'}(\bfu_1)
=\sum_{\bla=(\bla',\bla'')\in  \Lambda_{\bfv',\bfw'}\times\Lambda_{\bfv'',\bfw''}
}
{R}_{T^{1/2}}(\bfu_1,\bfu_2;\bgz)^\bla_{{\bom}}
\ 
\Phi^{(\blmu')}_{\bla'}(\bfu_1)\Phi^{(\blmu'')}_{\bla''}(\bfu_2),
\en
where  ${\bom}=(\bom',\bom'')\in  \Lambda_{\bfv',\bfw'}\times \Lambda_{\bfv'',\bfw''}$. 
\end{prop}
\noindent
{\it Proof.}\ 
Let us consider the RHS. 
By using \eqref{Phi1Phi2}, \eqref{comPhiempty2} and 
\be
&&\gz_j\Phi_{\bla'}(\bfu_1)=\hbar^{-\wt(\bla')_j}\Phi_{\bla'}(\bfu_1)\gz_j \qquad \mbox{etc.},
\en
move $\bar{R}_{{T^{1/2}}}(u_1,u_2;\gz)^{\bla}_{\bom}$ to the right end. One obtains
\be
&&\sum_{\bla=(\bla',\bla'')}\prod_{a\in \bla'}\left(\int_0^{\varphi^{\bmu'}_{\iota(a)}}d_px'_a\right)
\prod_{b\in \bla''}\left(\int_0^{\varphi^{\bmu''}_{\iota(b)}}d_px''_b\right)
\prod_{i=0}^{N-1}\left( :\prod_{a\in \bla\atop c_a\equiv i}x^-_i(x_a): \right)\ \Phi_{\blempty}(\bfu_2)\Phi_{\blempty}(\bfu_1)\\
&&\qquad\qquad \times  \prod_{a, b\in \bla\atop {c_a\equiv c_b \atop \rho_a<\rho_b}}<x^-_{c_a}(x_a)x^-_{c_b}(x_b)>^{sym}
\ \widehat{\Stab}_{\gC,T^{1/2}}(\bla;\bgz)\bar{R}_{T^{1/2}}(\bfu_1,\bfu_2;\bgz\hbar^{-\wt(\bla)}
)^\bla_{{\bom}}.
\en
Note that $\ds{:\prod_{a\in \bla\atop c_a\equiv i}x^-_i(x_a): }$ and  $\ds{\prod_{a, b\in \bla\atop {c_a\equiv c_b \atop \rho_a<\rho_b}}<x^-_{c_a}(x_a)x^-_{c_b}(x_b)>^{sym}}$   are symmetric in $x_a$, $a\in \bla$ for each fixed  $c_a\equiv i\in I$ so that they are independent from a choice $\bla\in \Lambda_{\bfv',\bfw'}\times \Lambda_{\bfv'',\bfw''}$. Evaluating the Jackson integral and  using a similar argument to the proof of Proposition \ref{QPStab} and \eqref{bRmat},  
one obtains
\be
&&\hspace{-1cm}\sum_{{\bf d}\in \N^{|\bullet|}}\ \prod_{i\in I}\gz_i^{-\sum_{a\in \bullet \atop c_a\equiv i}d_a}
\prod_{i\in I}\left(:\prod_{a\in \bullet 
\atop c_a\equiv i}x^-(x_a): \right)
\Biggr|_{x_a=\varphi^\bmu_{\iota(a)}p^{d_a}}\Phi_{\blempty}(\bfu_2)\Phi_{\blempty}(\bfu_1)
\\&&\hspace{-1cm} \quad\times  
\prod_{a, b\in \bullet
\atop  {c_a\equiv c_b \atop \rho_a<\rho_b}}<x^-_{c_a}(x_a)x^-_{c_b}(x_b)>^{sym}\Biggr|_{x_a=\varphi^\bmu_{\iota(a)}p^{d_a}}\hspace{-0.5cm}
\times \sum_{\bla}\widehat{\Stab}_{\gC,\tT^{1/2}}(\bla;\bgz
)\bigr|_{\bmu}
\bar{R}_{T^{1/2}}(\bfu_1,\bfu_2;\bgz
)^\bla_{{\bom}}
\\
&=&\sum_{{\bf d}\in \N^{|\bullet|}}\ 
\prod_{i\in I}\gz_i^{-\sum_{a\in \bullet \atop c_a\equiv i}d_a}
\prod_{i\in I}\left(:\prod_{a\in \bullet 
\atop c_a\equiv i}x^-(x_a): \right)
\Biggr|_{x_a=\varphi^\bmu_{\iota(a)}p^{d_a}}\Phi_{\blempty}(\bfu_2)\Phi_{\blempty}(\bfu_1)
\\&& \qquad\times  
\prod_{a, b\in \bullet
\atop  {c_a\equiv c_b \atop \rho_a<\rho_b}}<x^-_{c_a}(x_a)x^-_{c_b}(x_b)>^{sym}\Biggr|_{x_a=\varphi^\bmu_{\iota(a)}p^{d_a}}
\times \widehat{\Stab}_{\overline{\gC},\tT^{1/2}}(\bar{\bom};\bgz
)\bigr|_{\bar{\bmu}}\\
&=&\Phi^{(\blmu'')}_{\bom''}(\bfu_2)\Phi^{(\blmu')}_{\bom'}(\bfu_1).
\en
Here $\bullet$ denotes any partition in $\Lambda_{\bfv',\bfw'}\times \Lambda_{\bfv'',\bfw''}$. We also set $\{x_a=\varphi^\bmu_{\iota(a)}p^{d_a}\}:=\{x'_a=\varphi^{\bmu'}_{\iota(a)}p^{d_a}\ (a\in \bla')\}\cup \{x''_b=\varphi^{\bmu''}_{\iota(b)}p^{d_b}\ (b\in \bla'')\}$.
\qed

Similarly, letting 
\be
&&\hspace{-0.5cm}\F^{(0,-1)}_{\overline{\bfu_1}}:=\F^{(0,-1)(0)}_{\overline{\bfu^{(0)'}}}\tot \cdots \tot
\F^{(0,-1)(N-1)}_{\overline{\bfu^{(N-1)'}}},\qquad 
\F^{(0,-1)}_{\overline{\bfu_2}}:= \F^{(0,-1)(0)}_{\overline{\bfu^{(0)''}}}\tot\cdots \tot
\F^{(0,-1)(N-1)}_{\overline{\bfu^{(N-1)''}}}, 
\en
we define the elliptic dynamical  $R$-matrix $R^*_{T^{1/2}_{opp}}(\bfu_1,\bfu_2;\bgz^{*-1})\in\End_\FF(\F^{(0,-1)}_{\overline{\bfu_1}}\tot\F^{(0,-1)}_{\overline{\bfu_2}})$ 
by
\bea
&&R^*_{T^{1/2}_{opp}}(\bfu_1,\bfu_2;\bgz^{*-1})\ket{{\bbe}'}_{\bfu_1}\tot \ket{\bbe''}_{\bfu_2}=\sum_{\bal
}R^{*}_{T^{1/2}_{opp}}(\bfu_1,\bfu_2;\bgz^{*-1})^{\bbe}_{\bal}\ket{\bal'}_{\bfu_1}\tot \ket{\bal''}_{\bfu_2},\\
&&{R}^*_{T^{1/2}_{opp}}(\bfu_1,\bfu_2;\bgz^{*-1})_\bal^\bbe
=\mu^*(\bfu_1/\bfu_2)\bar{R}^*_{T^{1/2}_{opp}}(\bfu_1,\bfu_2;\bgz^{*-1})_\bal^\bbe
,\lb{Rmats}\\
&&\widetilde{\Stab}^*_{\overline{\gC},T^{1/2}_{opp}}(\bar{\bal};\bgz^{*-1})\bigr|_{\bar{\bga}}=\sum_{\bbe
}\widetilde{\Stab}^*_{\gC,T^{1/2}_{opp}}({\bbe};\bgz^{*-1})\bigr|_{{\bga}}\ \bar{R}^*_{T^{1/2}_{opp}}(\bfu_1,\bfu_2;\bgz^{*-1}
)_\bal^\bbe.\lb{bRmats}
\ena
Here $\bal=(\bal',\bal''),\ \bbe=(\bbe',\bbe''), \bga=(\bga',\bga'') \in\Lambda_{\bfv',\bfw'}\times \Lambda_{\bfv'',\bfw''}$ We also set $\bar{\bal}=(\bal'',\bal')$ etc. The function $\mu^*_0(\bfu)$ is  defined by 
\bea
&&\Psi^*_{\blempty}(\bfu_1)\Psi^*_{\blempty}(\bfu_2)=\mu^*(\bfu_1/\bfu_2)\Psi^*_{\blempty}(\bfu_2)\Psi^*_{\blempty}(\bfu_1).\lb{comPsisempty2}
\ena
It is explicitly given by
\be
&&\mu^*(\bfu_1/\bfu_2)=\prod_{k,l\in I}\prod_{i=1}^{w'_k}\prod_{j=1}^{w_l''}\mu^*({u^{(k)'}_i}/{u^{(l)''}_j})_{kl},\\
&&\hspace{-1cm}\mu^*(u/v)_{kl}=\left(u/v\right)^{\eta_{kl}}
\frac{\Gamma(\hbar t_2^Nt_1^{-(k-l)}u/v;t_1^N,t^N_2,\hbar)\Gamma(\hbar t_1^Nt_2^Nt_2^{k-l}u/v;t^N_1,t_2^N,\hbar)}{\Gamma(\hbar t_1^{-(k-l)}t_1^N,t_2^N,\hbar)\Gamma(\hbar t_2^Nt_2^{k-l}u/v;t_1^N,t^N_2,\hbar)}\\
&&\qquad\quad  \times \frac{\Gamma(t_2^Nt_1^{-(k-l)}u/v;t_1^N,t^N_2,p^*)\Gamma( t_1^Nt_2^Nt_2^{k-l}u/v;t^N_1,t_2^N,p^*)}{\Gamma(t_1^{-(k-l)}t_1^N,t_2^N,p^*)\Gamma( t_2^Nt_2^{k-l}u/v;t_1^N,t^N_2,p^*)}\qquad (k\leq l),\\
&&\hspace{-1cm}\mu^*(u/v)_{kl}=\frac{1}{\mu^*(v/u)_{lk}}\qquad (k>l).
\en
By definition, we have
\bea
&&[\hbar^{h^{(1)}_j+h^{(2)}_j}, R^*(\bfu_1, \bfu_2;\bgz^{*-1})]=0\qquad \forall j\in I.
\ena
We assume 
\bea
&&{R}^*_{T^{1/2}_{opp}}(\bfu_1,\bfu_2;\bgz^{*-1}\hbar^{\bfh^{(1)}+\bfh^{(2)}})={R}^*_{T^{1/2}_{opp}}(\bfu_1,\bfu_2;\bgz^{*-1}).
\lb{Rsgzshiftinv}
\ena

\begin{prop}\lb{comPsis2}
The type II dual vertex operators satisfy the following exchange relation.
\be
&&\Psi^{*(\blmu')}_{\bom'}(\bfu_1)\Psi^{*(\blmu'')}_{\bom''}(\bfu_2)
=\sum_{\bla=(\bla',\bla'')\in  \Lambda_{\bfv',\bfw'}\times\Lambda_{\bfv'',\bfw''}
}
\ 
\Psi^{*(\blmu'')}_{\bla''}(\bfu_2)\Psi^{*(\blmu')}_{\bla'}(\bfu_1)
{R}^{*}_{T^{1/2}}(\bfu_1,\bfu_2;\bgz^{*})^{{\bom}}_\bla,
\en
where  ${\bom}=(\bom',\bom'')\in  \Lambda_{\bfv',\bfw'}\times \Lambda_{\bfv'',\bfw''}$. 
\end{prop}
Here we use the formula ( Proposition 6.2 in \cite{KS})
\be
&&R^*_{T^{1/2}}(\bfu_1,\bfu_2;\bgz^*)={}^tR^*_{T^{1/2}_{opp}}(\bfu_1,\bfu_2;\bgz^{*-1}).
\en

Finally, the type I and the type II dual vertex operators exchange by a scalar function. 
\begin{prop}\lb{comPhiPsis}
In the level $(1,L)$ representation, one has 
\be
&&\Phi_\bla(\bfu)\Psi^*_\bmu(\bfv)=\chi(\bfu/\bfv)\Psi^*_\bmu(\bfv)\Phi_\bla(\bfu)\qquad \forall \bla\in 
\La_{\bfv',\bfw'}, \ \bmu\in \La_{\bfv'',\bfw''}
,\nn\\
&&\chi(\bfu/\bfv)=\prod_{k,l\in I}\prod_{i=1}^{w_k}\prod_{j=1}^{w_l''}\chi(u^{(k)}_i/v^{(l)}_j)_{kl},\\
&&\chi(u/v)_{kl}=(u/v)^{-\eta_{kl}}
\frac{\Gamma(\hbar^{1/2}t_2^Nt_2^{k-l} u/v;t_1^N,t_2^N,\hbar)\Gamma(\hbar^{1/2}t_1^{-(k-l)}u/v;t_1^N,t_2^N,\hbar)}{\Gamma(\hbar^{1/2}t_2^Nt_1^{-(k-l)} u/v;t_1^N,t_2^N,\hbar)\Gamma(\hbar^{1/2}t_1^Nt_2^Nt_2^{k-l}u/v;t_1^N,t_2^N,\hbar)}\qquad (k\leq l),\\
&&\chi(u/v)_{kl}=(u/v)^{-\eta_{kl}}
\frac{\Gamma(\hbar^{1/2}t_2^{k-l} u/v;t_1^N,t_2^N,\hbar)\Gamma(\hbar^{1/2}t_1^Nt_1^{-(k-l)}u/v;t_1^N,t_2^N,\hbar)}{\Gamma(\hbar^{1/2}t_1^Nt_2^Nt_1^{-(k-l)} u/v;t_1^N,t_2^N,\hbar)\Gamma(\hbar^{1/2}t_1^Nt_2^{k-l}u/v;t_1^N,t_2^N,\hbar)}\qquad (k> l).
\en
\end{prop}

\subsection{Dynamical Yang-Baxter equation}
As an  associativity condition for the composition of the  vertex operators, one obtains the following dynamical Yang-Baxter equation. 
\bea
&&R^{(12)}_{T^{1/2}}(\bfu_1,\bfu_2;\bgz\hbar^{\bfh^{(3)}}) R^{(13)}_{T^{1/2}}(\bfu_1,\bfu_3;\bgz) R^{(23)}_{T^{1/2}}(\bfu_2,\bfu_3;\bgz\hbar^{\bfh^{(1)}})\nn\\
&& \qquad= R^{(23)}_{T^{1/2}}(\bfu_2,\bfu_3;\bgz)R^{(13)}_{T^{1/2}}(\bfu_1,\bfu_3;\bgz\hbar^{\bfh^{(2)}}) 
R^{(12)}_{T^{1/2}}(\bfu_1,\bfu_2;\bgz).\lb{DYBE}
\ena
The dynamical shift in $\bgz=\{\gz_i\}_{i\in I}$ follows from 
\be
&&\Phi_\bla(\bfu)\gz_i=\hbar^{-\langle h_i,\sum_{j\in I}(w_j\La_j-v_j\al_j)\rangle}\gz_i
\Phi_\bla(\bfu),\\
&&\hbar^{h_i}\ket{\bla}_\bfu=\hbar^{\wt(\bla)_i}\ket{\bla}_\bfu
\en
and \eqref{wtla=wtPhi}.

\section{$L$-operator of $U_{t_1,t_2,p}(\mathfrak{gl}_{N,tor})$}\lb{sec:Lop}

For the elliptic dynamical $R$-matrices introduced in the last section, we construct a $L$-operator $L^+$ on $\F^{(1,L)}_{\bullet,\blLa_L,\blLa_M}$ satisfying the $RLL=LLR^*$ relation as a composition of the type I and the type II dual vertex operators for $\cM(\bfv,(\delta_{i,k}))$. 
We then consider its universal form $\cL^+$ and define the standard comultiplication $\Delta$ in terms of it and show that the vertex operators in Sec.\ref{sec:VObasic} 
are intertwining operators of the $U_{t_1,t_2,p}(\mathfrak{gl}_{N,tor})$-modules w.r.t.  $\Delta$.  

\subsection{$L$-operator on $\F^{(1,L)}_{\bullet,\blLa_L,\blLa_M}$}\lb{sec:LoponF1N}
Let  $k\in I$ and $\sigma : \xi \tot \eta \to \eta \tot \xi$ and  
 consider the following composition of the type I and the type II dual vertex operators for $\cM(\bfv,(\delta_{i,k}))$.  
\be
&&\F^{(0,-1)(k)}_u \tot \F^{(1,L)}_{(\hbar^{1/2})_k\cdot\bxi,\blLa_L,\blLa_M} \stackrel{\sigma}{\longrightarrow} \F^{(1,L)}_{(\hbar^{1/2})_k\cdot\bxi,\blLa_L,\blLa_M} \tot \F^{(0,-1)(k)}_u \\
 &&\qquad \stackrel{\Phi^{(k)}(\hbar^{1/2} u)\tot \id}{\longrightarrow} \F^{(0,-1)(k)}_{\hbar^{1/2} u} \tot \F^{(1,L+1)}_{\bxi_{(\delta_{i,k})},\blLa_L+\La_k,\blLa_M} \tot \F^{(0,-1)(k)}_u \stackrel{\id\tot \Psi^{*(k)}(u)}{\longrightarrow} \F^{(0,-1)(k)}_{\hbar^{1/2} u} \tot \F^{(1,L)}_{\bxi,\blLa_L,\blLa_M+\La_k}.
\en
Here $(\hbar^{1/2})_k\cdot\bxi=(\xi_0,\cdots,\hbar^{1/2}\xi_k,\cdots,\xi_{N-1})$ and $\bxi_{(\delta_{i,k})}=(\xi_0,\cdots,-\xi_k/\hbar^{1/2}u,\cdots,\xi_{N-1})$. Hence we have the operator 
\be
&&L^{+(k)}(u):=g_k(\id\tot \Psi^{*(k)}(u))\circ (\Phi^{(k)}(\hbar^{1/2} u)\tot \id)\sigma\\
&&\quad  : \ 
\F^{(0,-1)(k)}_u \tot \F^{(1,L)}_{(\hbar^{1/2})_k\cdot \bxi,\blLa_L,\blLa_M} \to \F^{(0,-1)(k)}_{\hbar^{1/2} u} \tot \F^{(1,M)}_{\bxi,\blLa_L,\blLa_M+\La_k}.  
\en
Here we set
\be
&&g_k=(-\hbar^{1/2})^{\eta_{kk}}(\hbar ;\hbar )_\infty(\hbar t_1^N;t_1^N,\hbar )_\infty(\hbar t_2^N;t_2^N,\hbar )_\infty.
\en

For $\mu, \nu\in \cP^{(k)}$, we define the components of $L^{+(k)}(u)$ by
\be
L^{+(k)}(u)\cdot \ket{\nu}^{(k)}_u\tot \eta &=&\sum_{\mu}\ket{\mu}^{(k)}_{\hbar^{1/2}  u}\tot L^{+(k)}_{\mu\nu}(u)\eta,
\en
for $\ket{\mu}^{(k)}_u\tot \eta \in \F^{(0,-1)(k)}_u \tot \F^{(1,L)}_{(\hbar^{1/2})_k\cdot \bxi,\blLa_L,\blLa_M} $. 
On the other hand one has
\be
L^{+(k)}(u)\cdot \ket{\nu}^{(k)}_u\tot \eta 
&=&g_k(\id\tot \Psi^{*(k)}(u))\circ (\Phi^{(k)}(\hbar^{1/2} u)\tot \id)
\cdot \eta \tot  \ket{\nu}^{(k)}_u\\
&=&g_k(\id\tot \Psi^{*(k)}(u))\cdot \sum_{\mu}\ket{\mu}^{(k)}_{\hbar^{1/2} u}\tot \Phi^{(k)}_\mu(\hbar^{1/2} u)\eta  \tot  \ket{\nu}^{(k)}_u\\
&=&g_k\sum_{\mu}\ket{\mu}^{(k)}_{\hbar^{1/2} u} \tot \Psi^{*(k)}_\nu(u)\Phi^{(k)}_\mu(\hbar^{1/2} u)\eta.
\en
Hence we obtain 
\bea
&&L^{+(k)}_{\mu\nu}(u)=g_k\Psi^{*(k)}_\nu(u)\Phi^{(k)}_\mu(\hbar^{1/2} u).
\ena
In particular,  setting $k^{+(k)}_\emptyset(u)=L^{+(k)}_{\emptyset\emptyset}(u)$, one has
\bea
k^{+(k)}_\emptyset(u)&=&g_k\Psi^{*(k)}_\emptyset(u)\Phi^{(k)}_\emptyset(\hbar^{1/2} u)\nn\\
&=&\hbar^{\frac{1}{2}h_{\La_k}}e^{Q_{\La_k}}:\exp\left\{
\sum_{m\not=0}\frac{1}{[m]_q}\frac{1-\hbar^m}{1-p^{m}}A_{k,m}(p^{*-1}\hbar^{1/2} u)^{-m}\right\}:.
\ena
By using \eqref{ComAA}, one obtains the following commutation relation. 
\begin{prop}
\be
&&k^{+(k)}_\emptyset(u)k^{+(k)}_\emptyset(v)=\rho(u/v)k^{+(k)}_\emptyset(v)k^{+(k)}_\emptyset(u),\\
&&\rho(u)=\frac{\rho^{+*}(u)}{\rho^+(u)},\qquad 
\rho^+(u)=\frac{\Gamma(1/u;t_1^N,t_2^N,p)}{\Gamma(u;t_1^N,t_2^N,p)},
\qquad \rho^{+*}(u)=\rho^+(u)\Bigr|_{p\mapsto p^*}. 
\en
\end{prop}
Note that $\rho^+(u)$  is independent of the color $k$. 

Let us now consider the following elliptic dynamical $R$-matrices. 
\bea
&&R^+_{T^{1/2}}(u,v;\bgz)^{\bbe}_{\bal}:=\rho^+(u/v)\bar{R}_{T^{1/2}}(u,v;\bgz)^{\bbe}_{\bal},
\ena
where $\bar{R}_{T^{1/2}}
$ 
is given in \eqref{bRmat} with $\bfw'=\bfw''=(\delta_{i,k})$ and $R^{+*}_{T^{1/2}}=R^+_{T^{1/2}}|_{p\mapsto p^*}$. 
\begin{prop}
The $L^{+(k)}$ operator satisfies the following relation.
\bea
&&\hspace{-1cm}\sum_{\mu',\nu'}R^{+}_{T^{1/2}}(u,v;\bgz)^{\mu'\nu'}_{\mu\nu}L^{+(k)}_{\mu'\mu''}(u)
L^{+(k)}_{\nu'\nu''}(v)
= \sum_{\mu',\nu'}L^{+(k)}_{\nu\nu'}(v)L^{+(k)}_{\mu\mu'}(u)
R^{+*}_{T^{1/2}}(u,v;\bgz^*)^{\mu''\nu''}_{\mu'\nu'}\lb{RLL}
\ena
\end{prop}

\noindent
{\it Proof.}\ The proof is the same as Proposition 7.2 in \cite{KS} by using \eqref{Rauau}, Propositions \ref{comPhi2}, \ref{comPsis2}, \ref{comPhiPsis} and 
\be
 &&\frac{\rho^+(u)}{\rho^{+*}(u)}=\frac{\chi(\hbar^{1/2} /u)_{kk}}{\chi(\hbar^{1/2} u)_{kk}}
 \frac{\mu(u)_{kk}}{\mu^*(u)_{kk}}. 
\en
\qed

\subsection{Intertwining relations}\lb{sec:intrel}
We next derive the exchange relation between the vertex operators and  $L^+(u)$. 
They turns out to be the intertwining relations w.r.t. the standard comultiplication defined in terms of the universal form of  $L^+(u)$. 
\begin{prop}\lb{Intertwinrel}
The type I and the type II dual vertex operators satisfy the following relations.
\bea
&&\Phi^{(k)}_{\nu}(\hbar^{1/2} v)L^{+(k)}_{\mu\mu''}(u)=\sum_{\mu'\nu'}
R^+_{T^{1/2}}(u,v;\bgz)_{\mu\nu}^{\mu'\nu'}L^{+(k)}_{\mu'\mu''}(u)\Phi^{(k)}_{\nu'}(\hbar^{1/2} v),\lb{intTypeIcomp}\\
&&L^{+(k)}_{\mu\mu''}(u)\Psi^{*(k)}_{\nu''}(v)=\sum_{\mu'\nu'}\Psi^{*(k)}_{\nu'}(v)L^{+(k)}_{\mu\mu'}(u)R^{+*}_{T^{1/2}}(u,v;\bgz^*)^{\mu''\nu''}_{\mu'\nu'}.\lb{intTypeIIdualcomp}
\ena
\end{prop}
The proof is the same as Proposition 7.3 in \cite{KS}. 

Let $\cL^{+(k)}(u)\in \End_\FF(\F^{(0,-1)(k)}_\bullet)\tot\; \cU$ be the universal $L$-operator with 
\be
&&\cL^+(u)\cdot\ket{\nu}^{(k)}_u\tot \eta=\sum_{\mu\in \cP^{(k)}}\ket{\mu}^{(k)}_{\hbar^{1/2}u}\tot \cL^+_{\mu\nu}(u)\eta,\qquad \eta\in \cU
\en
satisfying 
\bea
&&\hspace{-1cm}\sum_{\mu',\nu'}R^{+}_{T^{1/2}}(u,v;\bgz)^{\mu'\nu'}_{\mu\nu}\cL^{+(k)}_{\mu'\mu''}(u)
\cL^{+(k)}_{\nu'\nu''}(v)
= \sum_{\mu',\nu'}\cL^{+(k)}_{\nu\nu'}(v)\cL^{+(k)}_{\mu\mu'}(u)
R^{+*}_{T^{1/2}}(u,v;\bgz^*)^{\mu''\nu''}_{\mu'\nu'}.\lb{univRLL}
\ena
We assume $\cL^+_{\mu\nu}(u)= L^+_{\mu\nu}(u)$ on $\F^{(1,L)}_{\bxi,\blLa_L,\blLa_M}$. 

In order to make a connection to the dynamical Yang-Baxter equation, let us define $\cL^+(u;\bgz^*)$ by  
\bea
&&\cL^+(u)=\cL^+(u;\bgz^*)e^{-c\otimes Q_{\La_0} -\sum_{j=1}^Nh_{\bep_j}\tot Q_{\bep_j}}.
\ena
We expect that $\cL^+(u;\bgz^*)$ is independent from $e^{\pm Q_i}$, $e^{\pm Q_{\La_0}}$  and $e^{\pm Q_{\bep_{j}}}$ \cite{JKOS,KK03}.  One has
\be
\cL^+(u)\cdot \ket{\nu}^{(k)}_u\tot\eta&=&\cL^+(u;\bgz^*)\cdot\ket{\nu}^{(k)}_u\tot e^{ Q_{\La_0}-\sum_{j=1}^N\wt(\nu)^{(k)}_{\bep_j}Q_{\bep_j}}\eta\\
&=&\sum_{\mu\in \cP^{(k)}}\ket{\mu}^{(k)}_{\hbar^{1/2}u}\tot \cL^+_{\mu\nu}(u;\bgz^*)e^{Q_{\La_k}-\sum_{i\in I}v_i(\nu)Q_i}\eta. 
\en
Here we used  $Q_i=Q_{\bep_i}-Q_{\bep_{i+1}}$ and 
\be
&&e^{-c}\ket{\nu}^{(k)}_u=e^1\ket{\nu}^{(k)}_u,\quad \wt(\nu)^{(k)}_{\bep_i}=-\sum_{j=1}^k\delta_{i,j}+v_i(\nu)-v_{i-1}(\nu)
\en
with $v_i(\nu):=\#\{\square \in \nu\ | c_{\square}\equiv i \}$. Note that $\wt(\nu)^{(k)}_i=-\delta_{i,0}+\wt(\nu)^{(k)}_{\bep_i}-\wt(\nu)^{(k)}_{\bep_{i+1}}$.  Hence one has 
\bea
&&\cL^{+(k)}_{\mu\nu}(u)=\cL^{+(k)}_{\mu\nu}(u;\bgz^*)e^{Q_{\La_k}-\sum_{i\in I}v_i(\nu)Q_i}. 
\ena
Substituting this into \eqref{univRLL}, 
 one obtains the following dynamical $RLL$-relation.
 \bea
 &&\sum_{\mu'\nu'}R^{+}_{T^{1/2}}(u,v;\bgz^*\hbar^{\bfh^{(3)}})^{\mu'\nu'}_{\mu\nu}\cL^{+(k)}_{\mu'\mu''}(u;\bgz^*)
\cL^{+(k)}_{\nu'\nu''}(v;\bgz^* \hbar^{\bfh^{(1)}})\nn\\
&&\qquad\quad=\sum_{\mu'\nu'} \cL^{+(k)}_{\nu\nu'}(v;\gz^*)\cL^{+(k)}_{\mu\mu'}(u;\gz^*\hbar^{\bfh^{(2)}})
R^{+*}_{T^{1/2}}(u,v;\bgz^*)^{\mu''\nu''}_{\mu'\nu'}\lb{DRLL}
\ena
with $\bgz^*\hbar^{\bfh^{(3)}}=\{\gz^*_i\hbar^{h^{(3)}_i}\}$ etc.. Here we used
\be
&&e^{Q_{\La_k}-\sum_{i\in I}v_i(\mu'')Q_i}\gz^*_j=\hbar^{\wt(\mu'')^{(k)}_j}\gz^*_j
e^{Q_{\La_k}-\sum_{i\in I}v_iQ_i},\\
&&\wt(\mu'')^{(k)}_j=-\delta_{k,j}+\sum_{i\in I}v_i(\mu'')a_{ij}=-\langle h_j,\La_k-\sum_{i\in I}v_i(\mu'')\al_i\rangle=h_j\bigl|_{\ket{\mu''}^{(k)}_u}\quad \mbox{etc.}.
\en
Comparing this with \eqref{DYBE}, we assume 
\be
&&
\cL^{+(k)}_{\mu\nu}(u;\bgz^*)\cdot\ket{\la}^{(k)}_v=\sum_{\la'}R^{+}_{T^{1/2}}(u,v;\bgz^*)_{\mu\la'}^{\nu\la}\ket{\la'}^{(k)}_v, \quad \ket{\la}^{(k)}_v\in \F^{(0,-1)(k)}_v.
\en
 
Now we define a new comultiplication $\Delta$ by
\be
&&\Delta(\cL^{+(k)}_{\mu\nu}(u))=\sum_{\la} \cL^{+(k)}_{\mu\la}(u\gamma^{(2)})\tot \cL^{+(k)}_{\la\nu}(u),\lb{DeltaL}\\
&&\Delta(g(\bgz,p))=g(\bgz,p)\tot1,\quad \Delta(g(\bgz^*,p^*))=1\tot g(\bgz^*,p^*),\quad \forall g(\bgz,p), g(\bgz^*,p^*)\in \FF,\\
&&\Delta(\gamma^{1/2})=\gamma^{1/2}\tot\gamma^{1/2},\quad \Delta(K)=K\tot K.
\en
We call this the standard comultiplication.  Defining further a counit $\vep$ and an antipode $S$ 
in terms of $\cL^{+(k)}(u)$ in the same way as in Sec.4.2 of \cite{Konno16}, one can show that $(U_{t_1,t_2,p}(\mathfrak{gl}_{N,tor}), \Delta, \vep,\mu_l,\mu_r,S)$ is a Hopf algebroid.  One then show the following statement in the same way as in 
Theorem 4.2 in \cite{Konno08}. 

\begin{prop}
Relations in Proposition \ref{Intertwinrel} are the intertwining relations w.r.t. $\Delta$ for 
$\Phi^{(k)}(u)$ and $\Psi^{*(k)}(u)$, respectively i.e. 
\bea
&&\hspace{-0.5cm}\Phi^{(k)}(v)\cL^{+(k)}_{\mu\nu}(u)=\Delta(\cL^{+(k)}_{\mu\nu}(u))\Phi^{(k)}(v),\lb{intTypeI}\\
&&\hspace{-0.5cm}\cL^{+(k)}_{\mu\nu}(u)\Psi^{*(k)}(v)(\eta\tot\ket{\la}^{(k)}_v)=\Psi^{*(k)}(v)\Delta(\cL^{+(k)}_{\mu\nu}(u))(\eta\tot\ket{\la}^{(k)}_v),\quad 
 \lb{intTypeIIdual}
\ena
for $\eta\in \F^{(1,L)}_{\bullet,\blLa_L,\blLa_M},\ \ket{\la}^{(k)}_v\in \F^{(0,-1)(k)}_v$. 
\end{prop}

\noindent
{\it Proof.}\ From \eqref{intTypeI}, one has
\be
\mbox{LHS}\cdot \eta&=&\sum_{\la}\ket{\la}^{(k)}_v\tot \Phi^{(k)}_\la(v)L^{+(k)}_{\mu\nu}(u)\eta,
\\
\mbox{RHS}\cdot \eta&=&\sum_{\mu'}\cL^{+(k)}_{\mu\mu'}(\gamma^{(2)} u)\tot \cL^{+(k)}_{\mu'\nu}(u)\cdot\sum_{\la'}\ket{\la'}^{(k)}_v\tot \Phi^{(k)}_{\la'}(v)\eta\nn\\
\qquad&=&\sum_{\mu',\la',\la}R^+_{T^{1/2}}(\hbar^{1/2} u,v;\bgz^*)^{\mu'\la'}_{\mu\la}\ket{\la}^{(k)}_v\tot L^{+(k)}_{\mu'\nu}(u) \Phi^{(k)}_{\la'}(v)\eta\nn\\
\qquad&=&\sum_\la\ket{\la}^{(k)}_v\tot \sum_{\mu',\la'}R^+_{T^{1/2}}(\hbar^{1/2} u,v;\bgz)^{\mu'\la'}_{\mu\la}L^{+(k)}_{\mu'\nu}(u) \Phi^{(k)}_{\la'}(v)\eta.
\en
In the last equality we used \eqref{reltot}. 
Similarly, from \eqref{intTypeIIdual}, one obtains
\be
\mbox{LHS}&=&L^{+(k)}_{\mu\nu}(u)\Psi^{*(k)}_\la(v)\eta,
\\
\mbox{RHS}&=&\Psi^{*(k)}(v)\cdot\sum_{\mu'}\cL^{+(k)}_{\mu\mu'}(\gamma^{(2)} u)\tot \cL^{+(k)}_{\mu'\nu}(u)\cdot (\eta\tot\ket{\la}^{(k)}_v)\nn\\
\qquad&=&\Psi^{*(k)}(v)\cdot\sum_{\mu',\la'}L^{+(k)}_{\mu\mu'}(u)\eta\tot R^+_{T^{1/2}}(u,v;\bgz^*)^{\nu\la}_{\mu'\la'}
\ket{\la'}^{(k)}_v
\nn
\\
\qquad&=&\Psi^{*(k)}(v)\cdot\sum_{\mu',\la'}R^{+*}_{T^{1/2}}(u,v;\bgz^*\hbar^{-\bfh^{(3)}})^{\nu\la}_{\mu'\la'}L^{+(k)}_{\mu\mu'}(u)\eta\tot \ket{\la'}^{(k)}_v
\nn\\
\qquad&=&\Psi^{*(k)}(v)\cdot\sum_{\mu',\la'}L^{+(k)}_{\mu\mu'}(u)R^{+*}_{T^{1/2}}(u,v;\bgz^*{\hbar^{-\bfh^{(1)}-\bfh^{(3)}}})^{\nu\la}_{\mu'\la'}\eta\tot \ket{\la'}^{(k)}_v\nn\\
\qquad&=&\sum_{\mu',\la'}\Psi^{*(k)}_{\la'}(v)L^{+(k)}_{\mu\mu'}(u)
R^{+*}_{T^{1/2}}(u,v;\bgz^*
)^{\nu\la}_{\mu'\la'}\eta.
\en
To obtain the fourth line, we used $\gz_i^*=\gz_i\hbar^{-h_i}$ 
and \eqref{reltot}. Note also 
$p^*=p\hbar^{-1}$ on  
$\F^{(1,L)}_{\bullet,\blLa_L,\blLa_M}$. 
\qed

\section*{Acknowledgements}
The authors would like to thank Tatsuyuki Hikita and Taro Kimura for useful discussions. 
A part of the results was obtained during H.K. stays at the Mittag-Leffler Institute, 19 March - 25 April, 2025.  He would like to thank the organizers of the program ``Cohomological Aspects of Quantum Field Theory'' and the staffs of the institute for their hospitality.  H. K. is supported by the Grant-in-Aid for Scientific Research (C) 23K03029 JSPS, Japan. Work of A. S. is supported in part by the NSF under grant DMS - 2054527.

\bigskip

\renewcommand{\baselinestretch}{0.7}

\bigskip
\begin{flushleft}

{\it Department of Mathematics, Tokyo University of Marine Science and 
Technology, \\Etchujima, Koto, Tokyo 135-8533, Japan\\[-3mm]
E-mail address:        hkonno0@kaiyodai.ac.jp}\\
       
 \bigskip
{\it Department of Mathematics, University of North Carolina at Chapel Hill,  \\
Chapel Hill, NC 27599-3250, USA\\[-3mm]
E-mail address:             asmirnov@unc.edu}

\end{flushleft}

\end{document}